\documentclass[smallcondensed]{svjour3}       
\smartqed  
\usepackage{geometry}
\usepackage{xcolor}
\usepackage{multicol}
\usepackage{graphicx}
\usepackage{subfigure}
\usepackage{verbatim}
\usepackage{amssymb,amsmath,sansmath,lmodern}
\usepackage{color}

\usepackage{accents}
\usepackage{stmaryrd}
\usepackage{tabu}
\usepackage{booktabs}

\usepackage[numbers]{natbib}
\usepackage{textcomp}

\newcommand{\aver}[1]{\ensuremath{\left\{\!\left\{#1\right\}\!\right\}}}
\newcommand{\jump}[1]{\ensuremath{\left\llbracket #1 \right\rrbracket}}

\newcommand*\diff{\mathop{}\!\mathrm{d}}

\DeclareMathAccent{\svec}{\mathord}{letters}{126}

\newcommand\stvec[1]{\mathbf #1}				
\newcommand\stvecg[1]{\boldsymbol #1}				
\newcommand\ssvec[1]{\svec{\stvec{#1}}}	
\newcommand\ssvecg[1]{\svec{\stvecg{#1}}}	
\newcommand\cssvec[1]{\svec{\tilde{\stvec{#1}}}} 
\newcommand\cssvecg[1]{\svec{\tilde{\stvecg{ #1}}}}

\usepackage{xargs}    
\usepackage[colorinlistoftodos,prependcaption,textsize=tiny]{todonotes}
\newcommandx{\unsure}[2][1=]{\todo[linecolor=blue,backgroundcolor=blue!25,bordercolor=blue,#1]{#2}}
\newcommandx{\change}[2][1=]{\todo[linecolor=red,backgroundcolor=red!25,bordercolor=red,#1]{#2}}
\newcommandx{\info}[2][1=]{\todo[linecolor=OliveGreen,backgroundcolor=OliveGreen!25,bordercolor=OliveGreen,#1]{#2}}
\newcommandx{\improvement}[2][1=]{\todo[linecolor=Plum,backgroundcolor=Plum!25,bordercolor=Plum,#1]{#2}}
\newcommandx{\thiswillnotshow}[2][1=]{\todo[disable,#1]{#2}}

\begin{document}
	
\title{High--order discontinuous Galerkin approximation for a three--phase incompressible Navier--Stokes/Cahn--Hilliard model}

\titlerunning{High--order DG approximation for a three--phase iNS/CH model}        

\author{Juan Manzanero  \and Carlos Redondo \and Miguel Chávez--Módena \and Gonzalo Rubio \and Eusebio Valero \and Susana Gómez--Álvarez \and Ángel Rivero--Jiménez
}


\institute{Juan Manzanero (\email{juan.manzanero@upm.es}), Carlos Redondo, Miguel Chávez--Módena, Gonzalo Rubio, Eusebio Valero \at
	ETSIAE-UPM - School of Aeronautics, Universidad Polit\'ecnica de Madrid. Plaza Cardenal Cisneros 3, E-28040 Madrid, Spain. //
	Center for Computational Simulation, Universidad Polit\'ecnica de Madrid, Campus de Montegancedo, Boadilla del Monte, 28660, Madrid, Spain. \\
Susana Gómez--Álvarez, Ángel Rivero--Jiménez \at
    Repsol Technology Lab Agust\'in de Betancourt S/N, 28935, M\'ostoles, Madrid, Spain
}

\date{Received: date / Accepted: date}

\maketitle

\begin{abstract}

In this work we introduce the development of a three--phase incompressible Navier--Stokes/Cahn--Hilliard numerical method to simulate three--phase flows, present in many industrial operations. The numerical method is then applied to successfully solve oil transport problems, such as those found in the oil and gas industry. 
The three--phase model adopted in this work is a Cahn--Hilliard diffuse interface model, which was derived by Boyer and Lapuerta \cite{boyer2006study}. The Cahn--Hilliard model is coupled to the entropy--stable incompressible Navier--Stokes equations model derived by Manzanero et al. \cite{2019:Manzanero-iNS}.
The spatial discretization uses a high--order discontinuous Galerkin spectral element method which yields highly accurate results in arbitrary geometries, while an implicit--explicit (IMEX) method is adopted as temporal discretization.
The developed numerical tool is tested for two and three dimensional problems, including a convergence study, a two--dimensional jet, a three--dimensional annular flow, and realistic geometries like T--shaped pipe intersections. 
\end{abstract}

\section{Introduction}\label{sec-xpipe:Introduction}

The transportation of hydrocarbons from the reservoir to the processing facilities is characterized by the modification on flowing pressure and temperature conditions. These changes in operational conditions lead to a transition from typical one--phase behavior to a more complex multiphase flow with different number of phases present along the production system (wells, flowlines, export lines,...).
The vast majority of the reservoirs present fluids which evolve into a mixture of liquid crude oil, natural gas and water. In case flowing bottom whole pressures at the near well location are below certain values, even sand particles or fines can be dragged and produced with the stream, leading to a more complex flow.

The physical phenomena associated with hydrocarbon multiphase flow transport (e.g., change in the flow pattern or phase change) will impact the production process or even lead to safety issues (e.g., liquid overflooding in process facilities due to an underestimation of liquid surges caused by slug flow). 
Therefore, an accurate prediction of the flow distribution and behaviour is mandatory to ensure reliable and continuous transport of the production fluid.
In the last years, the numerical simulation of multiphase flows has become more frequent as it permits to conduct numerical experiments for many industries, such as energy, automotive or aerospace. This is due both to the improvements in multiphase flow models and the increase in computational power with HPC facilities (which, at the same time, limits its application in daily or routine engineering analysis). 
The petroleum industry can take advantage from high fidelity multiphase flow simulation tools to minimize the cost of production system design  as well as to support the optimization of its operation.\\ 


In the oil and gas industry, multiphase flows are usually modelled with one--dimensional (1D) simulations tools (e.g. OLGA\textsuperscript{\textregistered}, Pipesim\textsuperscript{\textregistered} or LedaFlow\textsuperscript{\textregistered}\footnote{OLGA\textsuperscript{\textregistered}/ Pipesim\textsuperscript{\textregistered} are registered trademarks of Schlumberger Inc. and LedaFlow\textsuperscript{\textregistered} is a registered trademark of Kongsberg A/S.}). These models rest on a large number of experimental databases, which results in a high accuracy of the predicted results with a low computational cost \cite{belt2011comparison}. However these 1D modeling tools are limited as they cannot capture some physical details, specially where three--dimensional effects are important. This means that resolving fast transient phenomena (e.g., slug flow) still presents limitations in current 1D modeling tools \cite{gharaibah2015overview}. 
A different approach is followed by multiphase flow simulation based on Computational Fluid Dynamics (CFD), that permits detailed three-dimensional (3D) simulations of immiscible fluids including effects of pressure, temperature and liquid-gas heat and mass transfer. The main objectives of these simulations is to provide detailed qualitative and quantitative evaluation of flow assurance issues such as erosion or slugging, supporting designers and operators to solve flow problems or to extend the life of the flow lines.


Interface capturing methods are among the simplest approaches to treat multiphase flows in CFD. In these models, the governing equations are the continuity and momentum equations for a divergence-free velocity field, in conjunction with a convective equation that tracks the interface. 
Amongst interface capturing methods, such as Volume of Fluid (VOF) \cite{hirt1981volume} or Level Set \cite{sussman1994level}; Diffuse Interface (DI) methods (also known as phase field methods) \cite{jacqmin1999calculation,anderson1998diffuse,badalassi2003computation} provide a useful alternative that does not seem to suffer from problems with either mass conservation or the accurate computation of surface tension. 
Although there are examples of three phase flows (or in general, N-phase flows, with N greater or equal than 3) with Level Set or VOF methods \cite{bonhomme2012inertial,inoue2004mesoscopic,merriman1994motion,saye2011voronoi,smith2002projection,villa2010implicit,zhao1996variational,zheng2007visual,zlotnik2009hierarchical}, most of the work in three phase flows is based on DI methods \cite{boyer2006study,boyer2010cahn,boyer2011numerical,kim2007phase,kim2009generalized,kim2012phase,kim2004conservative,lee2012practically,dong2014efficient,dong2017wall,yang2018multiphase,dong2018multiphase}.
In this paper we focus in the phase field approach. 

In DI methods, a phase--field function that describes the N--phase system is defined. The sharp fluid interface is replaced by a smooth transition layer that connects the two immiscible fluids. 
The free--energy, which represents the effect of the surface tension between the different fluids, is used to characterize the system. The free--energy presents two terms whose effect tend to mix the fluids and separate the fluids respectively \cite{liu2003phase,lowengrub1998quasi}.
The evolution of the phase--field function in our work is modelled by means of the convective Cahn--Hilliard (CH) equation \cite{1958:Cahn}. The use of the CH equation for the evolution of the phase--field function permits an accurate computation of the surface tension and the simulation of phase separation processes. In particular, in this work we use the model of Boyer et al. \cite{boyer2006study} to describe the three--phase system coupled to the incompressible Navier--Stokes (iNS) equations with variable density and artificial (or \textit{pseudo}) compressibility \cite{shen1996new}. A review of alternative iNS/CH models can be found in \cite{hosseini2017isogeometric}.
\\

The three--phase model is numerically approximated in space with a high--order Discontinuous Galerkin Spectral Element Method (DGSEM) \cite{2009:Kopriva} 
that uses the Symmetric Interior Penalty (SIP) method \cite{1978:Wheeler,Ferrer2010,Ferrer2012,2017:Ferrer,manzanero2018bassi}. 
The DGSEM has been used in the past to discretize  multiphase (two phase) flows \cite{fraysse2016upwind,redondo2017artificial,gomez2019novel,2019:Manzanero-MU,manzanero2018high}, and it is popular for its arbitrary order of accuracy \cite{2007:Hesthaven,2009:Kopriva}, low dissipative and dispersive errors \cite{gassner2011comparison,moura2015linear,manzanero2018dispersion,manzanero2020design}, the representation of arbitrary three--dimensional complex geometries through the use of unstructured meshes with curvilinear elements \cite{2006:Kopriva}, efficient mesh adaptation techniques \cite{Kompenhans201636,Kompenhans2016216,2019:Rueda} and the design of provably stable schemes 
\cite{2016:Gassner,2016:Winters,manzanero2018insights,2018:GassnerBR1,2019:Manzanero-CH,2019:Manzanero-iNS,2019:Manzanero-MU}.
Previously, three component Cahn--Hilliard models have been discretized by means of the finite element method \cite{boyer2011numerical}, local discontinuous Galerkin method \cite{xia2007local} or spectral element method \cite{2014:Dong}. 
The DGSEM has been used in the past to discretize the two component Cahn--Hilliard equation \cite{2019:Manzanero-CH} and the three component Cahn--Hilliard equation \cite{2020:Manzanero-UR-CaF}. To the authors' knowledge, this is the first implementation of the three component Cahn--Hilliard model \cite{2006:Boyer} coupled with the Navier--Stokes equations in a discontinuous Galerkin framework. Even though the DGSEM provides us a framework to construct stable schemes, in this work we have not included a stability analysis and it is left for future work. Nevertheless, our results suggest a stable formulation that provides a robust solver.

Finally, for the discretization of time we consider a first order IMplicit--EXplicit (IMEX) time integrator. The linear fourth order spatial operator of the Cahn--Hilliard equation is solved implicitly while the non--linear second order spatial operator is treated explicitly. The solution of the fully--discrete system involves the solution of one linear system for each of the Cahn--Hilliard equations (two for the three phase system). As detailed in \cite{2020:Manzanero-UR-CaF}, the two linear systems are decoupled such that the Jacobian matrices are constant in time and identical for both Cahn--Hilliard equations. Therefore this method permits a resolution in which only one LU factorization is performed for the two equations.\\ 

The rest of this work is organized as follows: we write the governing equations in Sec.~\ref{sec-xpipe:GovEqns}, and we construct its DG approximation and the IMEX time discretization in Sec.~\ref{sec-xpipe:DGSEM}. Finally, we present numerical experiments for the two--phase version of the model in Sec.~\ref{sec-xpipe:Validation}, and for the three--phase model in Sec.~\ref{sec-xpipe:Validation3PH}.

\section{Governing equations}\label{sec-xpipe:GovEqns}

In this work, we couple the three--phase Cahn--Hilliard model of Boyer et al. \cite{boyer2006study,2020:Manzanero-UR-CaF} 
to the incompressible Navier--Stokes with artificial compressibility \cite{2019:Manzanero-iNS}. We define the concentration of Phase~$j$ 
as the relative volume occupied by that phase. Thus, for three--phase flows we have that
\begin{equation}
  c_{1} + c_{2} + c_{3} = 1.
  \label{eq-xpipe:continuous:concentration-sum-1}
\end{equation}
Henceforth, without loss of generality, we consider that the concentrations of 
Phases~1 and 2 can freely vary, and we compute the concentration of Phase~3 from 
\eqref{eq-xpipe:continuous:concentration-sum-1}. For Phases~1 and 2, the concentration 
is computed from the Cahn--Hilliard equation,
\begin{equation}
  c_{i,t}+ \svec{\nabla}\cdot\left(c_{i}\svec{u}\right) = 
  \frac{M_0}{\Sigma_{i}}\svec{\nabla}\cdot\left(\svec{\nabla}\mu_{i}\right),~~i=1,2,
  \label{eq-xpipe:continuous:CHE}
\end{equation}
where $\svec{u}=\left(u,v,w\right)$ is the velocity field, $M_0$, is the mobility and $\mu_{i}$ is the 
chemical potential of Phase~$i$,
\begin{equation}
\mu_{i} = \frac{12}{\varepsilon}\Sigma_i f_{i}- 
\frac{3}{4}\varepsilon\Sigma_{i}\svec{\nabla}^2 c_i,\quad i=1,2,3,
\label{eq-3ph:continuous:chemical-potential}
\end{equation}
with
\begin{equation}
f_i =  \frac{\Sigma_{T}}{3\Sigma_i} \sum_{\substack{j=1 \\ j\neq i}}^{3}
\left(\frac{1}{\Sigma_{j}}\left[\frac{\partial F_0^{\stvecg{\sigma}}}{\partial c_{i}}-
\frac{\partial F_0^{\stvecg{\sigma}}}{\partial c_{j}}\right]\right),\quad\frac{3}{\Sigma_{T}} = 
\frac{1}{\Sigma_{1}}+\frac{1}{\Sigma_{2}}+\frac{1}{\Sigma_{3}}.
\end{equation}
The chemical potentials are algebraically constrained \cite{boyer2006study},
\begin{equation}
  \frac{\mu_{1}}{\Sigma_{1}}+  \frac{\mu_{2}}{\Sigma_{2}}+  
  \frac{\mu_{3}}{\Sigma_{3}}=0,
  \label{mu_constrained}
\end{equation}
where $\Sigma_{i}$ and $\varepsilon$ are positive constants called \textit{spreading factors} and interface width, respectively. The \textit{spreading 
factors} are computed from the interfacial tension between the two phases,
\begin{equation}
  \Sigma_{i} = \sigma_{ij} + \sigma_{ik} - 
  \sigma_{jk},\quad(i,j,k)~\mathrm{cyclical}.
  \label{sigma_relation}
\end{equation}
Finally, as in \cite{boyer2006study}, the \textit{chemical free--energy} $F_0^{ \boldsymbol{\sigma} }$
is a polynomial function on the concentrations,
\begin{equation}
  F_0^{\stvecg{\sigma}} =  \sigma_{12}c_1^2c_2^2+\sigma_{13}c_1^2c_3^2+\sigma_{23}c_2^2c_3^2+c_1c_2c_3\left(\Sigma_1 c_1 +\Sigma_2c_2 + 
  \Sigma_3c_3\right).
  \label{eq-3ph:continuous:F0}
\end{equation}

The density (and all the thermodynamic variables) is computed from the concentration of the three phases,
\begin{equation}
  \rho\left(c_1,c_2,c_3\right) = \rho_1 c_1 + \rho_2 c_2 + \rho_3 c_3 = \rho_1 c_1 + \rho_2 c_2 + \rho_3\left(1-c_1-c_2\right),
  \label{eq-3ph:governing:density}
\end{equation}
where $\rho_{1,2,3}$  are the 
densities of fluids 1, 2 and 3, respectively, assumed constant in space and time. 
The velocity field is given by the momentum equation,
\begin{equation}
\left(\rho \svec{u}\right)_{t} + \svec{\nabla}\cdot\left(\rho \svec{u}\svec{u}\right) 
= - \svec{\nabla}p + \sum_{m=1}^{3}\mu_{m} \svec{\nabla} c_{m} + 
\svec{\nabla}\cdot\left(\eta \left(\svec{\nabla}\svec{u} + 
\svec{\nabla}\svec{u}^{T}\right)\right)+ \rho \svec{g},
\label{eq-xpipe:governing:iNS-momentum}
\end{equation}
where $\eta$ is the viscosity, computed from the (constant) equilibrium phases viscosities $\eta_{1,2,3}$ 
in a similar fashion to the density (see \eqref{eq-3ph:governing:density}). The sum of 
the $\mu_{m}\nabla c_{m}$ products is the phase field approximation of the capillary pressure, and $\svec{g}$ 
is the gravity acceleration.

The pressure is computed with an artificial compressibility model \cite{1997:Shen,2010:Shen-energy},
\begin{equation}
  p_{t}+\rho_0c_0^2\svec{\nabla}\cdot\svec{u}=0,
  \label{eq-xpipe:governing:ACM}
\end{equation}
with $\rho_0=\max\left(\rho_1,\rho_2,\rho_3\right)$ and $c_0$ the artificial compressibility sound speed.

The governing equations \eqref{eq-xpipe:continuous:CHE}, \eqref{eq-xpipe:governing:iNS-momentum} 
and \eqref{eq-xpipe:governing:ACM} are written as a general advection--diffusion 
equation,
\begin{equation}
\stvec{q}_{t}+{\nabla}\cdot\ssvec{f}_{e}\left(\stvec{q}\right)={\nabla}\cdot\ssvec{f}_{v}\left(\stvec{q},{\nabla}\stvec{w}\right)+\stvec{s}\left(\stvec{q},\nabla\stvec{w}\right),
\label{eq-xpipe:continuous:non-conservative-pde}
\end{equation}
with the state vector 
$\stvec{q}=\left(c_1,c_2,\rho\svec{u},p\right)$, gradient variables 
$\stvec{w}=\left(\mu_{1}/\Sigma_1,\mu_{2}/\Sigma_2,\rho\svec{u},p\right)$,
inviscid and viscous fluxes,
\begin{equation}
  \stvec{f}_{e}=\left(\begin{array}{ccc} c_1 u & c_1 v & c_1 w \\ c_2 u & c_2 v & c_2 w \\ \rho u^2+p & \rho uv & \rho uw \\ \rho uv & \rho v^2+p & \rho vw \\ 
  \rho uw & \rho vw & \rho w^2 + p \\ \rho_0c_0^2 u & \rho_0c_0^2 v & \rho_0c_0^2 
  w
\end{array}\right),\quad
  \stvec{f}_{v} = \left(\begin{array}{ccc}\frac{M_0}{\Sigma_1}\mu_{1,x} &\frac{M_0}{\Sigma_1}\mu_{1,y} &\frac{M_0}{\Sigma_1}\mu_{1,z} \\
  \frac{M_0}{\Sigma_2}\mu_{2,x} &   \frac{M_0}{\Sigma_2}\mu_{2,y} &   \frac{M_0}{\Sigma_2}\mu_{2,z}\\
    2\eta \tens{s}_{11} &     2\eta \tens{s}_{12} &     2\eta \tens{s}_{13}\\
  2\eta \tens{s}_{21} &   2\eta \tens{s}_{22} &   2\eta \tens{s}_{23}\\
  2\eta \tens{s}_{31} &   2\eta \tens{s}_{32} &   2\eta \tens{s}_{33}\\ 0 & 0 & 
  0
   \end{array}\right)
   \label{eq:governing:fluxes}
\end{equation}
and source term,
\begin{equation}
  \stvec{s}\left(\stvec{q},\nabla\stvec{w}\right) = \left(\begin{array}{c}
  0 \\
  \rho g_{1} + \mu_{1}c_{1,x}+ \mu_{2}c_{2,x}+ \mu_{3}c_{3,x}\\
    \rho g_{2} + \mu_{1}c_{1,y}+ \mu_{2}c_{2,y}+ \mu_{3}c_{3,y}\\
      \rho g_{3} + \mu_{1}c_{1,z}+ \mu_{2}c_{2,z}+ \mu_{3}c_{3,z}\\
      0
  \end{array}\right).
\end{equation}
In equation
\eqref{eq:governing:fluxes}, $\tens{s}=\frac{1}{2}\left(\nabla\svec{u}+\nabla\svec{u}^{T}\right)$ is the strain tensor. 


\subsection{Reduction of the model to a two--phase flow model}\label{sec-xpipe:3ph-as-2ph}

Constructed this way, the chemical potential satisfies an important consistency
property: when one phase is not initially present (e.g. Phase~2), the chemical potential associated to 
that phase is zero \cite{boyer2006study},
\begin{equation}
  \mu_{2}\bigr|_{c_{2}=0} = 0,
  \label{eq:governing:consistency-mu-zero}
\end{equation}
and the chemical potential associated to the other two phases, which satisfy 
$c_{1}+c_{3}=1$, reduces to that of a two--phase model \cite{boyer2006study,2020:Manzanero-UR-CaF}. We compute the chemical free--energy derivatives for
$c_2=0$, $c_1=c$, and $c_3=1-c$,
  \begin{equation}
\begin{split}
\frac{\partial F_0^{\boldsymbol{\sigma}}}{\partial c_1} &= 2\sigma_{13}c_1c_3^2=\left(\Sigma_{1}+\Sigma_{3}\right)c_1c_3^2=\left(\Sigma_{1}+\Sigma_{3}\right)c\left(1-c\right)^2 ,\\
\frac{\partial F_0^{\boldsymbol{\sigma}}}{\partial c_2} &= c_1 c_3\left(\Sigma_1 c_1 + \Sigma_3 c_3\right)=c\left(1-c\right)\left(\Sigma_1 c + \Sigma_3\left(1-c\right)\right), \\
\frac{\partial F_0^{\boldsymbol{\sigma}}}{\partial c_3} &= 2\sigma_{13}c_1^2 c_3 = \left(\Sigma_{1}+\Sigma_{3}\right) c_1^2c_3=\left(\Sigma_{1}+\Sigma_{3}\right) c^2\left(1-c\right),
\end{split}
\end{equation}
which are replaced into the chemical potential (see \eqref{eq-3ph:continuous:chemical-potential}), defining ${c=c_j}$ and ${c_k=1-c}$,
\begin{equation}
\begin{split}
  \mu_{1}\biggr|_{\substack{c_1=c\\c_2=0\\c_3=1-c}} =& \frac{12\Sigma_1}{\varepsilon}\left(c(1-c)^2 - c^2(1-c)\right)-\frac{3}{4}\Sigma_1\varepsilon\nabla^2 
c,\\
=&\frac{\Sigma_{1}}{2\sigma_{13}}\left(\frac{12\sigma_{13}}{\varepsilon}\frac{\mathrm{d}}{\mathrm{d}c}\left(c^2\left(1-c^2\right)\right)-\frac{3}{2}\sigma_{13}\varepsilon\nabla^2 
c\right)=\frac{\Sigma_{1}}{2\sigma_{13}}\mu^{2\mathrm{ph}}.
\end{split}
\label{eq:governing:consistency-mu-2ph}
\end{equation}

The three--phase flow model also reduces to a two--phase model if one of the 
concentrations is not present initially. If Phase~2 is not initially present,
$c_2(\svec{x},0)=0$, the second Cahn--Hilliard equation is
\begin{equation}
  c_{2,t} =- \svec{\nabla}\cdot\left(c_2\svec{u}\right) +
  \frac{M_0}{\Sigma_{2}}\svec{\nabla}^2\mu_{2}=0,
\end{equation}
since $\mu_2=0$ when $c_{2}=0$, as described in equation  
\eqref{eq:governing:consistency-mu-zero}. Therefore, $c_2\left(\svec{x},t\right)=0$ 
for $t>0$, and Phase~2 will not show in later times. Whereas for Phase~1, its 
Cahn--Hilliard equation is simplified to
\begin{equation}
  c_{t}+\nabla\cdot\left(c\svec{u}\right) = M_0 
  \frac{\Sigma_{1}}{2\sigma_{13}}\nabla^2 \mu^{2\mathrm{ph}},
\end{equation}
which corresponds to a two--phase Cahn--Hilliard model with an adjusted mobility, $M_0\frac{\Sigma_{1}}{2\sigma_{13}}$.

Lastly, we check that the capillary pressure term of the three--phase flow also reduces to that 
of a two--phase flow. For three phases,
\begin{equation}
\begin{split}
  \svec{p}_{c} &= \mu_1\svec{\nabla}c_1 + \mu_2\svec{\nabla}c_{2}+\mu_3\svec{\nabla}c_{3} 
=  \mu_1\svec{\nabla}c_1 +\frac{\Sigma_3}{\Sigma_{1}}\mu_1\svec{\nabla}c_{1}   \\
&=\frac{\Sigma_{1}+\Sigma_{3}}{\Sigma_{1}}\mu_1\svec{\nabla}c_1 = 
\frac{2\sigma_{13}}{\Sigma_{1}}\mu_1\svec{\nabla}c_1 = 
\mu^{2ph}\svec{\nabla}c.
  \end{split}
\end{equation}
Therefore, we confirm that the approximation for the capillary pressure is 
identical in both models. In the second line, we used the property ${\Sigma_i+\Sigma_j=2\sigma_{ij}}$ 
of the spreading factors, ${\nabla c_3=-\nabla c_1}$, and ${\mu_3/\Sigma_3=-\mu_1/\Sigma_1}$ (see  \eqref{eq-xpipe:continuous:concentration-sum-1} \eqref{mu_constrained}  and \eqref{sigma_relation}).

We conclude that the three--phase model is a valid two--phase model
if one of the three phases is not present initially.

\section{Discontinuous Galerkin method and IMEX scheme}\label{sec-xpipe:DGSEM}

The evolution of the three--phase flow is driven by the two Cahn--Hilliard equations (see
\eqref{eq-xpipe:continuous:CHE}), the momentum equation \eqref{eq-xpipe:governing:iNS-momentum}, and the artificial compressibility equation \eqref{eq-xpipe:governing:ACM}. The 
discretization of the spatial differential operators is performed using a nodal 
Discontinuous Galerkin Spectral Element Method (DGSEM), and the 
discretization of the time derivatives is performed using an IMplicit--EXplicit (IMEX)
method. The latter combines a third order low--storage explicit Runge--Kutta RK3 
method, and a first order backward and forward Euler method.

\subsection{Spatial discretization using the DGSEM}

\begin{figure}[h]
  \centering
  \includegraphics[width=0.7\textwidth]{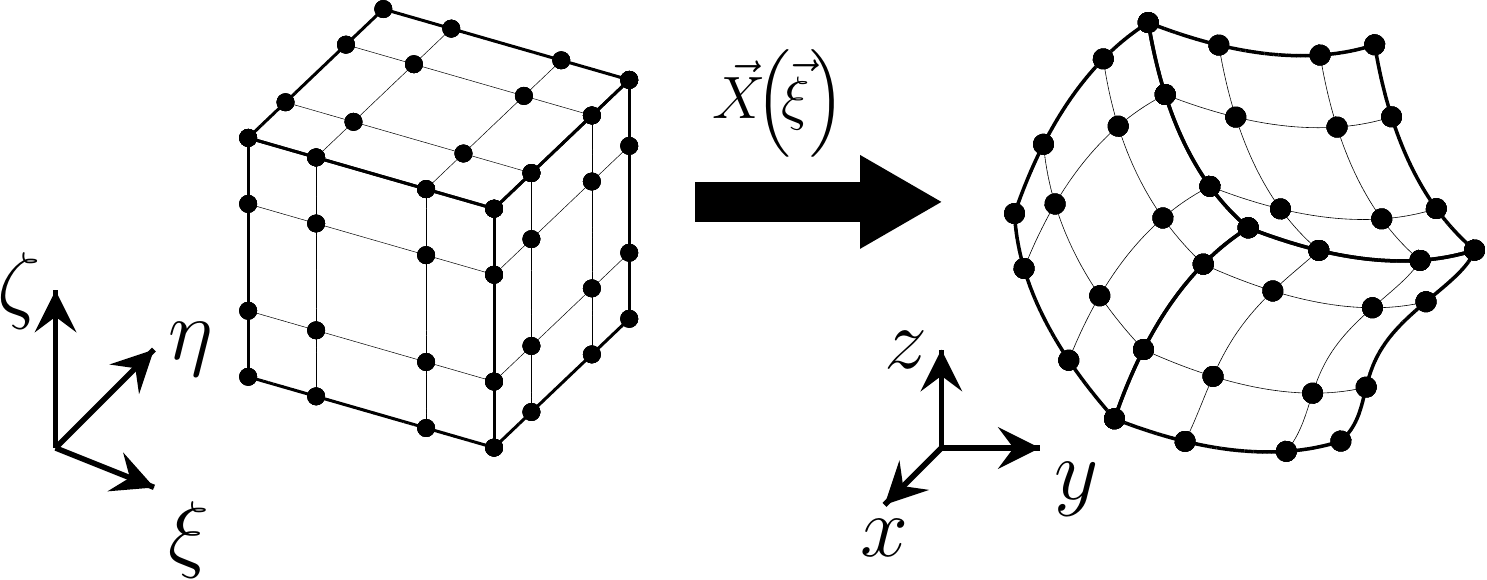}
  \caption{Elements geometrical transformation from the reference element $E=[-1,1]^{3}$ to their final shape and position on the physical space. 
  The transformation uses a transfinite order $N$ mapping $\svec{X}\left(\svec{\xi}\right)$. The tensor product Gauss--Lobatto points are 
  also transformed}
  \label{fig:DG:mapping}
\end{figure}

The system of equations \eqref{eq-xpipe:continuous:non-conservative-pde} 
is approximated in space with a high--order discontinuous Galerkin method. The 
computational domain $\Omega$ is tessellated in non--overlapping hexahedral 
elements $e$. Both the solution and the geometry are approximated by order $N$ 
polynomials, and thus the elements can be curvilinear. We establish a transfinite 
mapping as in \cite{2006:Kopriva,2019:Manzanero-CH} that transforms the unit cube $E=[-1,1]^{3}$ (known as the \textit{reference element}) 
to an arbitrarily--shaped hexahedral element $e$ (see Fig.~\ref{fig:DG:mapping}). 
Thus, we work on the refence space $(\xi,\eta,\zeta)$ where the 
reference element $E$ is defined, and then its coordinates are related to the physical 
space with the mapping $(x,y,z)=\svec{X}\left(\xi,\eta,\zeta\right)$. 

We define a set of tensor 
product Gauss--Lobatto (GL) points $\left(\xi_i,\eta_j,\zeta_k\right)_{i,j,k=0}^{N}$ \cite{2009:Kopriva}, 
which we use to approximate the solution by order $N$ polynomials,
\begin{equation}
\stvec{q}\bigr|_{e}\approx \mathrm{I}^{N}\left(\stvec{q}\right)=\stvec{Q}=\sum_{i,j,k=0}^{N}\stvec{Q}_{ijk}(t)l_i\left(\xi_{i}\right)l_j\left(\eta_j\right)l_k\left(\zeta_k\right),
\label{eq:DG:N-order-interpolation}
\end{equation}
where $\stvec{Q}_{ijk}$ are the nodal coefficients, and $l_i$ are the Lagrange 
interpolating polynomials. The GL points also define quadrature rules that approximate the integrals on the reference 
element,
\begin{equation}
  \left\langle f,g\right\rangle_{E}\approx \left\langle 
  F,G\right\rangle_{E,N}=\sum_{i,j,k=0}^{N}w_{ijk}F_{ijk}G_{ijk},
\end{equation}
where $w_{ijk}=w_{i}w_{j}w_{k}$ are the tensor product quadrature weights 
\cite{2009:Kopriva}.

We define the covariant and contravariant basis to relate the derivatives in the reference ($\nabla_{\xi}$) and physical  ($\nabla$)
spaces, 
\begin{equation}
  \svec{a}_{i} =\frac{\partial\svec{X}}{\partial \xi^{i}},\quad 
  \svec{a}^{i}=\nabla\xi^{i} = \frac{\svec{a}^{j}\times\svec{a}^{k}}{J},\quad 
  J=\svec{a}_{1}\cdot\left(\svec{a}_{2}\times\svec{a}_{3}\right),\quad 
  \left(i,j,k\right)~\mathrm{cyclic},
\end{equation}
which allow us to compute the gradient of a scalar and the divergence of a vector as,
\begin{equation}
  J\nabla\stvec{w}_i = \mathcal{M}\nabla_{\xi}\stvec{w}_i,\quad 
  J\nabla\cdot\ssvec{f}_i=\nabla_{\xi}\cdot\cssvec{{f}}_i,\quad 
  \cssvec{f}_i=\mathcal{M}^{T}\ssvec{f}_i,\quad
  \mathcal{M}=\left(J\svec{a}^{\xi},J\svec{a}^{\eta},J\svec{a}^{\zeta}\right).
  \label{eq:DG:grad-and-div}
\end{equation}

Discretely, the mapping is approximated with the order $N$ interpolation (see \eqref{eq:DG:N-order-interpolation}), which is then differentiated to get the 
discrete 
covariant basis and Jacobian, $\mathcal J$. The contravariant basis, however, is 
computed using a curl form \cite{2006:Kopriva},
\begin{equation}
  \mathcal J{a}^{i}_{n}=-\hat{x}^{i}\cdot\nabla_{\xi}\times 
  \mathrm{I}^{N}\left(X_l\nabla_{\xi}X_{m}\right),\quad 
  i,n=1,2,3,\quad\left(n,m,l\right)~\mathrm{cyclic},
\end{equation}
to fulfill the \textit{discrete metric identities} \cite{2006:Kopriva}, key to ensure the 
free--stream preservation (i.e. the derivative of a constant is zero) on 
curvilinear grids.\\

To construct the scheme we first cast the PDE \eqref{eq-xpipe:continuous:non-conservative-pde} 
as a second order system introducing the auxiliary variable $\ssvec{g}=\nabla 
\stvec{w}$, and transform the operators to the reference space,
\begin{equation}
\begin{split}
J\stvec{q}_{t}+{\nabla}_{\xi}\cdot\cssvec{f}_{e}\left(\stvec{q}\right)&={\nabla}_{\xi}\cdot\cssvec{f}_{v}\left(\stvec{q},\ssvec{g}\right)+J\stvec{s}\left(\stvec{q},\ssvec{g}\right),\\
  J\ssvec{g}&=\mathcal{M}{\nabla}_{\xi}\stvec{w}.
  \end{split}
  \label{eq:DG:system-cast}
\end{equation}

Next, we multiply equation \eqref{eq:DG:system-cast} by two order $N$ polynomial arbitrary test 
functions, integrate over the reference element $E$, and apply the Gauss law on 
the inviscid and viscous fluxes, and on the gradient of $\stvec{w}$,
\begin{equation}
\begin{split}
&\left\langle J\stvec{q}_{t},\stvecg{\phi}\right\rangle_{E}+\int_{\partial e}\stvecg{\phi}^{T}\left(\ssvec{f}_{e}-\ssvec{f}_{v}\right)\cdot\mathrm{d}\svec{S}  
-\left\langle\cssvec{f}_{e},\nabla_{\xi}\stvecg{\phi}\right\rangle_{E}=-\left\langle\cssvec{f}_{v},\nabla_{\xi}\stvecg{\phi}\right\rangle_{E} 
+ \left\langle J\stvec{s},\stvecg{\phi}\right\rangle_{E}, \\
& \left\langle J\ssvec{g},\ssvecg{\varphi}\right\rangle_{E}=\int_{\partial e} \stvec{w}^{T}\ssvecg{\varphi}\cdot\mathrm{d}\svec{S}-\left\langle \stvec{w},\nabla_{\xi}\cdot\cssvecg{\varphi}\right\rangle_{E}.
  \end{split}
  \label{eq:DG:system-weak}
\end{equation}

Now we replace the polynomial ansatz. The functions are approximated by 
polynomials, and the integrals by quadratures. As a result of the disconnection
between adjacent elements, the solution can be discontinuous across the 
inter--element faces, and the fluxes at the surface integrals are not defined. 
Thus, we use a uniquely defined \textit{numerical flux} at the surface 
integrals, $\ssvec{f}\approx\ssvec{f}^{\star}\left(\stvec{q}_{\mathrm{L}},\stvec{q}_{\mathrm{R}}\right)$ that depend on the two 
neighbouring states,
\begin{equation}
\begin{split}
&\left\langle \mathcal J\stvec{Q}_{t},\stvecg{\phi}\right\rangle_{E,N}+\int_{\partial e,N}\stvecg{\phi}^{T}\left(\ssvec{F}^{\star}_{e}-\ssvec{F}^{\star}_{v}\right)\cdot\mathrm{d}\svec{S}  
-\left\langle\cssvec{F}_{e},\nabla_{\xi}\stvecg{\phi}\right\rangle_{E,N}=-\left\langle\cssvec{F}_{v},\nabla_{\xi}\stvecg{\phi}\right\rangle_{E,N} 
+ \left\langle \mathcal J\stvec{S},\stvecg{\phi}\right\rangle_{E,N}, \\
& \left\langle \mathcal J\ssvec{G},\ssvecg{\varphi}\right\rangle_{E,N}=\int_{\partial e,N} \stvec{W}^{\star,T}\ssvecg{\varphi}\cdot\mathrm{d}\svec{S}-\left\langle \stvec{W},\nabla_{\xi}\cdot\cssvecg{\varphi}\right\rangle_{E,N}.
  \end{split}
  \label{eq:DG:system-weak-disc}
\end{equation}

Lastly, we enhance the robustness of this implementation by using a split--form 
scheme \cite{2013:Gassner,2016:Gassner}. To do so, we apply a second time the Gauss law on the inviscid fluxes, 
\begin{equation}
\begin{split}
&\left\langle \mathcal J\stvec{Q}_{t},\stvecg{\phi}\right\rangle_{E,N}+\int_{\partial e,N}\stvecg{\phi}^{T}\left(\ssvec{F}^{\star}_{e}-\ssvec{F}_{e}-\ssvec{F}^{\star}_{v}\right)\cdot\mathrm{d}\svec{S}  
+\left\langle\mathbb{D}\left(\cssvec{F}_{e}\right),\stvecg{\phi}\right\rangle_{E,N}=-\left\langle\cssvec{F}_{v},\nabla_{\xi}\stvecg{\phi}\right\rangle_{E,N} 
+ \left\langle \mathcal J\stvec{S},\stvecg{\phi}\right\rangle_{E,N}, \\
& \left\langle \mathcal J\ssvec{G},\ssvecg{\varphi}\right\rangle_{E,N}=\int_{\partial e,N} \stvec{W}^{\star,T}\ssvecg{\varphi}\cdot\mathrm{d}\svec{S}-\left\langle \stvec{W},\nabla_{\xi}\cdot\cssvecg{\varphi}\right\rangle_{E,N},
  \end{split}
  \label{eq:DG:system-weak-disc-split}
\end{equation}
where $\mathbb{D}\left(\ssvec{F}_{e}\right)$ is a split--form approximation of the 
divergence $\nabla_{\xi}\cdot\cssvec{F}_{e}$ that uses a two--point flux 
$\ssvec{F}^{\#}_{e}$,
\begin{equation}  
  \mathbb{D}\left(\ssvec{F}_{e}\right)_{ijk}=2\sum_{m=0}^{N}\left(D_{im}\cssvec{F}_{e}\left(\stvec{Q}_{ijk},\stvec{Q}_{mjk}\right)+D_{jm}\cssvec{G}_{e}\left(\stvec{Q}_{ijk},\stvec{Q}_{imk}\right)+D_{km}\cssvec{H}_{e}\left(\stvec{Q}_{ijk},\stvec{Q}_{ijm}\right)\right),
\end{equation}
with $D_{ij}=l_{j}'\left(\xi_{i}\right)$. For this work, we adapt the 
two--point flux derived in \cite{2019:Manzanero-iNS} for the incompressible Navier--Stokes to 
the system solved herein. For the last four equations we simply copy the 
two--point flux from \cite{2019:Manzanero-iNS}, and then we perform the product of the averages 
for the first two equations:
\begin{equation}
  \ssvec{F}^{\#}_{e} = \left(\begin{array}{ccc} \aver{c_1}\aver{u} & \aver{c_1}\aver{v} & \aver{c_1}\aver{w} \\ 
  \aver{c_2}\aver{u} & \aver{c_2}\aver{v}  & \aver{c_2}\aver{w} \\ 
  \aver{\rho}\aver{u}^2 + \aver{p} &   \aver{\rho}\aver{u}\aver{v} &   \aver{\rho}\aver{u}\aver{v}\\ 
    \aver{\rho}\aver{u}\aver{v}  &   \aver{\rho}\aver{v}^{2} + \aver{p}&   \aver{\rho}\aver{v}\aver{w}\\ 
      \aver{\rho}\aver{u}\aver{w}  &   \aver{\rho}\aver{v}\aver{w} &   \aver{\rho}\aver{w}^2+ \aver{p}\\ 
\rho_0c_0^2\aver{u} & \rho_0c_0^2\aver{v} & \rho_0c_0^2\aver{w}
  \end{array}\right),\quad \cssvec{F}_{e,i}^{\#}=\aver{\mathcal 
  M}\ssvec{F}_{e,i}^{\#},
\end{equation}
where the brackets represent the average between the two states,
\begin{equation}
  \aver{u} = \frac{u_{i}+u_{m}}{2}.
\end{equation}

The approximation of the equations is completed with the approximation of 
the chemical potentials, which are the first two gradient variables $w_{1},w_{2}$.  
To do so, we proceed as in  the PDE: we cast the definition of the 
chemical potentials introducing auxiliary variables 
$\svec{g}_{c,i}=\nabla{c}_i$, we transform the 
differential operators to the reference space, we construct weak forms within 
the elements, and integrate the volume terms with differential operators to get,
\begin{equation}
\begin{split}
\left\langle \mathcal J\mu_{i},{\phi}\right\rangle_{E,N}&=\left\langle\mathcal J\frac{12}{\varepsilon}{\Sigma}_{i}{F}_{i},{\phi}\right\rangle_{E,N}-\frac{3}{4}\varepsilon\Sigma_{i}\int_{\partial e,N}{\phi}\left(\svec{G}^{\star}_{c,i}\right)\cdot\diff \svec{S} +\frac{3}{4}\varepsilon \Sigma_{i}\left\langle \svec{\tilde G}_{c,i},\svec{\nabla}_{\xi}{\phi}\right\rangle_{E,N},\\
\left\langle \mathcal J\svec{G}_{c,i},\svec{\varphi}\right\rangle_{E,N} &=\int_{\partial e,N}{C}^{\star,T}_{i}\svec{\varphi}\cdot\diff\svec{S}
-\left\langle {C}_i,\svec{\nabla}_{\xi}\cdot\svec{\tilde \varphi}_{G_{c}}\right\rangle_{E,N}.
\end{split}
\end{equation}

From the concentration field ($C_1,C_2$), we compute the chemical potentials 
($\mu_1,\mu_2$), which are then introduced in the entropy variables $\stvec{W}$, 
to compute their gradient $\ssvec{G}$, and the state vector time derivative $\stvec{Q}_{t}$. 
In the next sections, we describe the computation of the numerical fluxes for 
inter--element and physical boundary faces.

\subsubsection{Numerical fluxes}\label{sec-xpipe:DGSEM:NumericalFluxes}

The numerical fluxes couple the inter--element solutions through an uniquely 
defined value for the surface integrals. For the inviscid Riemann solver $\ssvec{F}^{\star}_{e}$ we compute the exact Riemann problem solution derived in \cite{2017:Bassi}
for the normal velocity $u^\star$ and 
the pressure $p^\star$, 
\begin{equation}
\begin{split}
&u^{\star} = \frac{p_R-p_L +\rho_R u_R \lambda^{-}_{R}-\rho_L u_L \lambda^{+}_L}{\rho_R 
	\lambda_R^{-}-\rho_L\lambda_L^{+}},~~p^{\star}=p_R + \rho_R 
\lambda_R^{-}(u_R-u^\star), \\
&c_1^{\star},c_2^{\star},v^{\star},w^{\star} = \left\{\begin{array}{lll}
c_{1L},c_{2L},v_L,w_L& \text{if} & u^{\star} \geqslant 0 
\\
c_{1R},c_{2R},v_R,w_R& \text{if} & u^{\star} < 
0
\end{array}\right.,
\end{split}
\label{eq-xpipe:Riemann:star-solution}
\end{equation}
and we compute the tangential velocities $v^\star,w^\star$ and concentrations $c_1^\star,c_2^\star$ from the appropriate 
element depending on the sign of $u^{\star}$.
Amongst the beneficial properties of this Riemann solver that led to this choice, we highlight 
that it is physical (it is the exact solution of the Riemann problem), efficient (e.g. does not need any iterative solution) and 
parameter--free. 

For viscous fluxes and gradient variables, $\ssvec{F}_{v}^{\star}$ and $\stvec{W}^{\star}$, 
we use the Symmetric Interior Penalty (SIP) method, 
\begin{equation}
\stvec{W}^{\star} = \aver{\stvec{W}},~~\ssvec{F}_{v}^{\star} = \aver{\ssvec{F}_{v}\left(\stvec{Q},\svec{\nabla}\stvec{W}\right)}+\beta\left(\begin{array}{c}\frac{M_0}{\Sigma_{1}}\jump{\mu_{1}}\\
\frac{M_0}{\Sigma_{2}}\jump{\mu_{2}}\\
\mu\jump{\rho \svec{u}}\\
\mu\jump{p}
\end{array}\right)\svec{n}_{L},
\end{equation}
which uses the local gradient \eqref{eq:DG:grad-and-div} for $\svec{\nabla}\stvec{W}$. We take the penalty parameter 
from \cite{Shahbazi:2005} 
\begin{equation}
\beta=\frac{\left(N+1\right)\left(N+2\right)}{2\bar{h}},~~\bar{h} 
= \frac{\min\left(V_{L},V_{R}\right)}{S},
\end{equation}
where $\bar{h}$ is yields an approximated measure of the minimum element size normal to the 
face, $V_{L}$ and $V_{R}$ are the volumes of the neighbouring elements, and $S$ 
is the area of the face. Finally, for the concentration and its gradient we 
also use the SIP method,
\begin{equation}
{C}_{i}^{\star} = \aver{{C}_i},~~\svec{G}_{c,i}^{\star} = \aver{\svec{\nabla} 
  {C}_i}-\beta\left({C}_{i,L}\svec{n}_{L}+{C}_{i,R}\svec{n}_{R}\right).
  \end{equation}

\subsection{Time discretization}\label{sec-xpipe:TimeDiscretization}

The  fourth order spatial derivative in the Cahn--Hilliard equation is too stiff to be solved explicitly.  
Unless the mobility parameter is low enough to allow reasonable time--steps, we use an IMplicit--EXplicit (IMEX)
method to integrate in time: Navier--Stokes terms are solved using a 
third--order explicit Runge--Kutta (RK3) method, the Cahn--Hilliard chemical free--energy 
is solved using forward Euler and the Cahn--Hilliard interface energy is solved 
using backward Euler. This first order IMEX scheme was adapted from that derived in \cite{dong2018multiphase} in the context
of $N$--phase flows. 

This IMEX procedure has two steps, which we describe in a semi--discrete fashion (continuous in space, discrete in 
time):

\begin{enumerate}
\item We perform the explicit RK3 step for the Navier--Stokes terms, without the contribution 
from the chemical--free energy and interfacial energy in the Cahn--Hilliard 
equation,

\begin{equation}
\resizebox{0.92\textwidth}{!}{$
\begin{split}
\left(\begin{array}{c}c_1 \\ c_2 \\ \rho \svec{u} \\ p\end{array}\right)_t =- \nabla\cdot\left(\begin{array}{c} c_1\svec{u} \\ c_2\svec{u}\\
\rho\svec{u}\svec{u} + p\tens{I}_{3} \\ \rho_0 c_0^2\svec{u}
\end{array}\right) + \nabla\cdot\left(\begin{array}{c}  0  \\ 0 \\ \eta\left(\nabla\svec{u}+\nabla\svec{u}^{T}\right)\\0
\end{array}\right)
+ \left(\begin{array}{c}
0 \\
0 \\
\rho \svec{g} +\displaystyle{\sum_{m=1}^{3}}\mu_m\nabla c_m\\
0
\end{array}  \right).
\end{split}  $}
\end{equation}
After the RK3 time--step, the variables are called
$\left(\hat{c}_{1},\hat{c}_{2},\rho\svec{u}^{n+1},p^{n+1}\right)$, since the 
concentrations need a correction step to include the chemical potential.

\item We compute a correction step on the two concentrations to solve the Cahn--Hilliard 
equations.
The chemical free--energy terms are solved explicitly (i.e. evaluated in 
$c_{i}^{n}$) and the stiff interfacial energy terms implicitly (i.e. evaluated in 
$c_{i}^{n+1}$). Besides, we also introduce the stabilizing term 
$S_0\left(c_{i}^{n+1}-c_{i}^{n}\right)$, being $S_0$ a constant.
Finally, the approximation of the time derivative is ${c_{t}\approx \left(c_{i}^{n+1}-\hat{c}_{i}\right)/\Delta 
t}$. As a result, the IMEX correction step for each concentration is:
\begin{equation}
  \frac{c^{n+1}_{i}-\hat{c}_{i}}{\Delta t} = M_0\svec{\nabla}^{2}\left(\frac{12}{\varepsilon}f_{i}(c^{n}_{1},c^{n}_{2},c^{n}_{3})+S_{0}\left(c_{i}^{n+1}-c_{i}^{n}\right)-\frac{3}{4}\varepsilon\svec{\nabla}^{2}c_{i}^{n+1}\right).
\end{equation}

Two decoupled linear systems for the two concentration parameters are solved. 
However, the linearity of the operator produces a constant in time Jacobian, 
which also are identical for both phases. This properties led us to compute the 
solution to the linear problem with LU factorization (performed only once at the preprocessing) 
and Gauss substitution. The latter, after the preprocessing, gets similar 
computational times to an iteration in a explicit method.

\end{enumerate}

\subsection{Boundary conditions}\label{sec-xpipe:DGSEM:BoundaryConditions}

In this section we describe the imposition of inflow, outflow, and no--slip wall 
boundary conditions. In this work, we prescribe the boundary conditions weakly. 
Hence, we create a ghost (exterior) state with the appropriate boundary 
information, and then compute the interface fluxes between the interior $\stvec{Q}^{i}$ and 
exterior $\stvec{Q}^{e}$  states.

\subsubsection{Inflow boundary condition}

For the inflow boundary condition, we specify the inflow concentration $c_{i,\mathrm{inflow}}(\svec{x};t)$ and 
the velocity $\svec{u}_{\mathrm{inflow}}\left(\svec{x};t\right)$. For the inviscid fluxes, we construct a ghost 
state,
\begin{equation}
\stvec{Q}^{e} = \left(\begin{array}{c}
C_{1,\mathrm{inflow}} \\
C_{2,\mathrm{inflow}}\\
\rho\left(C_{\mathrm{inflow}}\right)  \svec{U}_{\mathrm{inflow}}\\
P
\end{array}\right),
\label{eq-xpipe:bcs:inflow-bc}
\end{equation}
where we take the pressure $P$ from the interior, and compute the interface flux 
from the exact Riemann problem solution (see   \eqref{eq-xpipe:Riemann:star-solution}), 
$\ssvec{F}_{e}^{\star}\left(\stvec{Q}_{e}^{i},\stvec{Q}_{e}^{e}\right)$.

We defined the viscous fluxes as:
\begin{equation}
\stvec{W}^{\star} = 
\frac{\stvec{W}^{i}+\stvec{W}^{e}}{2},~~\stvec{F}_{v}^{\star}\cdot\svec{n}
= \left(\begin{array}{c} 0 \\ 0 \\ \eta\left(\nabla\svec{U}+\nabla\svec{U}^{T}\right)\cdot\svec{n} \\ 0 
\end{array}\right),
\label{eq-xpipe:DG:inflow-W}
\end{equation}
where we apply the Neumann boundary conditions to the chemical potential, 
and take the interior values for viscous stresses. In \eqref{eq-xpipe:DG:inflow-W}, 
we compute the gradient variables from the ghost state,
${\stvec{W}^{e}=\stvec{W}\left(\stvec{Q}^{e}\right)}$. 

The implementation features an automatic method for the distribution of the 
phases in a circular section, given the superficial velocities and the slip 
velocities. The superficial velocity of each phase is, for a given flow rate, the 
equivalent velocity obtained if the phase occupies the entire section,
\begin{equation}
  v_{s,i} = 
  \frac{1}{A_{\mathrm{inflow}}}\int_{\mathrm{inflow}}\svec{u}_{i}\cdot\diff\svec{S}.
\end{equation}

\begin{figure}[h]
\centering
\includegraphics[width = 0.5\textwidth]{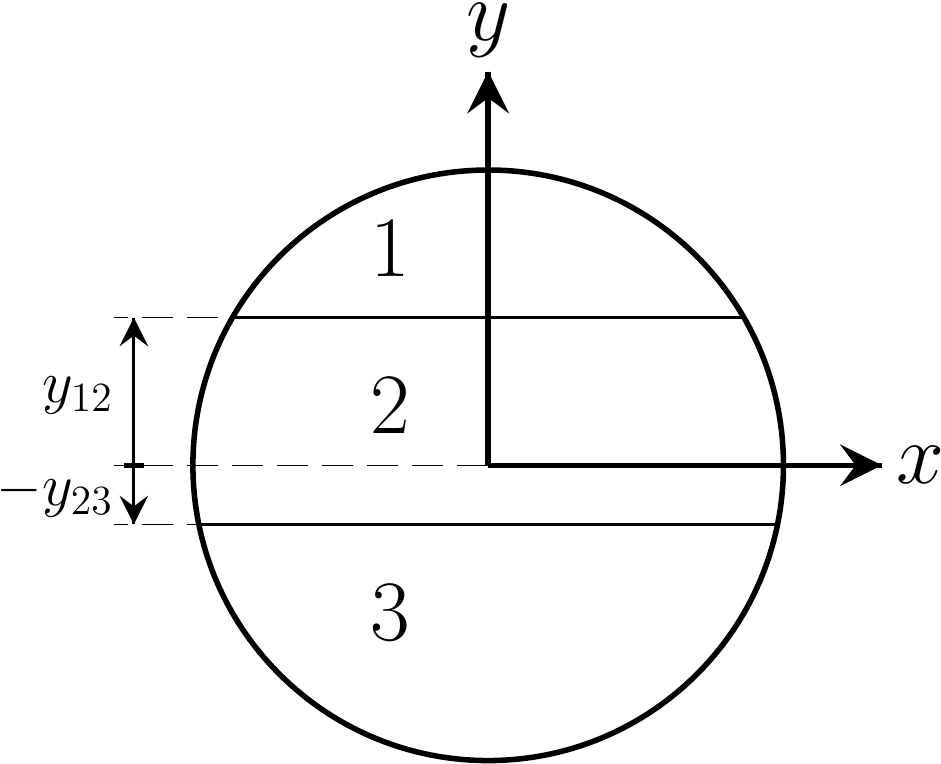}
\caption{Configuration of a layered inflow. The position of the two interfaces is provided by the values $y_{12}$ 
and $y_{23}$} 
\label{fig-xpipe:inflow-conf}
\end{figure}

Furthermore, because real flow configurations feature a large ratio of the 
superficial velocities that confine one of the phases to the near wall region, we 
allow the velocity to be discontinuous at the interface between the phases. By 
doing so, we can use lower velocities for the phases with smaller superficial 
velocity, so that they occupy a larger region of the cross section. The 
problem with a phase being very confined to the near wall region is that it might not 
be well captured by the numerical method. For the flow configuration given in 
Fig.~\ref{fig-xpipe:inflow-conf},
the concentration inflow boundary condition is
\begin{equation}
  \begin{split}
  c_{1,\mathrm{inflow}}(x,y) &= 
  \frac{1}{2}+\frac{1}{2}\tanh\left(\frac{y-y_{12}}{\varepsilon}\right),
  \\
  c_{2,\mathrm{inflow}}(x,y)&= 
  \frac{1}{2}\tanh\left(\frac{y-y_{23}}{\varepsilon}\right)-\frac{1}{2}\tanh\left(\frac{y-y_{12}}{\varepsilon}\right),
  \end{split}
\end{equation}
and the velocities are computed from a Poiseuille flow,
\begin{equation}
u_{\mathrm{inflow}}=\left(V_{1,\max}c_{1,\mathrm{inflow}} + V_{2,\max}c_{2,\mathrm{inflow}} + 
V_{3,\max}c_{3,\mathrm{inflow}}\right)\left(1-\left(\frac{r}{R}\right)^2\right),
\end{equation}
with the two slip velocities, which are user--input,
\begin{equation}
  V_{s,12} = V_{1,\max}-V_{2,\max},~~  V_{s,23} = V_{2,\max}-V_{3,\max}.
  \label{eq-xpipe:DG:slip-vels}
\end{equation}

Therefore, there are five unknowns ($V_{1,\max}$, $V_{2,\max}$, $V_{3,\max}$, $y_{12}$, and 
$y_{23}$), and five equations: the two slip--velocities definitions (see \eqref{eq-xpipe:DG:slip-vels}), and the 
superficial velocities,
\begin{equation}
  v_{s,i} = \frac{V_{i,\max}}{A_{\mathrm{inflow}}}\int_{\mathrm{inflow}}c_{i,\mathrm{inflow}}\left(1-\left(\frac{r}{R}\right)^2\right)\diff 
  S,
  \label{eq-xpipe:bcs:superficial-velocities}
\end{equation}
that are solved using a Newton--Rhapson method.

\subsubsection{Outflow boundary condition}

The outflow boundary condition specifies the ambient pressure at the exit of the 
domain $P_{o}$, and applies a Neumann boundary condition to the rest of the variables. 
Therefore, for the inviscid fluxes the exterior state is defined as:
\begin{equation}
  \stvec{Q}^{e} = \left(\begin{array}{c}C_1 \\ C_2 \\ \rho\svec{U} \\ 
  P_o\end{array}\right),
  \label{eq-xpipe:bcs:outflow}
\end{equation}
whereas for viscous fluxes we simply use a Neumann boundary condition for all 
the variables,
\begin{equation}
\stvec{W}^{\star} = 
\stvec{W}_{i},~~\stvec{F}_{v}^{\star}\cdot\svec{n} = 0.
\label{eq-xpipe:DG:outflow-W}
\end{equation}

\subsubsection{No--slip wall boundary condition}

We construct a ghost state with the same variables as the 
inside, but changing the sign of the normal velocity,

\begin{equation}
\stvec{Q}^{e}= \left(\begin{array}{c}
c_1 \\ c_2 \\
\rho \left(\svec{U}-2\left(\svec{U}\cdot\svec{n}\right)\svec{n}\right)\\
p
\end{array}\right).
\label{eq-xpipe:bcs:wall}
\end{equation}

For the viscous numerical fluxes, we apply Neumann boundary conditions in all 
variables except velocities, which take the interior values,

\begin{equation}
\stvec{W}^{\star} = \frac{\stvec{W}^{i} + \stvec{W}^{e}}{2},~~\ssvec{F}_{v}^{\star}\cdot\svec{n} = \left(\begin{array}{c}
0 \\
0\\
\eta\left(\nabla\svec{U}+\nabla\svec{U}^{T}\right)\cdot\svec{n} \\
0
\end{array}  \right).
\end{equation}

Finally, for the gradient of the concentrations $\svec{G}_{c,i}^{\star}$ the Neumann 
boundary condition is non--homogeneous if one wants to solve for arbitrary wall 
contact angles. Thus, we follow \cite{2014:Shi,2020:Manzanero-UR-CaF} and use the following expression:
\begin{equation}
  {C}_i^{\star}= {C}_i,~~\svec{G}_{c,i}^{\star}\cdot\svec{n}= 
  {F}_{w,i},
 \end{equation}
 where the boundary coefficients $F_{w,i}$ ar.e
 \begin{equation}
   \begin{split}
     F_{w,1} &= 
     -\frac{4}{\varepsilon}\left(\cos\theta_{12}^{w}C_1C_2\left(C_1+C_2\right)+\cos\theta_{13}^{w}C_1C_3\left(C_1+C_3\right)\right),\\
     F_{w,2}&=-\frac{4}{\varepsilon}\left(-\cos\theta_{12}^{w}C_1C_2\left(C_1+C_2\right)+\cos\theta_{23}^{w}C_2C_3\left(C_1+C_3\right)\right),
   \end{split}
 \end{equation}
 being the three wall contact angles $\theta_{ij}^{w}$ related by the wall equilibrium constraint  \cite{2014:Shi,2020:Manzanero-UR-CaF}, ${\sigma_{12} \cos \theta_{12}^{w} + \sigma_{23}\cos{\theta}_{23}^{w} = \sigma_{13}\cos{\theta}_{13}^{w}}$.
 For $90^{\circ}$ angles, the 
 coefficients $F_{w,i}$ are zero, and the boundary condition reduces to 
 homogeneous Neumann.

\section{Two--phase simulations}\label{sec-xpipe:Validation}

We perform the validation of the solver in the particular case of solving a 
two--phase flow (i.e. as described in Sec.~\ref{sec-xpipe:3ph-as-2ph}). We 
first solve a manufactured solution in Sec.~\ref{sec-xpipe:num2ph:Mansol}, and 
then we solve a two--phase horizontal pipe in Sec.~\ref{sec-xpipe:2PH-pipe}. The 
enhancement of the robustness provided by the split--form scheme is addressed in 
Sec.~\ref{sec-xpipe:num2ph:robustness-split}.

\subsection{Manufactured solution}\label{sec-xpipe:num2ph:Mansol}

We first study the convergence properties of the method. To do this,
we borrow the manufactured solution from a previous two--phase work \cite{2019:Manzanero-MU}. 
This two--phase manufactured solution is defined as:
\begin{equation}
\begin{split}
c_{1,0}(x,y;t) &= \frac{1}{2}\left(1 + \cos\left(\pi x\right)\cos\left(\pi y\right)\sin\left(t\right)\right), \\
c_{2,0}\left(x,y;t\right) &= 0, \\
u_0(x,y;t) &=2\sin\left(\pi x\right) \cos\left(\pi z\right)\sin\left(t\right), \\
v_0(x,y;t) &=-2\cos\left(\pi x\right)\sin\left(\pi y\right) \sin\left(t\right), \\
p_0(x,y;t) &=2\sin\left(\pi x\right)\sin\left(\pi z\right)\cos\left(t\right),
\end{split}
\label{eq-xpipe:num:conv:mansol}
\end{equation}
which we solve on the domain $(x,y)\in[-1,1]^2$~m. The final time is $t_F=0.1$~s, and all the physical parameters are presented in Table \ref{tab-xpipe:num:conv:param}. 
{\begin{table}[h]
		\centering
		\caption{List of the parameter values used with the manufactured solution (see  \eqref{eq-xpipe:num:conv:mansol})}
		\label{tab-xpipe:num:conv:param}
			\begin{tabular}{llllllll}
				\hline
				$\rho_1$& $\rho_3$ ($\text{kg}/\text{m}^3$)   & $\eta_1$ & $\eta_2$ (Pa$\cdot$s) & $\varepsilon$ (m)& $M_0$ (m/s)& $c_0^2$ (m/s$^2$)$^2$& $\sigma$ (N/m) \\ \hline
				1.0 & 2.0 & 1.0E-3 & 1.0E-3 & $1/\sqrt{2}$ & 1.134E-2 & 1.0E3 & 6.236E-3 \\
				\hline
		\end{tabular}
	\end{table}}

\begin{figure}[h]
  \centering
  \subfigure[$\Delta t=10^{-4}$~s]{\includegraphics[width=0.49\textwidth]{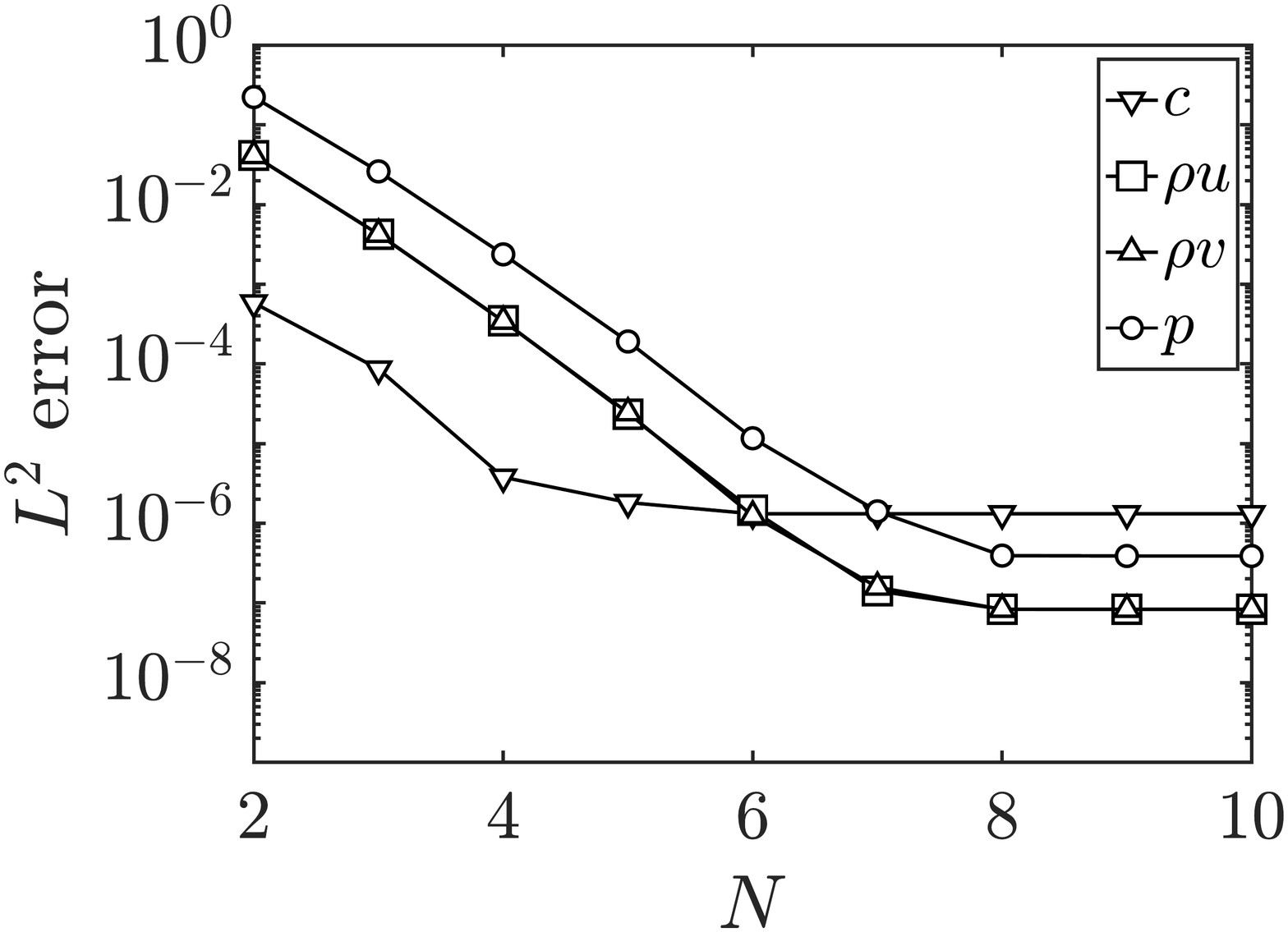}}
  \subfigure[$\Delta t=10^{-5}$~s]{\includegraphics[width=0.49\textwidth]{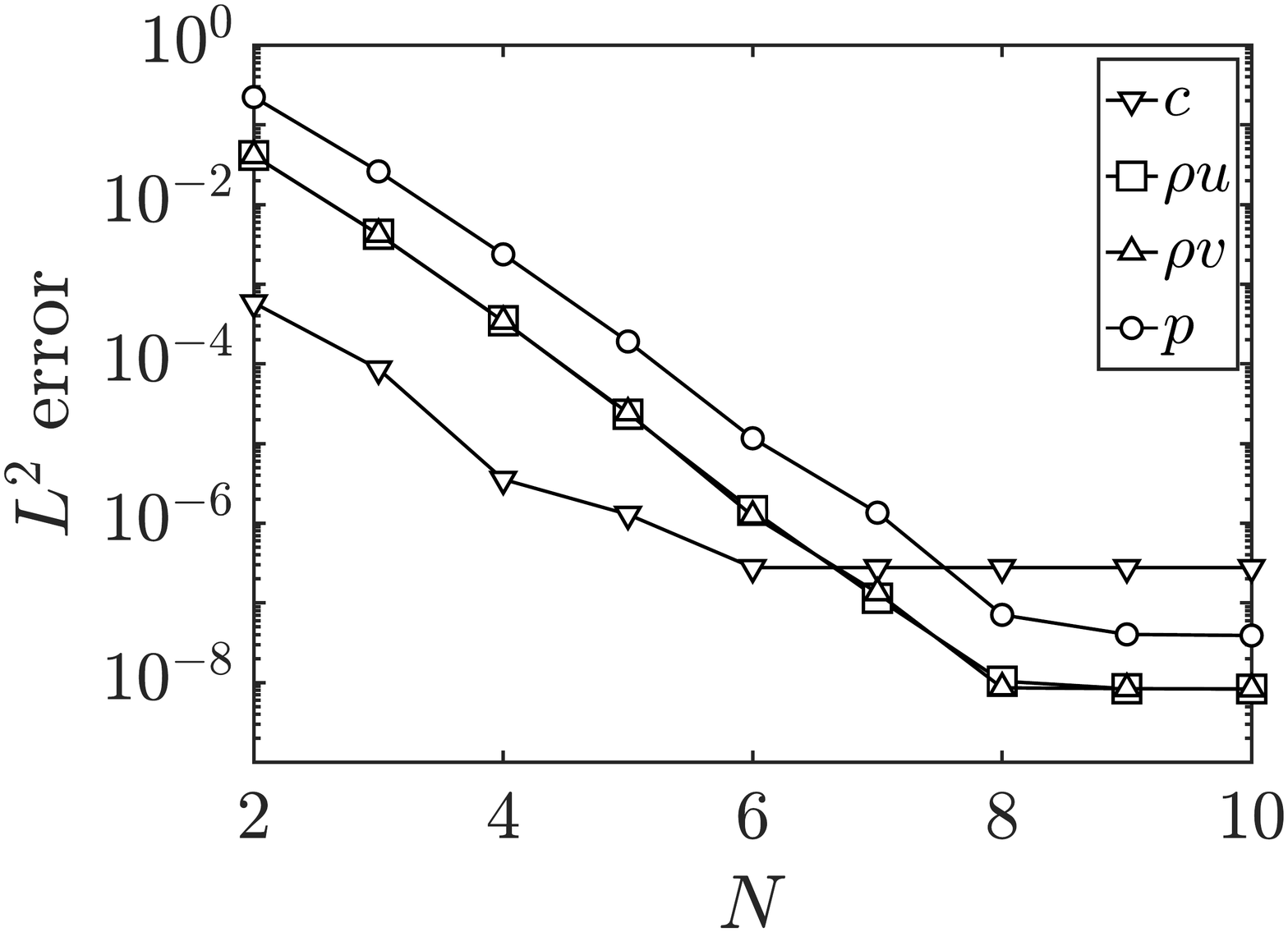}}
  \caption{Two--phase solver: polynomial order convergence study of the manufactured solution \eqref{eq-xpipe:num:conv:mansol}. 
  We represent the L$^2$ errors in concentration, $x$-- and $y$--momentum, and pressure. 
  The polynomial order ranges from 2 to 10, and we integrate in time until $t_{F}=0.1$~s with two time step sizes: $\Delta t=10^{-4}$~s and $10^{-5}$~s. All physical parameters are given in Table \ref{tab-xpipe:num:conv:param}}
  \label{fig-xpipe:num:conv}  
\end{figure}

We first present the polynomial order convergence analysis in 
Fig.~\ref{fig-xpipe:num:conv} for two values of the $\Delta t$ parameter, $\Delta t=10^{-4}$~s and $10^{-5}$~s. 
We 
consider a mesh with $4^2$ elements, and the polynomial order ranges from $N=2$ to $N=10$. We 
represent the L$^2$ errors for the concentration, momentum, and pressure. For 
both time--step values, the errors behave similarly, with a space 
under--resolved region for lower polynomial orders, where the errors decrease 
exponentially, and a time under--resolved region where the error stagnates with 
the polynomial order. We see that the use of a first order Euler scheme for the 
Cahn--Hilliard part impacts its accuracy, as the stagnation is reached in low 
polynomial orders (i.e. where temporal errors dominate).
However, we believe that for an industrial solver, to achieve exponential convergence in
the Cahn--Hilliard equation is not critical. The 
role of the Cahn--Hilliard equation for an industrial simulation is to advect the phases, to
introduce interface regularization that helps to solve under--resolved 
simulations, and to separate the phases (or avoid the mixing). We support this argument 
by adding that the mobility is adjusted in practice by trial and error. 

The flow variables are solved with the third order RK3 
method, which leads to more accurate solutions for the momentum and 
pressure equations. Overall, we find a satisfactory convergence behavior of the 
method for smooth solutions.

We also perform a mesh convergence study, presented in 
Table~\ref{tab-xpipe:num:conv:h-conv}. 
\begin{table}[h]
\centering
  \caption{Two--phase solver: manufactured solution \eqref{eq-xpipe:num:conv:mansol} convergence analysis: we use five meshes with $4^2$, $6^2$, $8^2$, $12^2$ and $16^2$  meshes, and $N=2,3,4$ and 5.
The final time is $t_{F}=0.1$~s, and we use the IMEX scheme
  with $\Delta t=5\cdot 10^{-5}$~s}
  \label{tab-xpipe:num:conv:h-conv}
\resizebox{\textwidth}{!}{%
\begin{tabular}{llllllllll}
\hline
      & Mesh   & $c$ error & order & ${\rho} u$ error & order & ${\rho} v$ error & order & $p$ error & order \\ \hline
N=2 & $4^2$ & 5.85E-04 & -- & 4.13E-02 & -- & 4.12E-02 & -- & 2.24E-01 & -- \\ 
 & $6^2$ & 2.61E-04 & 1.99 & 1.50E-02 & 2.49 & 1.50E-02 & 2.49 & 9.09E-02 & 2.22 \\ 
 & $8^2$ & 9.83E-05 & 3.39 & 7.17E-03 & 2.58 & 7.17E-03 & 2.58 & 4.71E-02 & 2.29 \\ 
 & $12^2$ & 2.45E-05 & 3.43 & 2.43E-03 & 2.67 & 2.43E-03 & 2.67 & 1.79E-02 & 2.39 \\ 
 & $16^2$ & 9.93E-06 & 3.13 & 1.10E-03 & 2.76 & 1.10E-03 & 2.76 & 8.71E-03 & 2.50 \\ 
N=3 & $4^2$ & 8.62E-05 & -- & 4.23E-03 & -- & 4.24E-03 & -- & 2.61E-02 & -- \\ 
 & $6^2$ & 1.43E-05 & 4.44 & 9.90E-04 & 3.58 & 9.90E-04 & 3.59 & 7.04E-03 & 3.23 \\ 
 & $8^2$ & 3.52E-06 & 4.87 & 3.40E-04 & 3.72 & 3.40E-04 & 3.72 & 2.69E-03 & 3.34 \\ 
 & $12^2$ & 1.41E-06 & 2.26 & 7.11E-05 & 3.86 & 7.12E-05 & 3.86 & 6.55E-04 & 3.48 \\ 
 & $16^2$ & 1.39E-06 & 0.05 & 2.27E-05 & 3.97 & 2.27E-05 & 3.97 & 2.31E-04 & 3.62 \\ 
N=4 & $4^2$ & 3.81E-06 & -- & 3.47E-04 & -- & 3.43E-04 & -- & 2.37E-03 & -- \\ 
 & $6^2$ & 1.66E-06 & 2.05 & 5.26E-05 & 4.65 & 5.26E-05 & 4.62 & 4.06E-04 & 4.36 \\ 
 & $8^2$ & 1.39E-06 & 0.60 & 1.35E-05 & 4.74 & 1.35E-05 & 4.73 & 1.12E-04 & 4.48 \\ 
 & $12^2$ & 1.39E-06 & 0.01 & 1.90E-06 & 4.83 & 1.91E-06 & 4.83 & 1.74E-05 & 4.59 \\ 
 & $16^2$ & 1.39E-06 & 0.00 & 4.63E-07 & 4.91 & 4.65E-07 & 4.91 & 4.52E-06 & 4.68 \\ 
N=5 & $4^2$ & 1.87E-06 & -- & 2.33E-05 & -- & 2.43E-05 & -- & 1.91E-04 & -- \\ 
 & $6^2$ & 1.39E-06 & 0.74 & 2.36E-06 & 5.64 & 2.37E-06 & 5.75 & 2.07E-05 & 5.48 \\ 
 & $8^2$ & 1.39E-06 & 0.00 & 4.60E-07 & 5.69 & 4.62E-07 & 5.68 & 4.21E-06 & 5.54 \\ 
 & $12^2$ & 1.39E-06 & 0.00 & 6.28E-08 & 4.91 & 6.29E-08 & 4.92 & 4.61E-07 & 5.45 \\ 
 & $16^2$ & 1.39E-06 & 0.00 & 4.26E-08 & 1.34 & 4.26E-08 & 1.35 & 2.10E-07 & 2.74 
 \\
\hline
\end{tabular}}
\end{table}
We use meshes with $4^2$, $6^2$, $8^2$, 
$12^2$, and $16^2$ elements, and we vary the polynomial order from $N=2$ to 
$N=5$. For lower polynomial orders, we find the expected order of convergence for the 
Cahn--Hilliard equation (i.e. even more than $N+1$), but then the early 
stagnation is found as a result of the first--order Euler scheme. For the 
Navier--Stokes part we find that the order of accuracy is systematically between $N$ 
and $N+1$, which is a similar error behavior to that seen previous works \cite{2017:Bassi}.

We conclude that the convergence study confirms that the error behavior of the scheme is the expected one for 
smooth solutions.

\subsection{Pipe simulations}\label{sec-xpipe:2PH-pipe}

In this section, we focus on more practical test cases, relevant for the oil and gas industry. In particular, 
we simulate the flow in horizontal pipes with different superficial velocities (see  \eqref{eq-xpipe:bcs:superficial-velocities}) at the inflow, U$_L^S$ and U$_L^S$, which correspond to the gas and liquid respectively. 
Depending on the superficial velocities values, we obtain different flow 
regimes. 
This test case follows the one proposed in \cite{xie2017direct,gomez2019novel} and agrees with the experimental results of \cite{Taitel1976}. 

\begin{figure}[h]
\centering
\includegraphics[width=0.55\textwidth]{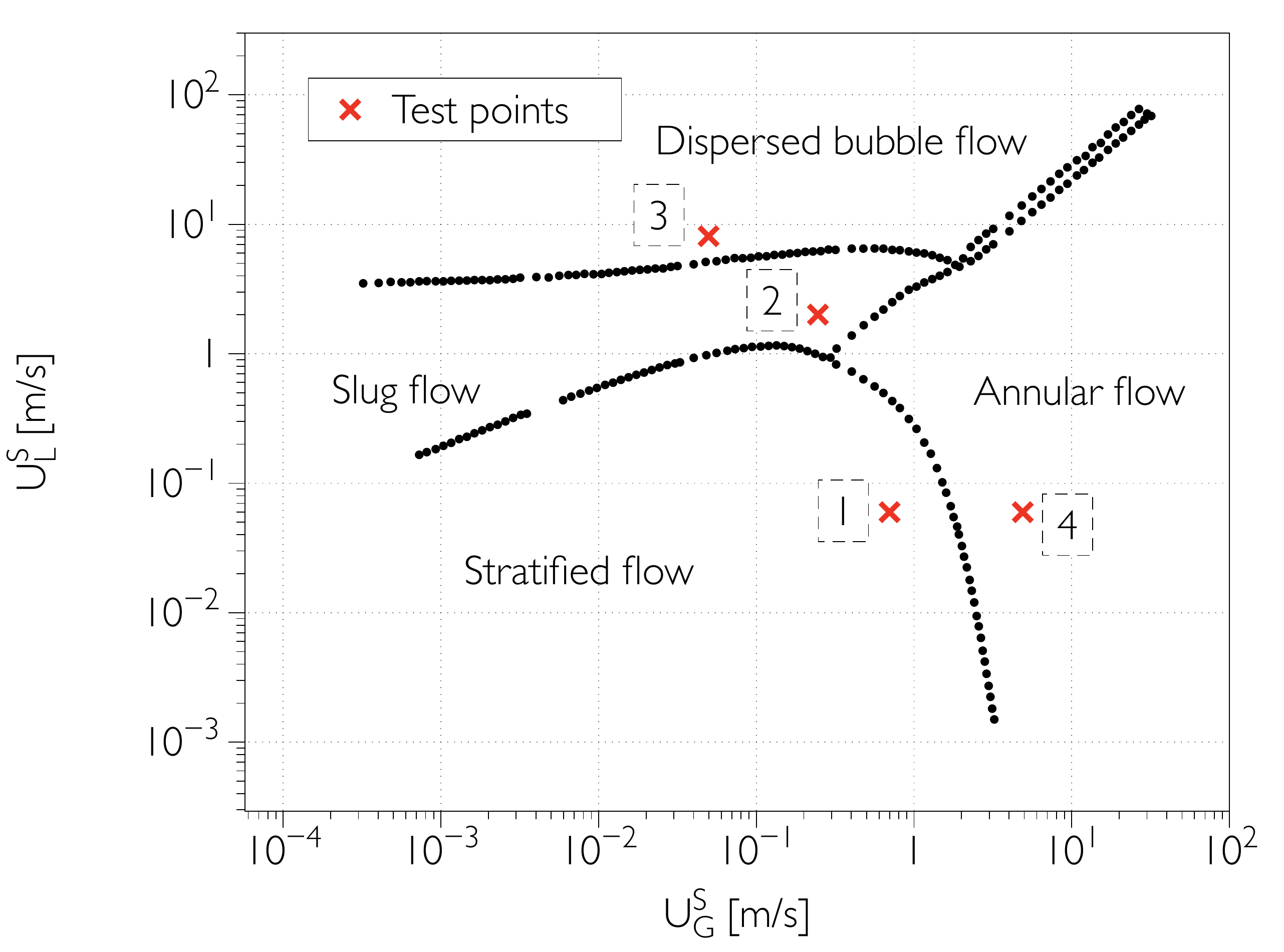}
\caption{Two--phase solver: flow pattern map of the two-phase flow in a horizontal pipe with a diameter of 1 meter. The data has been extracted from \cite{xie2017direct}}
\label{fig-xpipe:TaitelDukler}
\end{figure}

The flow pattern map of Taitel \& Dukler (see Fig.~\ref{fig-xpipe:TaitelDukler}) classifies the flow regimes as stratified flow, slug flow, dispersed bubble flow and annular flow. 
In a stratified flow the phases are completely separated with gas in the upper part and liquid in the lower part of the pipe. 
In a slug flow, the waves in the flow reach the top of the pipe, eventually closing the gas path in the top. 
In a dispersed bubble flow, small bubbles are present in the flow, and are dispersed everywhere in the cross section.
In an annular flow the liquid forms a coat all around the pipe walls. 
Following \cite{xie2017direct}, we define four test cases that reproduce the four flow regimes, denoted with red crosses in Fig.~\ref{fig-xpipe:TaitelDukler}, and whose superficial velocities  are given in Table~\ref{table:pipeRegime}.
%
\begin{table}[h]
\centering
\caption{Two--phase solver: superficial velocities (in m/s) and theoretical flow regime in horizontal pipe test case}
\label{table:pipeRegime}
\begin{tabular}{@{}llllll@{}}
	\toprule
	Test & Flow regime           & $U_G^S $ & $U_L^S$              \\ \midrule
	1         & Stratified flow       & 0.7      & 0.06                 \\ \midrule
	2         & Slug flow             & 0.25     & 2                    \\ \midrule
	3         & Dispersed bubble flow & 0.05     & 8                    \\ \midrule
	4         & Annular flow          & 4.9      & 0.06                 \\ \bottomrule
\end{tabular}
\end{table}
The domain for the simulations consists of a pipe with length $L=1$~m and diameter $D=0.1$~m. The physical domain is discretized using a mesh of 8220 hexahedral elements and the solution is approximated by order $N=3$ polynomials. The physical parameters are summarized in Table~\ref{tab-xpipe:num:pipeflow:params}. 
{\begin{table}[h]
		\centering
		\caption{Two--phase solver: physical parameters of the pipe flow}
		\label{tab-xpipe:num:pipeflow:params}
		\resizebox{\textwidth}{!}{%
		\begin{tabular}{llllllllll}
			\hline
			$\rho_1$& $\rho_2$ ($\text{kg}/\text{m}^3$)   & $\eta_1$ & $\eta_2$ (Pa$\cdot$s) & $\varepsilon$ (m)& $M_0$ (m/s)& $c_0^2$ (m/s$^2$)$^2$& $\sigma$ (N/m)& $g$ (m/s$^2$) \\ \hline
			1.0 & 5.0 & $5\cdot 10^{-3}$ &$ 10^{-2} $& $0.0424$ & 0.1886 & 1.0E3 & $2.5\cdot 10^{-4}$ & 1.0 \\
			\hline
		\end{tabular}}
\end{table}}
The mobility is taken from \cite{xie2017direct}, and the interface width covers 
approximately three high--order mesh points.  The time step chosen for the simulations is 
$\Delta t=10^{-5}$~s. 
Regarding the boundary conditions, a no--slip boundary condition (see \eqref{eq-xpipe:bcs:wall}) is enforced at the pipe walls 
(with a contact angle of $90^\circ$) while a velocity inflow boundary condition (see \eqref{eq-xpipe:bcs:inflow-bc})
and a constant pressure outflow boundary condition (see \eqref{eq-xpipe:bcs:outflow}) are used. It should be noticed that the flow regime inflow is considered layered (see Fig.~\ref{fig-xpipe:inflow-conf}). 
The initial condition for all the simulations is propagated along the $Z$ axis, with a small wave--like perturbation with $Z$ coordinate, to introduce asymmetry.

\begin{figure}[h]
	\centering
	\includegraphics[width=0.95\textwidth]{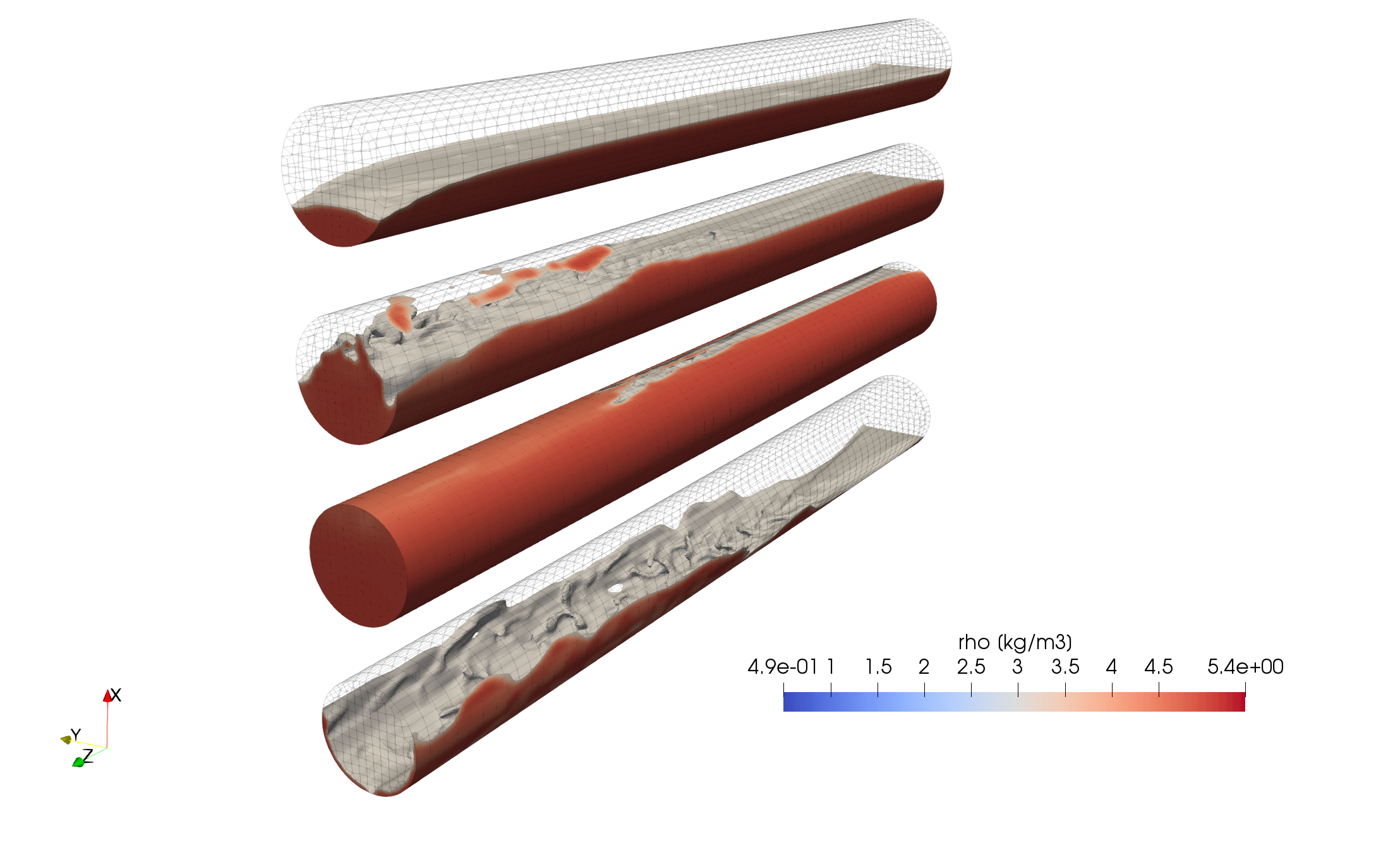}
	\caption{Two--phase solver: results of test cases 1-4 (organized top to bottom) at $t=4$s. Stratified, Slug, Dispersed Bubble and Annular Flow regimes}
	\label{fig-xpipe:PipeRegime3D}
\end{figure}

Fig.~\ref{fig-xpipe:PipeRegime3D}
shows an isosurface of $c$ values under $0.5$ colored by density for the four test cases shown in Fig.~\ref{fig-xpipe:TaitelDukler} (organized from 1 to 4 top to bottom) at $t=4$~s.  

\begin{figure}[h]
	\centering
	\includegraphics[width=0.65\textwidth]{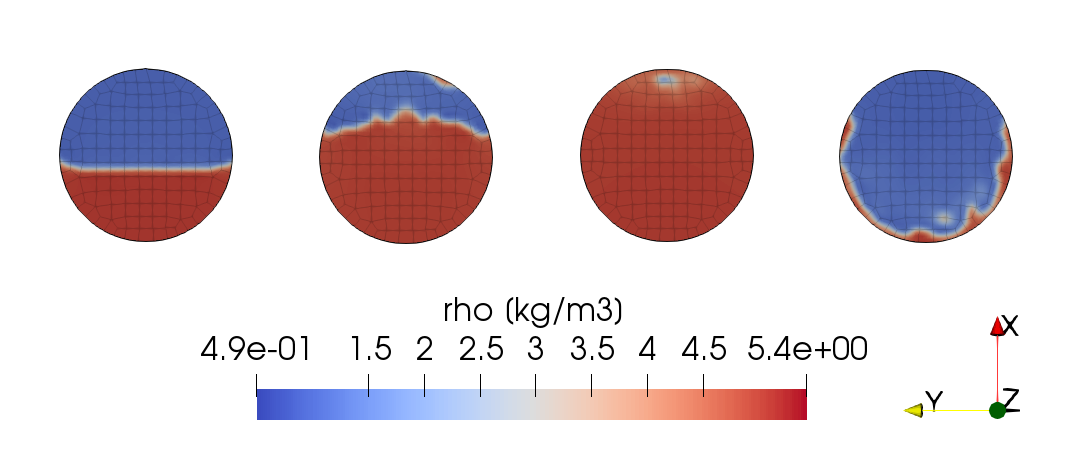}\label{fig-xpipe:sfig1}
	\caption{Two--phase solver: representation of the density contour in a pipe at ${z/D=8}$ cross section. The four regimes (stratified, slug, dispersed bubble, and annular) are represented.}
	\label{fig-xpipe:PipeRegimeSlice}
\end{figure}

In Fig.~\ref{fig-xpipe:PipeRegimeSlice} we show the density contour in a  $Z$--normal slice at $L/D=8$. As can be seen, the flow regimes are correctly predicted. Note that ignoring the effect of the hydrostatic pressure at the outflow when imposing a constant outlet pressure induces a velocity in the negative $X$ direction that curves the interface. This can be more easily seen in the stratified flow regime (Test~1). 

%

\subsection{Enhanced robustness of the discretization}\label{sec-xpipe:num2ph:robustness-split}

The simulations of the flow regime prediction in pipes are under--resolved with the mesh used (except the stratified flow).
Thus, this numerical experiment can be used to compare the robustness of the discretization developed in this work with the standard DG method (i.e. without the use of split--forms)
with the more traditional Gauss points. 
Although none of the two discretizations are entropy--stable (see \cite{2013:Fisher}), we see that the use of the split--form introduced herein also enhances the robustness of the multiphase solver. 

\begin{figure}[h]
	\centering
	\subfigure[Stratified flow]{\includegraphics[width=0.48\textwidth]{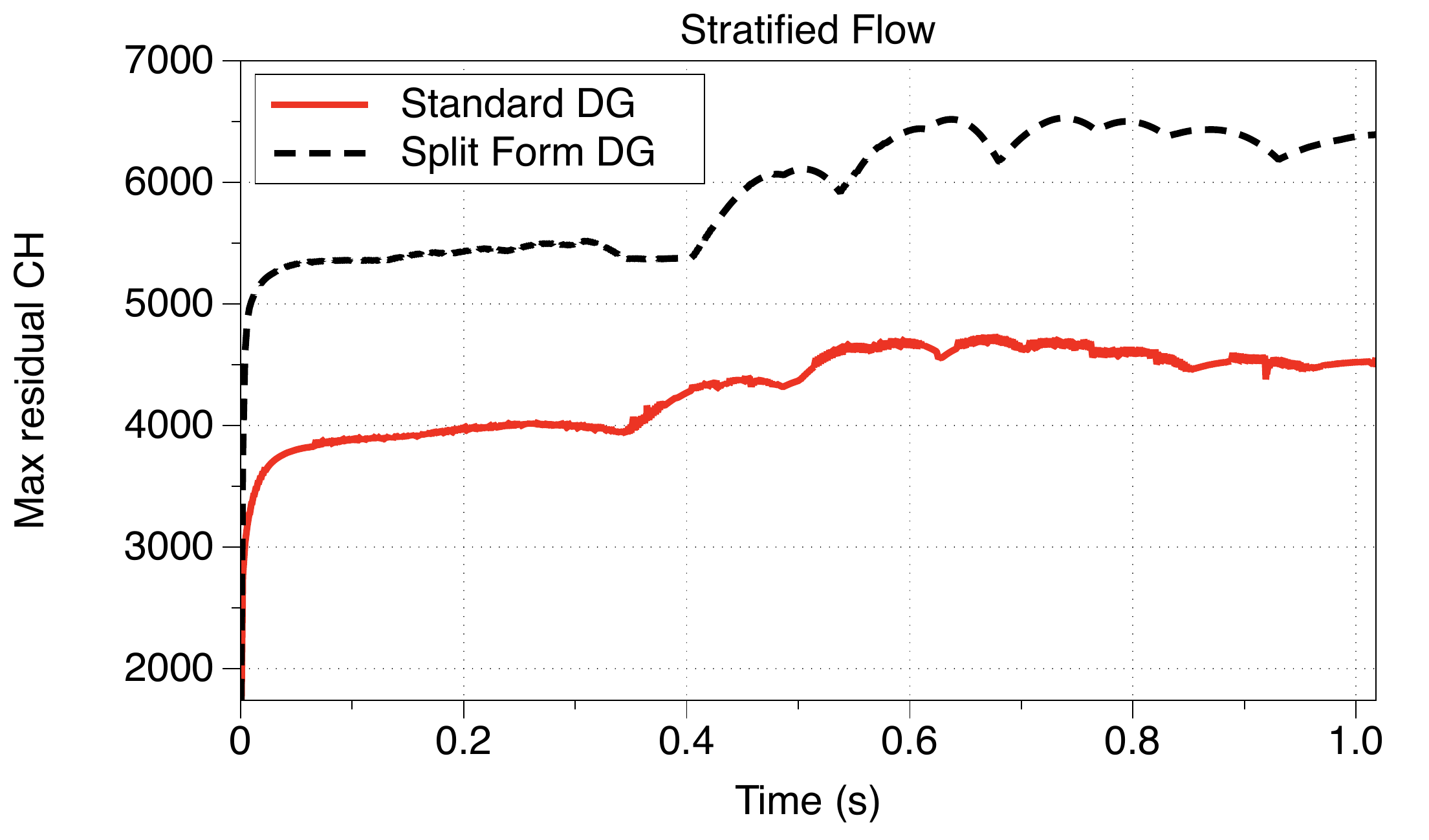}}
	\subfigure[Slug flow]{\includegraphics[width=0.48\textwidth]{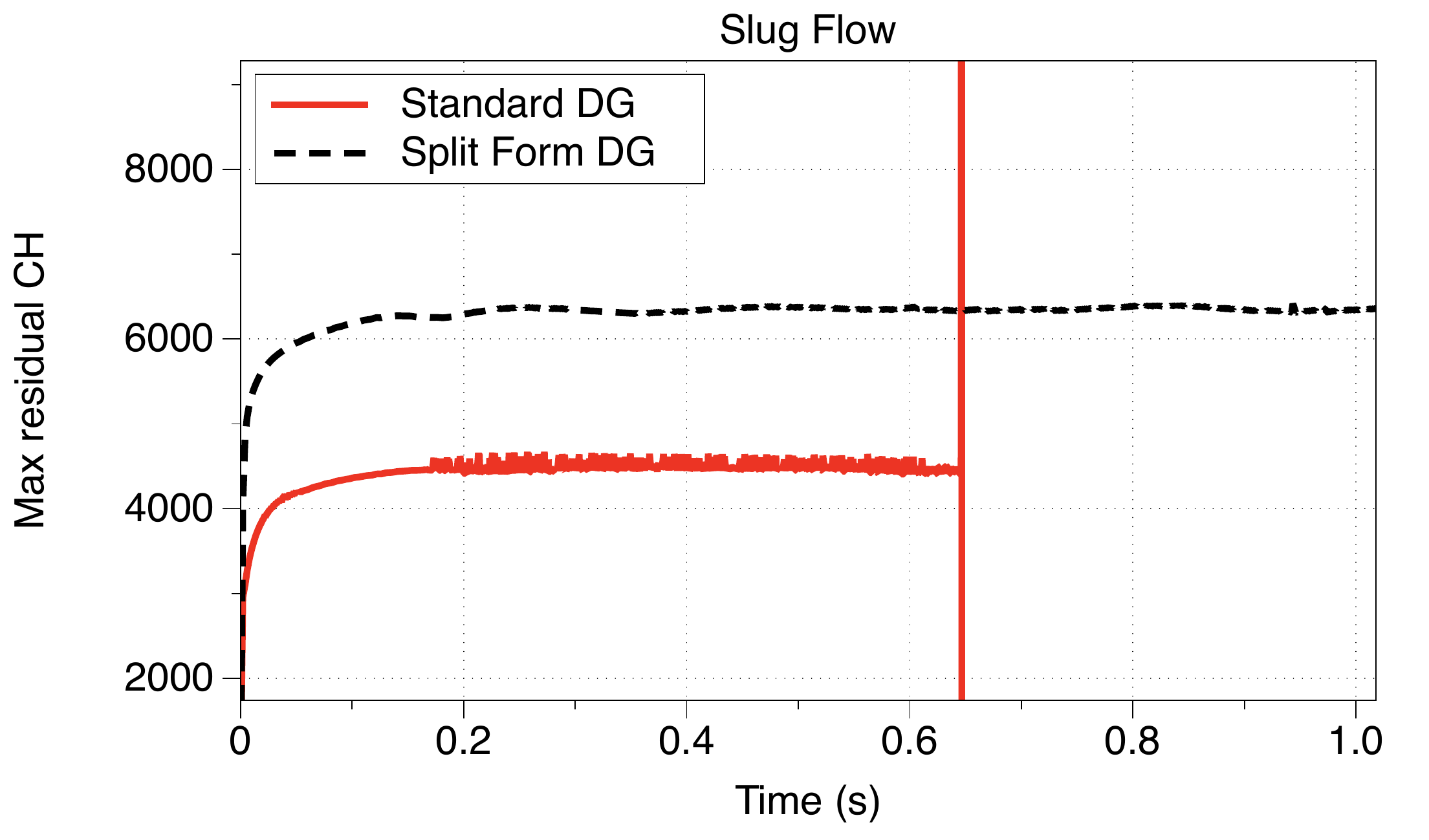}}
	\subfigure[Dispersed bubble flow]{\includegraphics[width=0.48\textwidth]{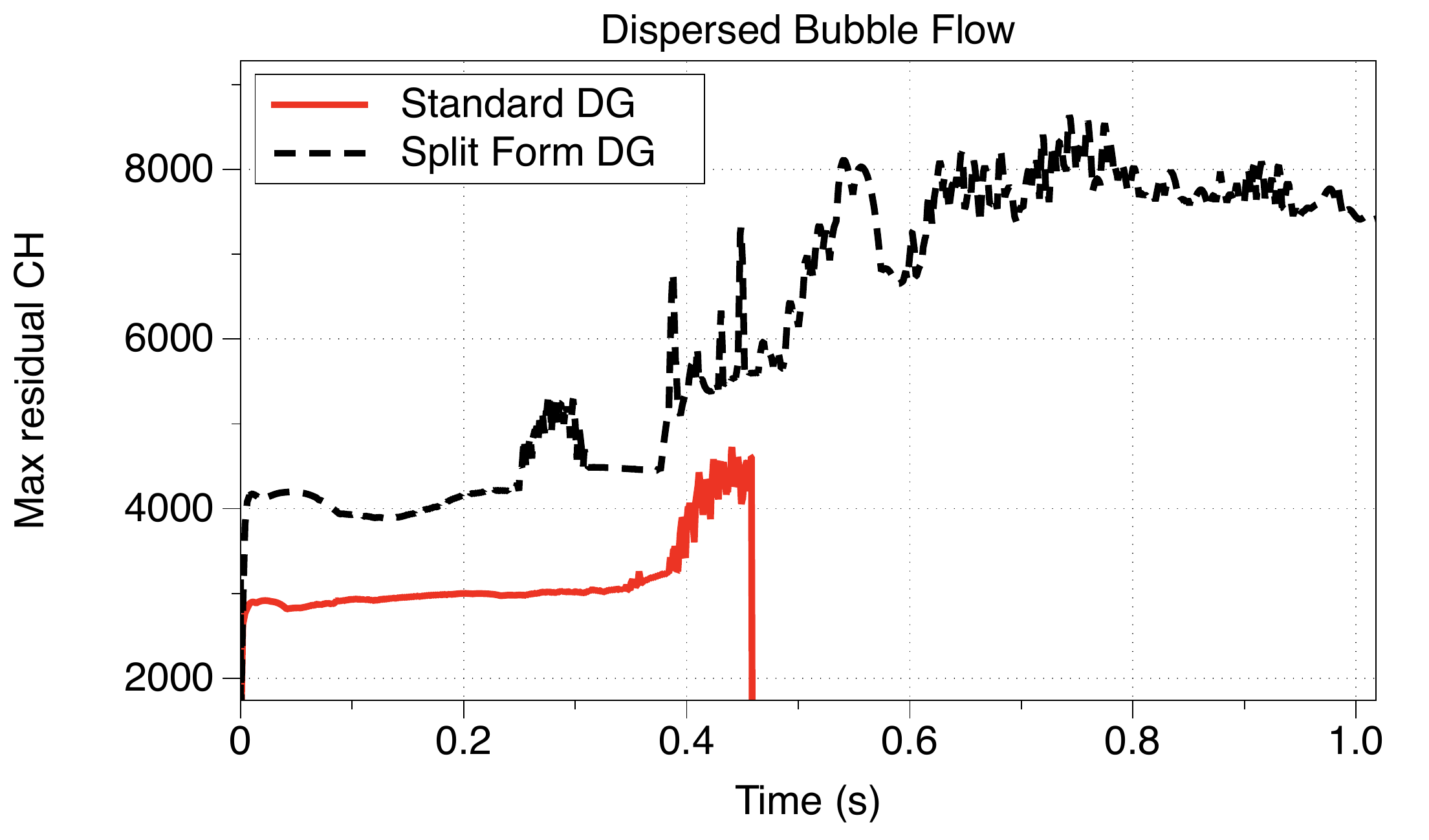}}
	\subfigure[Annular flow]{\includegraphics[width=0.48\textwidth]{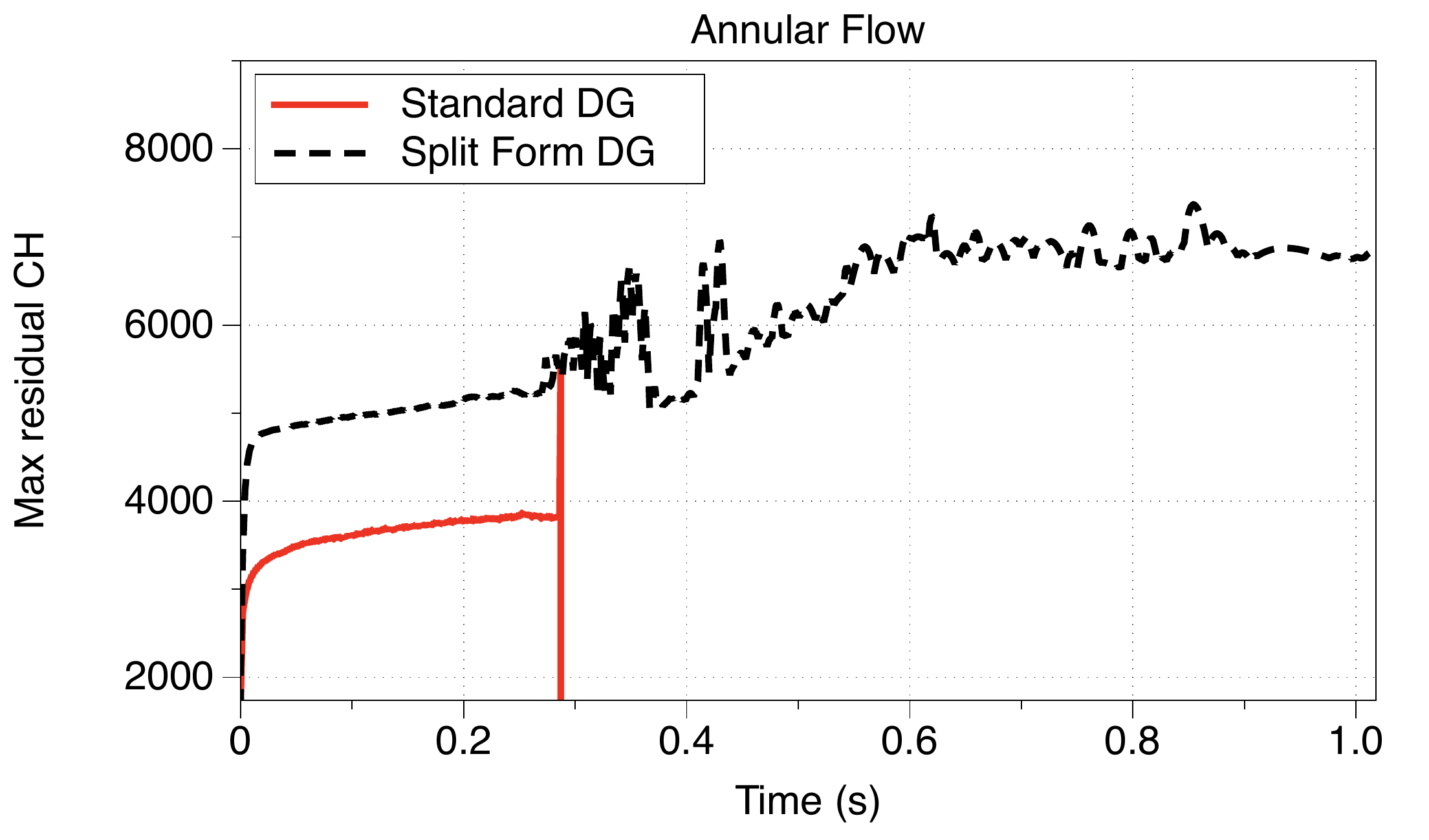}}
	\caption{Two--phase solver: maximum value of the concentration $c_1$ residual $\left(c_{1,t}\right)$ with time, for the four pipe flow regimes studied. 
	We compare the standard DG scheme with Gauss points (solid, red line) 
	to split--form DG (dashed, black)
	}
	\label{fig-xpipe:RobustnessSplitSTD}
\end{figure}

We solve the four regimes with the two schemes (split--form/Gauss--Lobatto and standard/Gauss). In Fig.~\ref{fig-xpipe:RobustnessSplitSTD} 
we represent the maximum residual of the concentration ($\max_{e}C_{t}$) with the number of iterations as a representative value of the simulation progress (as a point of reference, the steady--state is $\max_{e} C_{t}=0$). 
Although both methods are stable at the early stages, once the growth of fluid instabilities obtains smaller structures, the standard DG scheme with Gauss points crashes, while the split--form developed in this work is stable and allows us to obtain a final solution. 
Both approaches are stable for the smoother stratified flow.

We conclude that although there is no stability proof, this numerical experiment encourages the use of the split--form scheme.

\section{Three--phase simulations}\label{sec-xpipe:Validation3PH}

In this section we validate and test the three--phase solver. We solve a 
manufactured solution in Sec.~\ref{sec-xpipe:num3ph:mansol}, a 
two--dimensional horizontal channel in Sec.~\ref{sec-xpipe:num3ph:2DChannel}
and a T--shaped 
pipe intersection in Sec.~\ref{sec-xpipe:num3ph:Tee}.

\subsection{Manufactured solution}\label{sec-xpipe:num3ph:mansol}

First, we perform a manufactured solution analysis to assess the accuracy of the three--phase solver for smooth problems. We extend the two--dimensional manufactured solution used for the two--phase flow in Sec.~\ref{sec-xpipe:num2ph:Mansol} 
with the manufactured solution of the three--phase Cahn--Hilliard solver 
introduced in \cite{2020:Manzanero-UR-CaF}.

The manufactured solution now is defined as:
\begin{equation}
\begin{split}
c_{1,0}(x,y;t) &= \frac{1}{3}\left(1 + \cos\left(\pi x\right) \sin\left(\pi y\right)\sin\left(t\right) \right), \\
c_{2,0}\left(x,y;t\right) &=\frac{1}{3}\left(1 + \cos\left(\pi x\right) \sin\left(\pi y\right)\sin\left(1.2 t\right) \right),\\
u_0(x,y;t) &=2\sin\left(\pi x\right) \cos\left(\pi z\right)\sin\left(t\right), \\
v_0(x,y;t) &=-2\cos\left(\pi x\right)\sin\left(\pi y\right) \sin\left(t\right), \\
p_0(x,y;t) &=2\sin\left(\pi x\right)\sin\left(\pi z\right)\cos\left(t\right),
\end{split}
\label{eq-xpipe:num3ph:conv:mansol}
\end{equation}
which requires an appropriate source term to the right--hand side of the equation, 
not presented here for simplicity. 
The configuration is similar to the two--phase manufactured solution, 
where the domain is $\left(x,y\right)\in[-1,1]^2$~m, and the final time is $t_F=0.1$~s. 
The physical parameters have been adapted from \cite{dong2018multiphase} and they are given in 
Table~\ref{tab-xpipe:num3ph:mansol:param}.
{\begin{table}[h]
		\centering
		\caption{Three--phase solver: list of the parameter values used with the manufactured solution (see \eqref{eq-xpipe:num3ph:conv:mansol})}
		\label{tab-xpipe:num3ph:mansol:param}
			\begin{tabular}{lllllll}
				\hline
				$\rho_1$& $\rho_2$ & $\rho_3$ ($\text{kg}/\text{m}^3$)   & $\eta_1$ & $\eta_2$ & $\eta_3$  (Pa$\cdot$s) & $\varepsilon$ (m) \\ \hline
				1.0 & 3.0 & 2.0 & 1.0E-3 & 1.0E-3 & 1.0E-3 & $1/\sqrt{2}$   \\
				\hline \\ \hline 
				& $M_0$ (m/s) & $c_0$ (m/s$^2$) & $\sigma_{12}$& $\sigma_{13}$ & $\sigma_{23}$ (N/m)  \\ \hline
    & 1.134E-2 & 1.0E3 &  6.236E-3 & 7.265E-3 & 8.165E-3 \\
				\hline				
		\end{tabular}
	\end{table}}

\begin{figure}[h]
  \centering
  \subfigure[$\Delta t=10^{-4}$~s]{\includegraphics[width=0.49\textwidth]{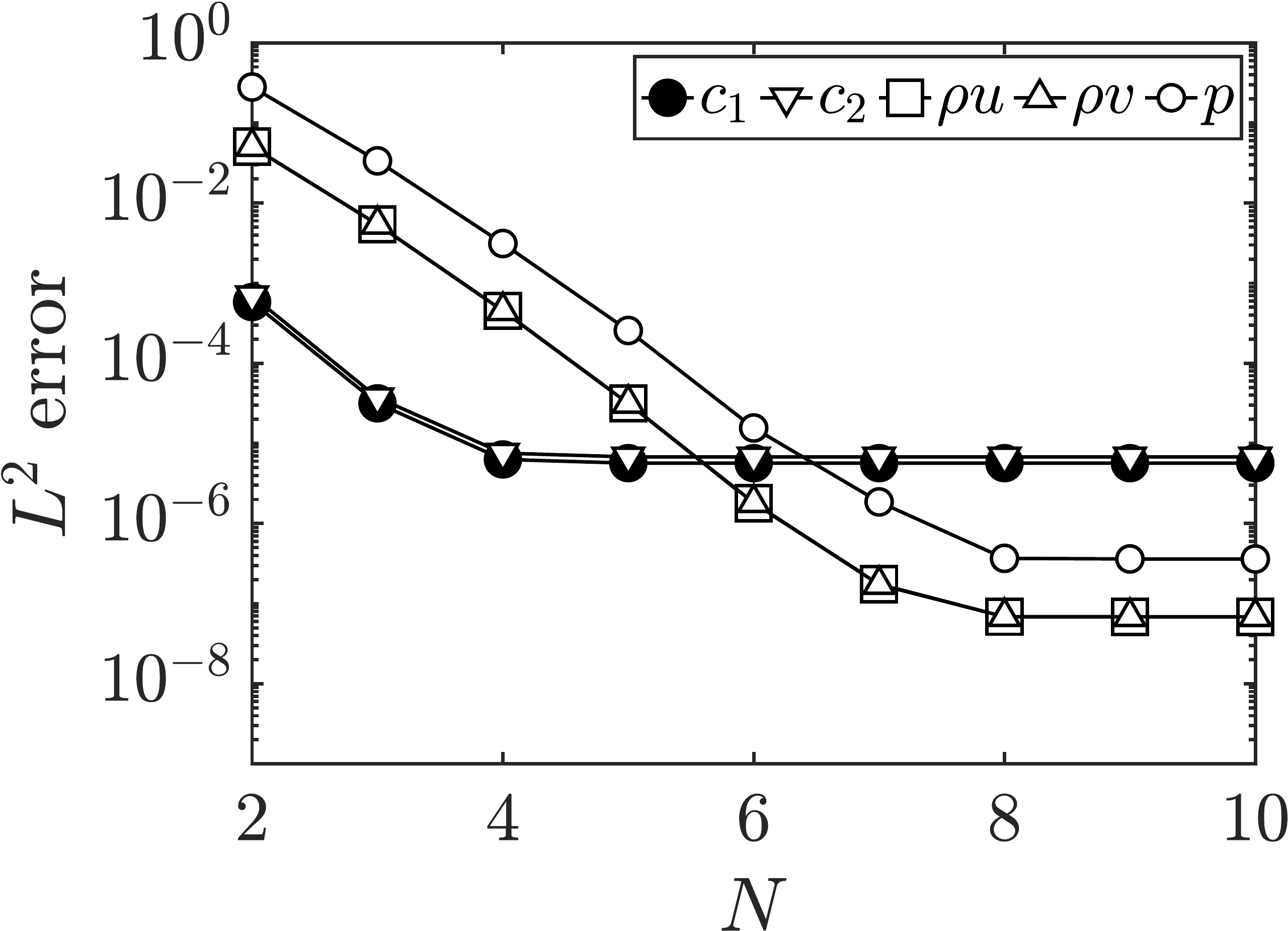}}
  \subfigure[$\Delta t=10^{-5}$~s]{\includegraphics[width=0.49\textwidth]{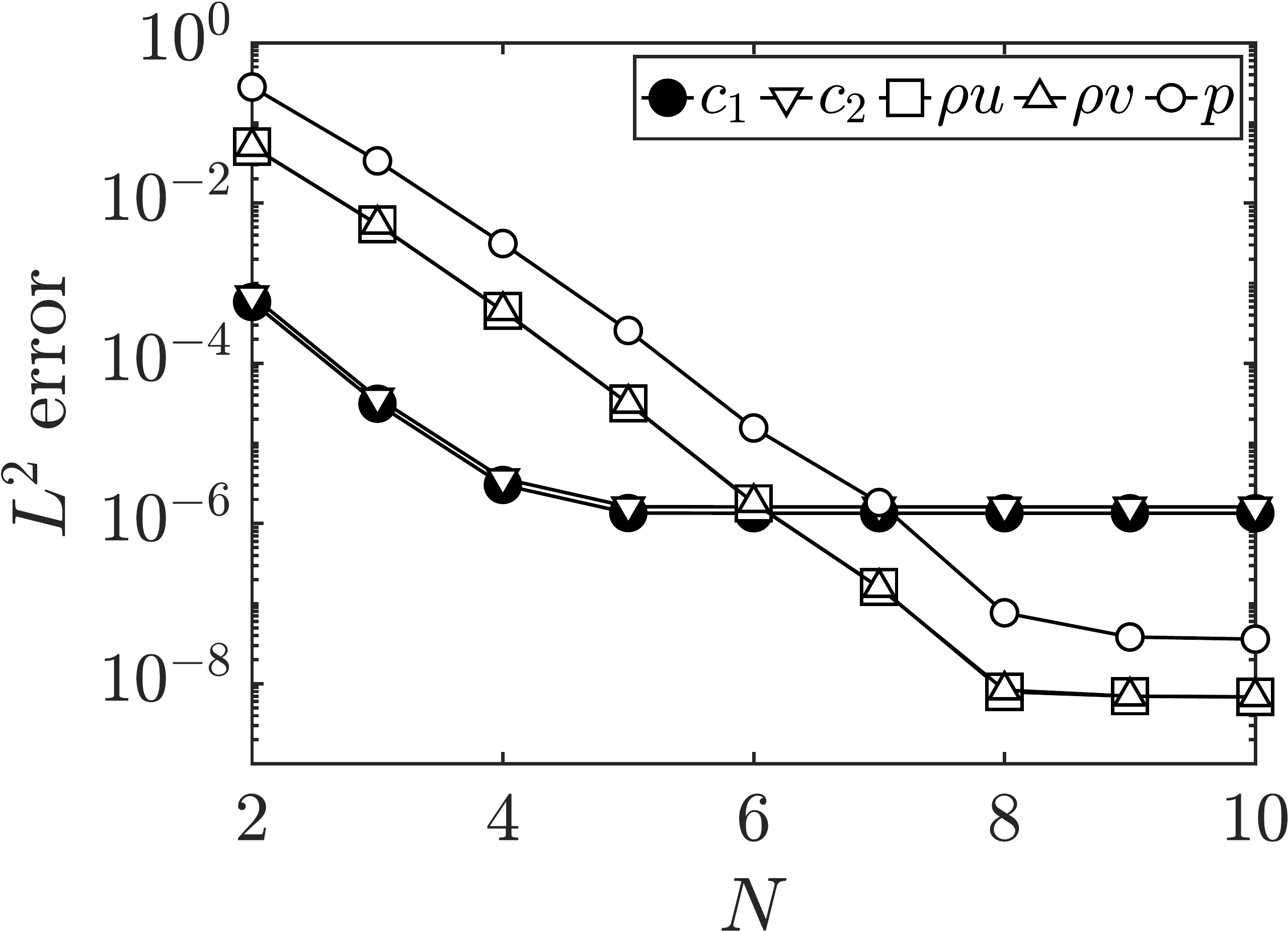}}
  \caption{Three--phase solver: polynomial order convergence study of the manufactured solution \eqref{eq-xpipe:num3ph:conv:mansol}. 
  We represent the L$^2$ errors of the two concentrations $c_1$ and $c_2$, $x$-- and $y$--momentum, and pressure. 
  The polynomial order ranges from 2 to 10, and we integrate in time until $t_{F}=0.1$~s with two time step sizes: $\Delta t=10^{-4}$~s and $10^{-5}$~s. All physical parameters are given in Table \ref{tab-xpipe:num3ph:mansol:param}}
  \label{fig-xpipe:num3ph:mansol}
\end{figure}

We perform first a polynomial order convergence study with a Cartesian mesh of $4^2$ elements, and with the polynomial order ranging from $N=2$ to 10. In Fig.~\ref{fig-xpipe:num3ph:mansol} we represent the L$^2$ errors on the five variables ($c_1$, $c_2$, $\rho u$, $\rho v$ and $p$), for two time--step sizes, $\Delta t = 10^{-4}$~s and $10^{-5}$~s. We find that the error behavior is similar to the two--phase solver, where exponential accuracy is obtained for lower polynomial orders, as expected, and then the error stagnation associated to the time discretization is anticipated for the two concentrations, as a result of the first order IMEX scheme.

Finally, we perform the mesh convergence study, where we use meshes with $4^2$, $6^2$, $8^2$, $12^2$ and $16^2$ elements, and vary the polynomial order from $N=2$ to $N=5$. The L$^2$ errors and the convergence rates are written in Table~\ref{tab-xpipe:num3ph:conv:h-conv}.
\begin{table}[h]
\centering
  \caption{Three--phase solver: manufactured solution \eqref{eq-xpipe:num3ph:conv:mansol} convergence analysis: we use $4^3$, $8^3$, and $16^3$  meshes, and $N=2,3,4$ and 5.
The final time is $t_{F}=0.1$~s, and we use the IMEX scheme
  with $\Delta t=5\cdot 10^{-5}$~s}
  \label{tab-xpipe:num3ph:conv:h-conv}
\resizebox{\textwidth}{!}{%
\begin{tabular}{llllllllllll}
\hline
      & Mesh   & $c_1$ error & order & $c_2$ error & order& ${\rho} u$ error & order & ${\rho} v$ error & order & $p$ error & order \\ \hline
N=2 & $4^2$ & 5.82E-04 & -- & 6.97E-04 & -- & 5.06E-02 & -- & 5.06E-02 & -- & 2.81E-01 & -- \\ 
 & $6^2$ & 1.60E-04 & 3.19 & 1.91E-04 & 3.19 & 1.92E-02 & 2.39 & 1.92E-02 & 2.39 & 1.17E-01 & 2.15 \\ 
 & $8^2$ & 5.97E-05 & 3.42 & 7.15E-05 & 3.41 & 9.38E-03 & 2.49 & 9.38E-03 & 2.49 & 6.19E-02 & 2.22 \\ 
 & $12^2$ & 1.61E-05 & 3.23 & 1.93E-05 & 3.23 & 3.28E-03 & 2.59 & 3.28E-03 & 2.59 & 2.40E-02 & 2.34 \\ 
 & $16^2$ & 6.65E-06 & 3.07 & 7.97E-06 & 3.07 & 1.51E-03 & 2.70 & 1.51E-03 & 2.70 & 1.18E-02 & 2.47 \\ 
N=3 & $4^2$ & 3.13E-05 & -- & 3.69E-05 & -- & 5.41E-03 & -- & 5.41E-03 & -- & 3.36E-02 & -- \\ 
 & $6^2$ & 3.67E-06 & 5.29 & 4.38E-06 & 5.26 & 1.30E-03 & 3.52 & 1.30E-03 & 3.52 & 9.41E-03 & 3.14 \\ 
 & $8^2$ & 1.49E-06 & 3.13 & 1.79E-06 & 3.11 & 4.51E-04 & 3.68 & 4.50E-04 & 3.68 & 3.65E-03 & 3.30 \\ 
 & $12^2$ & 1.35E-06 & 0.25 & 1.62E-06 & 0.25 & 9.55E-05 & 3.83 & 9.55E-05 & 3.83 & 9.01E-04 & 3.45 \\ 
 & $16^2$ & 1.35E-06 & 0.00 & 1.61E-06 & 0.00 & 3.08E-05 & 3.94 & 3.08E-05 & 3.94 & 3.20E-04 & 3.60 \\ 
N=4 & $4^2$ & 3.03E-06 & -- & 3.65E-06 & -- & 4.52E-04 & -- & 4.52E-04 & -- & 3.12E-03 & -- \\ 
 & $6^2$ & 1.39E-06 & 1.93 & 1.66E-06 & 1.94 & 7.10E-05 & 4.56 & 7.10E-05 & 4.57 & 5.47E-04 & 4.30 \\ 
 & $8^2$ & 1.35E-06 & 0.10 & 1.62E-06 & 0.10 & 1.85E-05 & 4.67 & 1.85E-05 & 4.67 & 1.52E-04 & 4.46 \\ 
 & $12^2$ & 1.35E-06 & 0.00 & 1.61E-06 & 0.00 & 2.67E-06 & 4.77 & 2.67E-06 & 4.77 & 2.37E-05 & 4.58 \\ 
 & $16^2$ & 1.35E-06 & 0.00 & 1.61E-06 & 0.00 & 6.59E-07 & 4.87 & 6.59E-07 & 4.87 & 6.20E-06 & 4.67 \\ 
N=5 & $4^2$ & 1.35E-06 & -- & 1.62E-06 & -- & 3.17E-05 & -- & 3.15E-05 & -- & 2.57E-04 & -- \\ 
 & $6^2$ & 1.35E-06 & 0.01 & 1.61E-06 & 0.01 & 3.18E-06 & 5.67 & 3.18E-06 & 5.66 & 2.83E-05 & 5.44 \\ 
 & $8^2$ & 1.35E-06 & 0.00 & 1.61E-06 & 0.00 & 6.23E-07 & 5.67 & 6.22E-07 & 5.67 & 5.79E-06 & 5.52 \\ 
 & $12^2$ & 1.35E-06 & 0.00 & 1.61E-06 & 0.00 & 6.48E-08 & 5.58 & 6.47E-08 & 5.58 & 5.81E-07 & 5.67 \\ 
 & $16^2$ & 1.35E-06 & 0.00 & 1.61E-06 & 0.00 & 1.50E-08 & 5.08 & 1.50E-08 & 5.08 & 1.15E-07 & 5.62 \\ 
 \\
\hline
\end{tabular}}
\end{table}
We observe that for the two concentrations, the convergence rates are always between $N$ and $N+2$, for $N=2$ and $N=3$, and due to the early stagnation, as in the two--phase solver, we cannot evaluate the convergence rates for $N=4$ and $N=5$. For the rest of the variables, as in the two--phase solver, we find that the convergence rates are always between $N$ and $N+1$, as expected.

Overall, we confirm that the scheme and its implementation are accurate for the industrial applications.

\subsection{Two--dimensional channel}\label{sec-xpipe:num3ph:2DChannel}

In this section, we study the three--phase flow obtained in a two--dimensional channel $\Omega=[0,10]{\times}[-0.5,0.5]$~m. The configuration is an extension to three--phase of the two--phase channel studied in \cite{xie2017direct}. At the inlet, we introduce a heavy fluid (Phase~2, red) on the top, and a light fluid (Phase~3, black) on the bottom, both immersed in Phase~1 (white), with medium density. The vertical gravity results in the heavy fluid falling to the bottom side of the domain, while the light fluid rising to the top. The physical parameters are given in Table~\ref{tab-xpipe:num3ph:2Dpipe:param}.
\begin{table}[h]
		\centering
		\caption{Three--phase solver: list of the parameter values used for the two--dimensional three--phase channel}
		\label{tab-xpipe:num3ph:2Dpipe:param}
		\begin{tabular}{lllllll}
			\hline
			$\rho_1$& $\rho_2$ & $\rho_3$ ($\text{kg}/\text{m}^3$)   & $\eta_1$ & $\eta_2$ & $\eta_3$  (Pa$\cdot$s) & $\varepsilon$ (m) \\ \hline
			1.0 & 5.0 & 0.8 & 5.0E-3 & 1.0E-2 & 1.0E-2 & $0.0424$   \\
			\hline \\ \hline
			&$M_0$ (m/s) & $c_0$ (m/s$^2$) & $\sigma_{12}$& $\sigma_{13}$ & $\sigma_{23}$ (N/m)  \\ \hline
			&1.0E-4 & 1.0E3 &  2.5E-4 & 2.5E-4 & 2.5E-4 \\
			\hline				
		\end{tabular}
\end{table}
At the inlet, we impose the inflow boundary condition, and the initial condition is $c_1=1$, $c_2=c_3=0$, $u=1-4z^2$, and $v=w=p=0$. The computational domain $\Omega$ is divided into 15$\times$60 elements (which is the coarsest configuration described in \cite{xie2017direct}), and we approximate the solution with $N=5$ polynomials. We use the split--form scheme and the IMEX time integrator with $S_0=8$, and $\Delta t = 3.0\cdot 10^{-5}$~s.
\begin{equation}
\begin{split}
c_{1} &= 1-c_2-c_3, \\
c_2 &= \frac{1}{2}\left(1+\tanh\left(\frac{z-0.3}{\varepsilon}\right)\right),\\
c_3 &= \frac{1}{2}\left(1-\tanh\left(\frac{z+0.3}{\varepsilon}\right)\right),\\
u &= 1-4z^2,\\
v &= w = 0,\\
p &= p_{in}.
\end{split}
\end{equation}

\begin{figure}
	\centering
	\includegraphics[clip,width=\textwidth,trim=0cm 11cm 1cm 10cm]{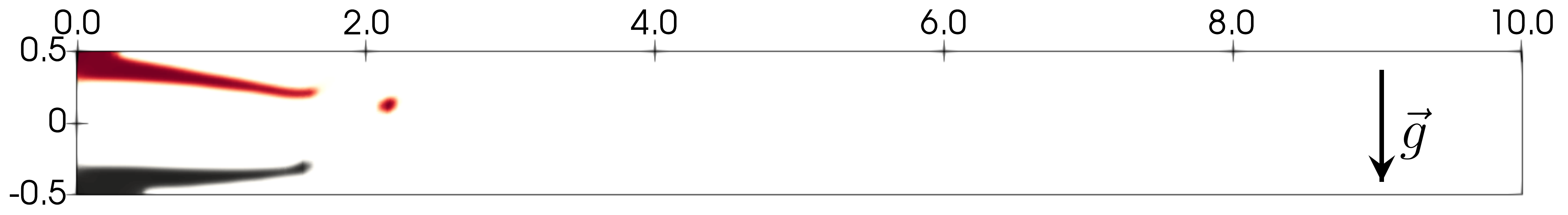}
	\includegraphics[clip,width=\textwidth,trim=0cm 11cm 1cm 11cm]{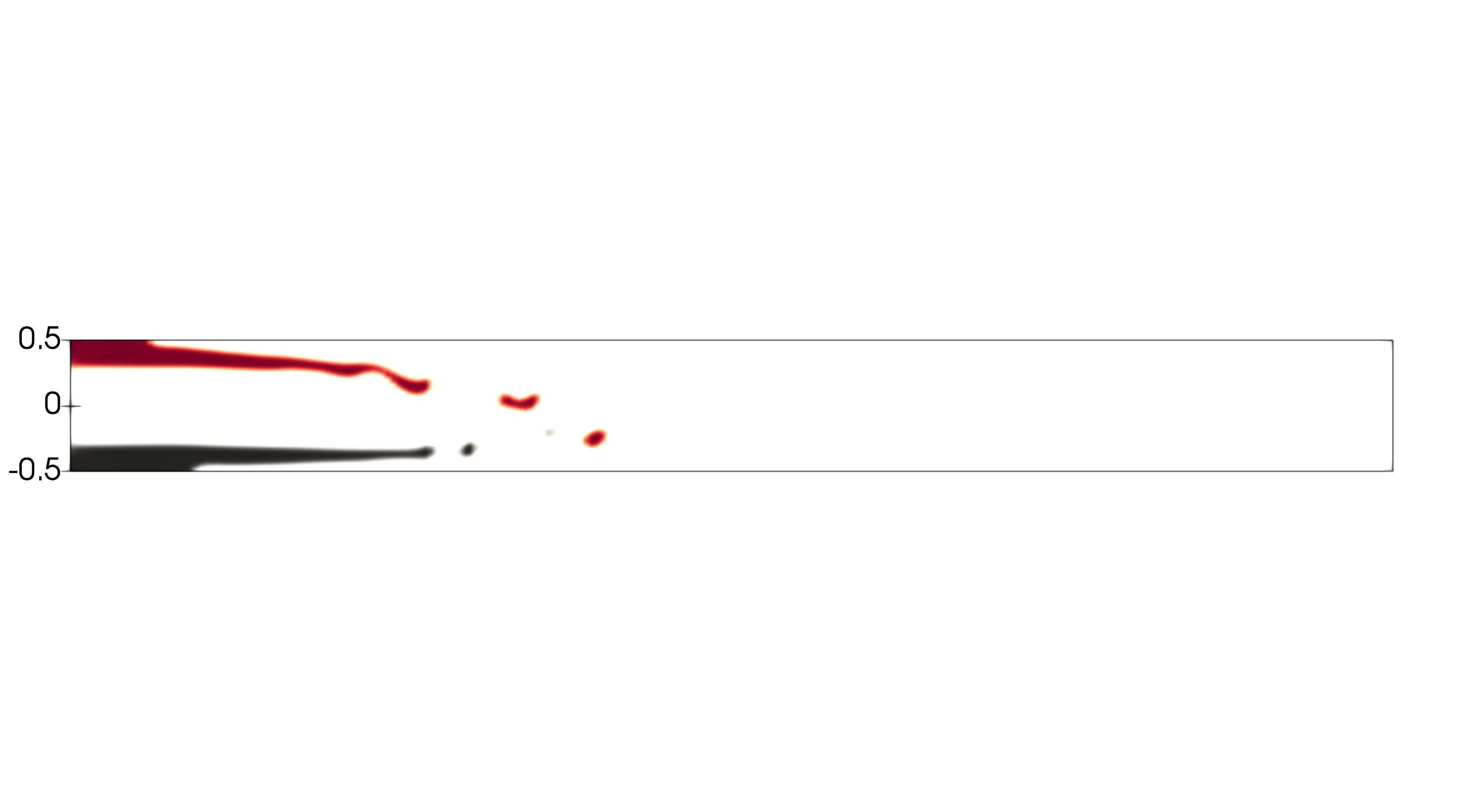}
	\includegraphics[clip,width=\textwidth,trim=0cm 11cm 1cm 11cm]{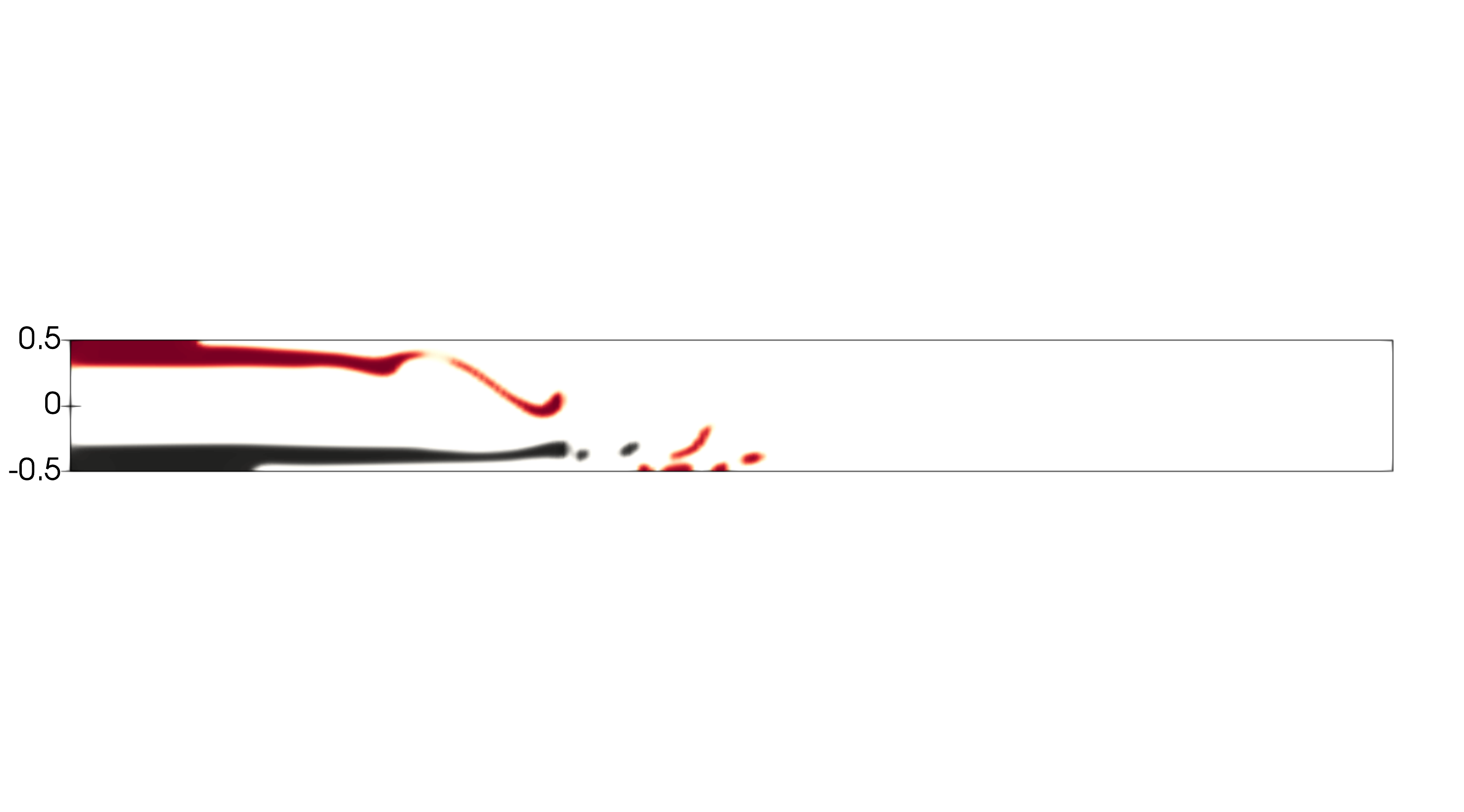}
	\includegraphics[clip,width=\textwidth,trim=0cm 11cm 1cm 11cm]{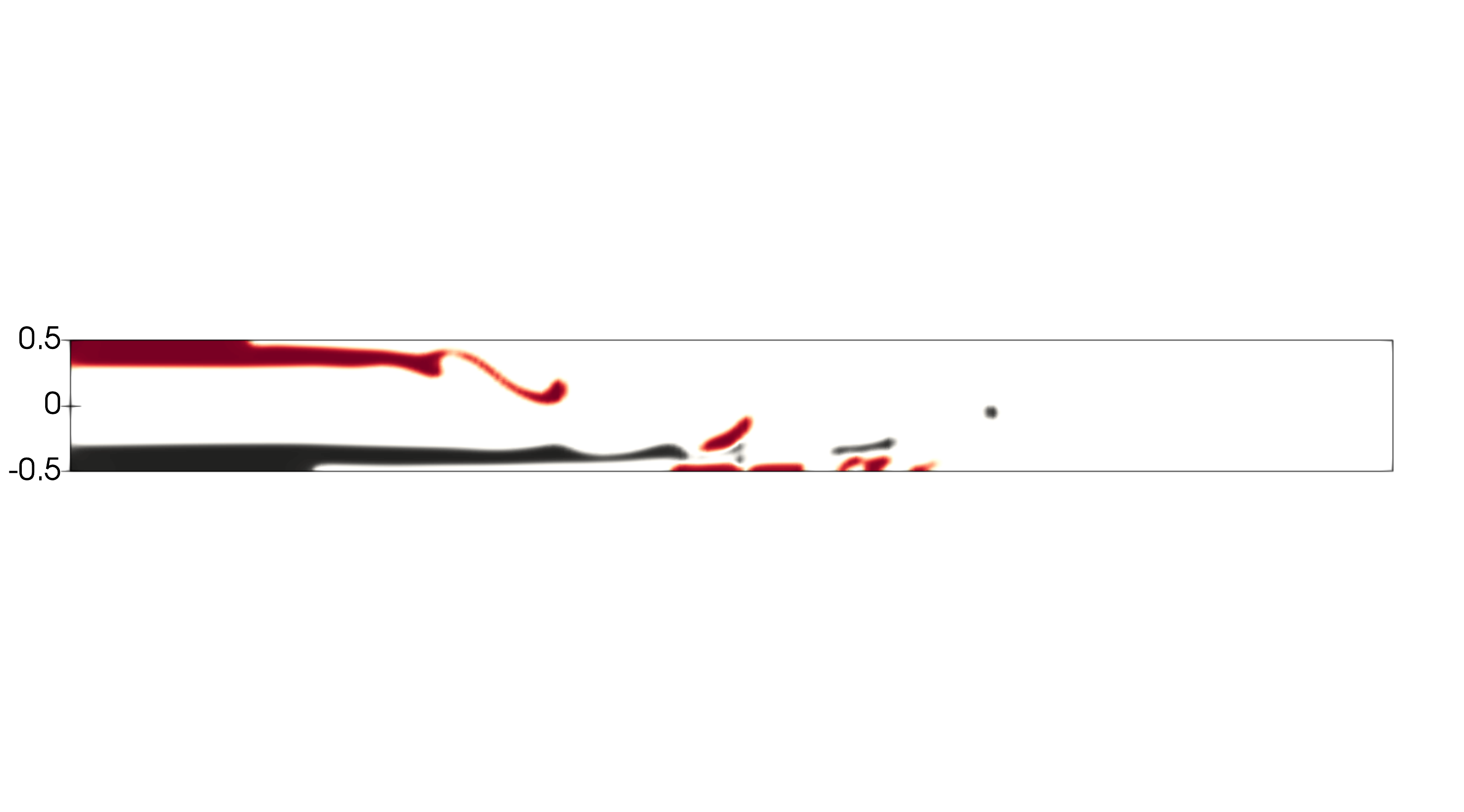}
	\includegraphics[clip,width=\textwidth,trim=0cm 11cm 1cm 11cm]{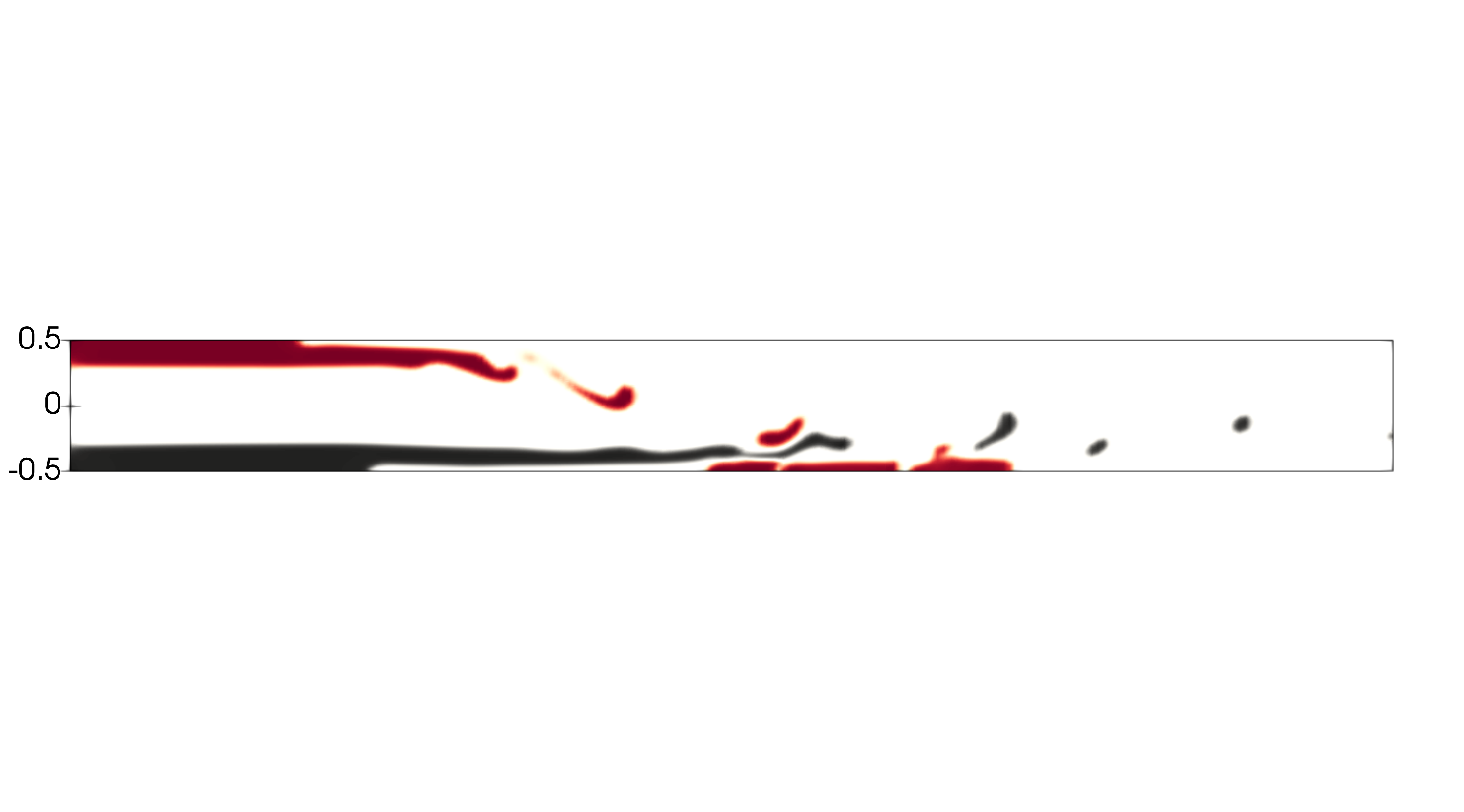}
	\includegraphics[clip,width=\textwidth,trim=0cm 11cm 1cm 11cm]{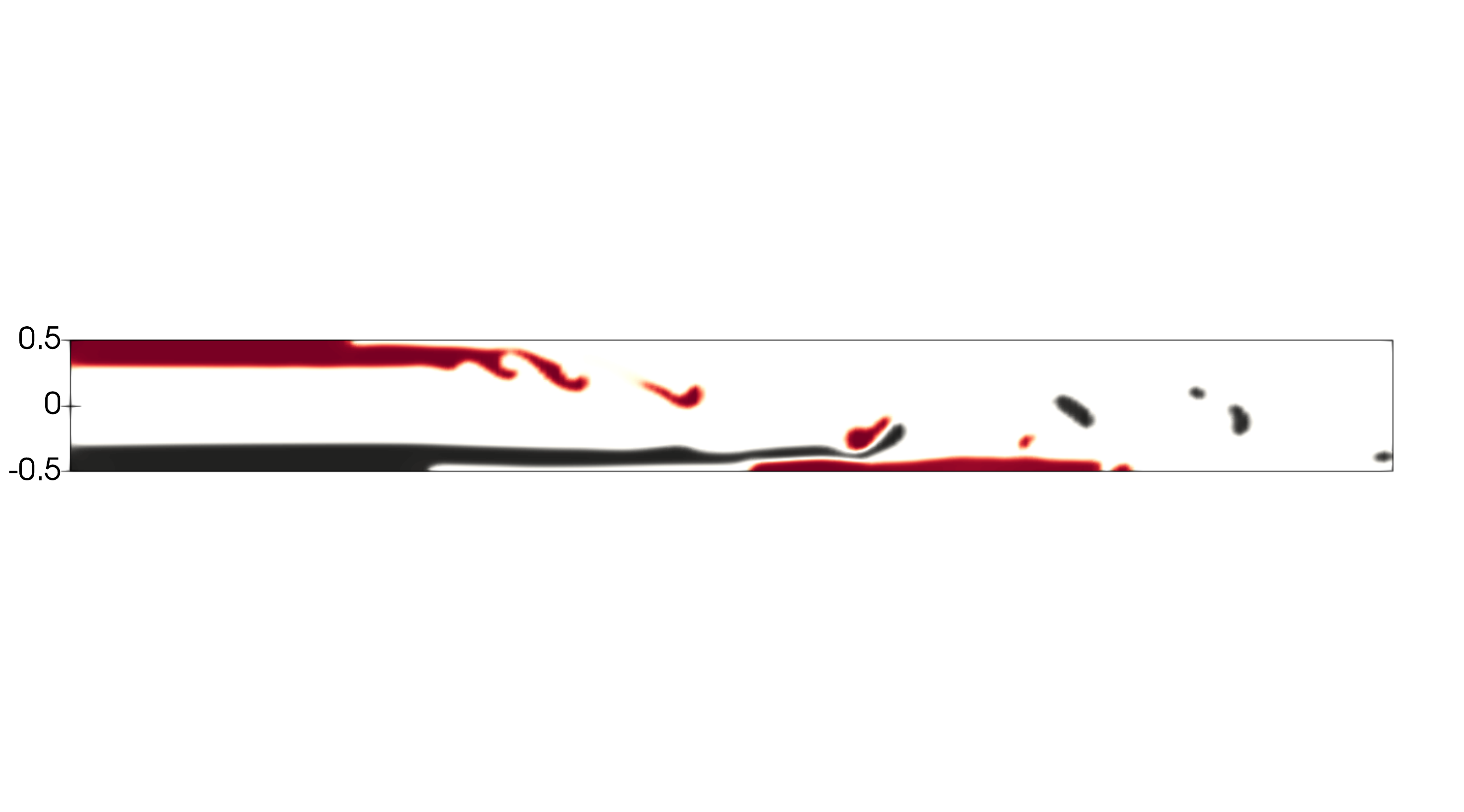}
	\includegraphics[clip,width=\textwidth,trim=0cm 11cm 1cm 11cm]{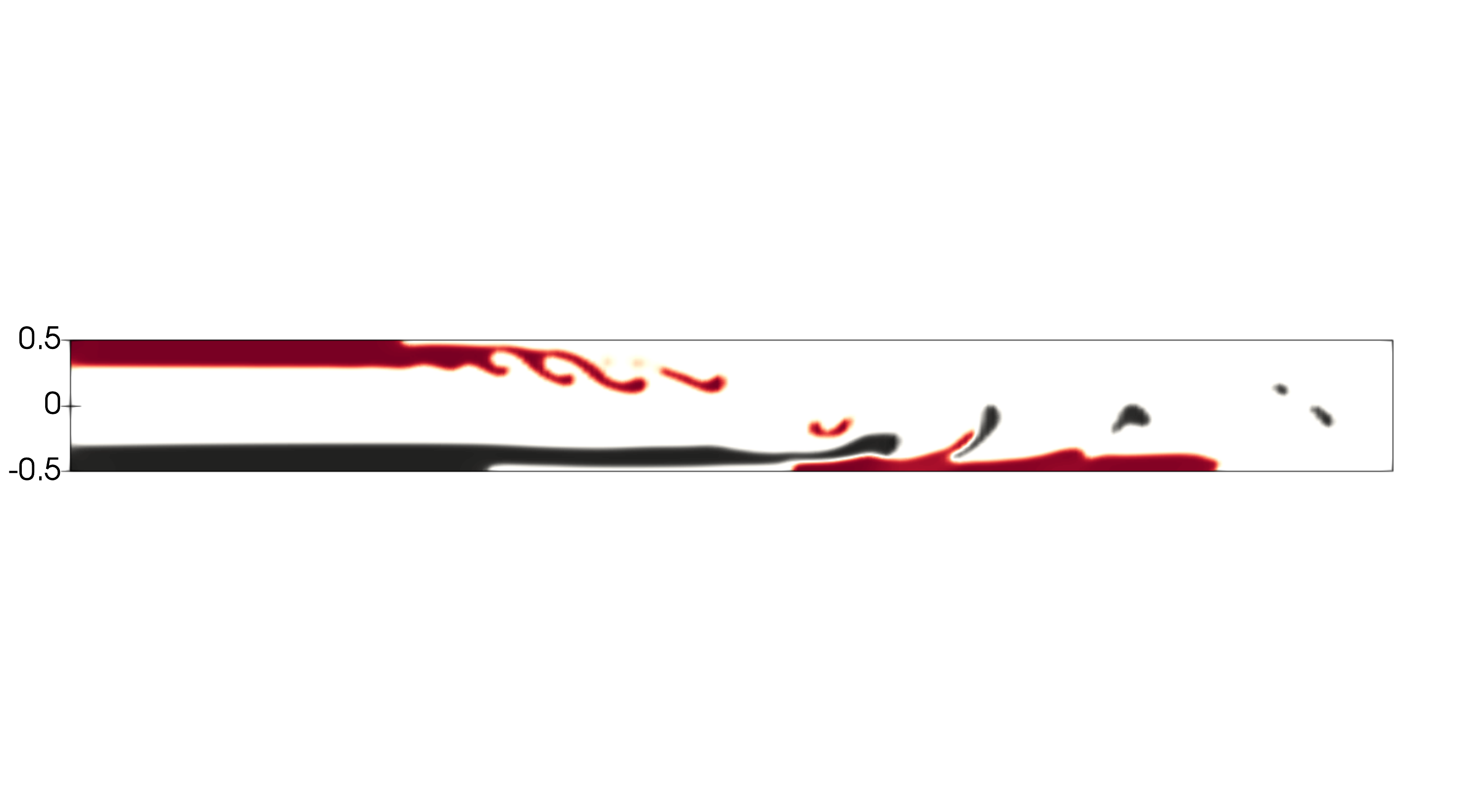}
	\includegraphics[clip,width=\textwidth,trim=0cm 11cm 1cm 11cm]{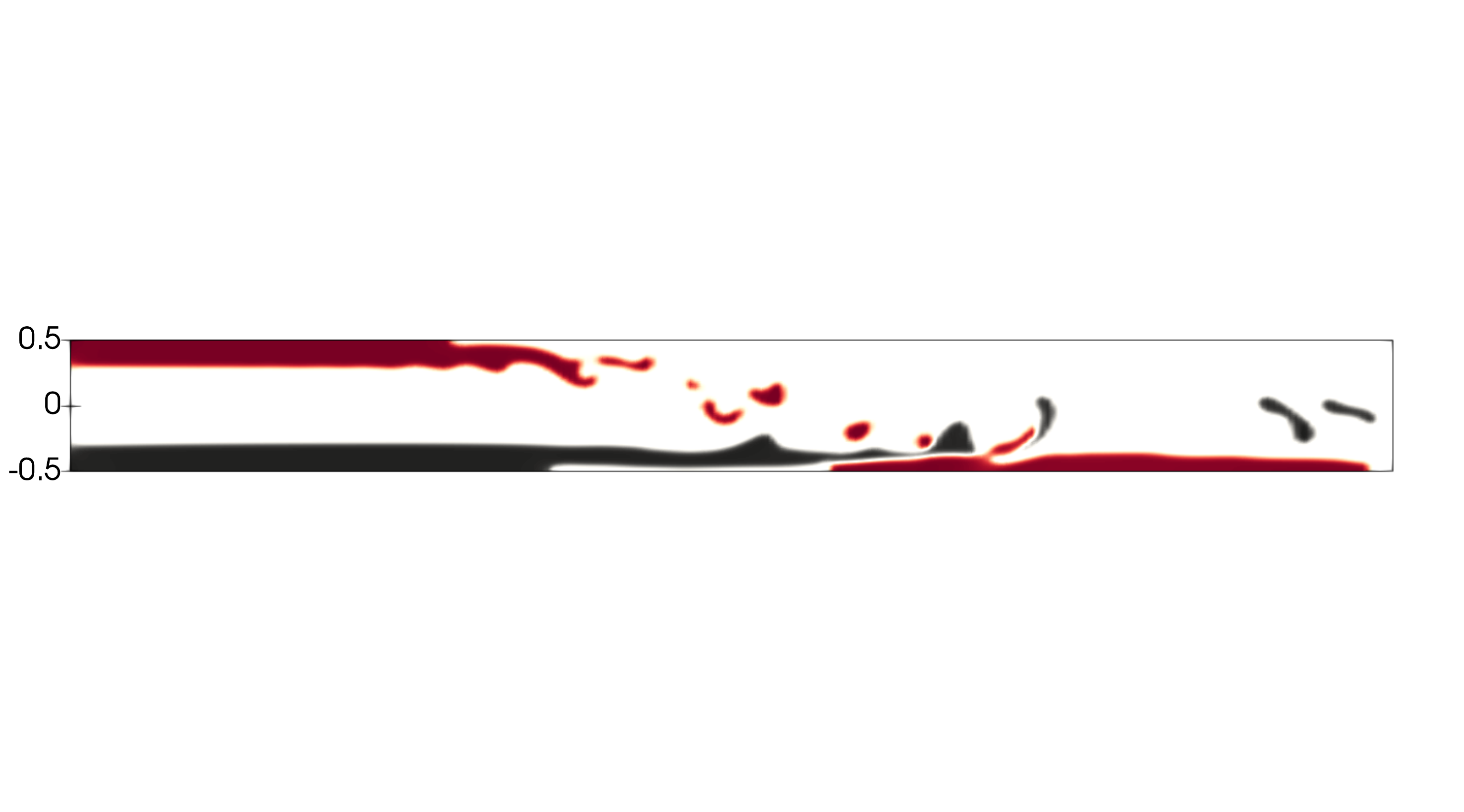}
	\includegraphics[clip,width=\textwidth,trim=0cm 11cm 1cm 11cm]{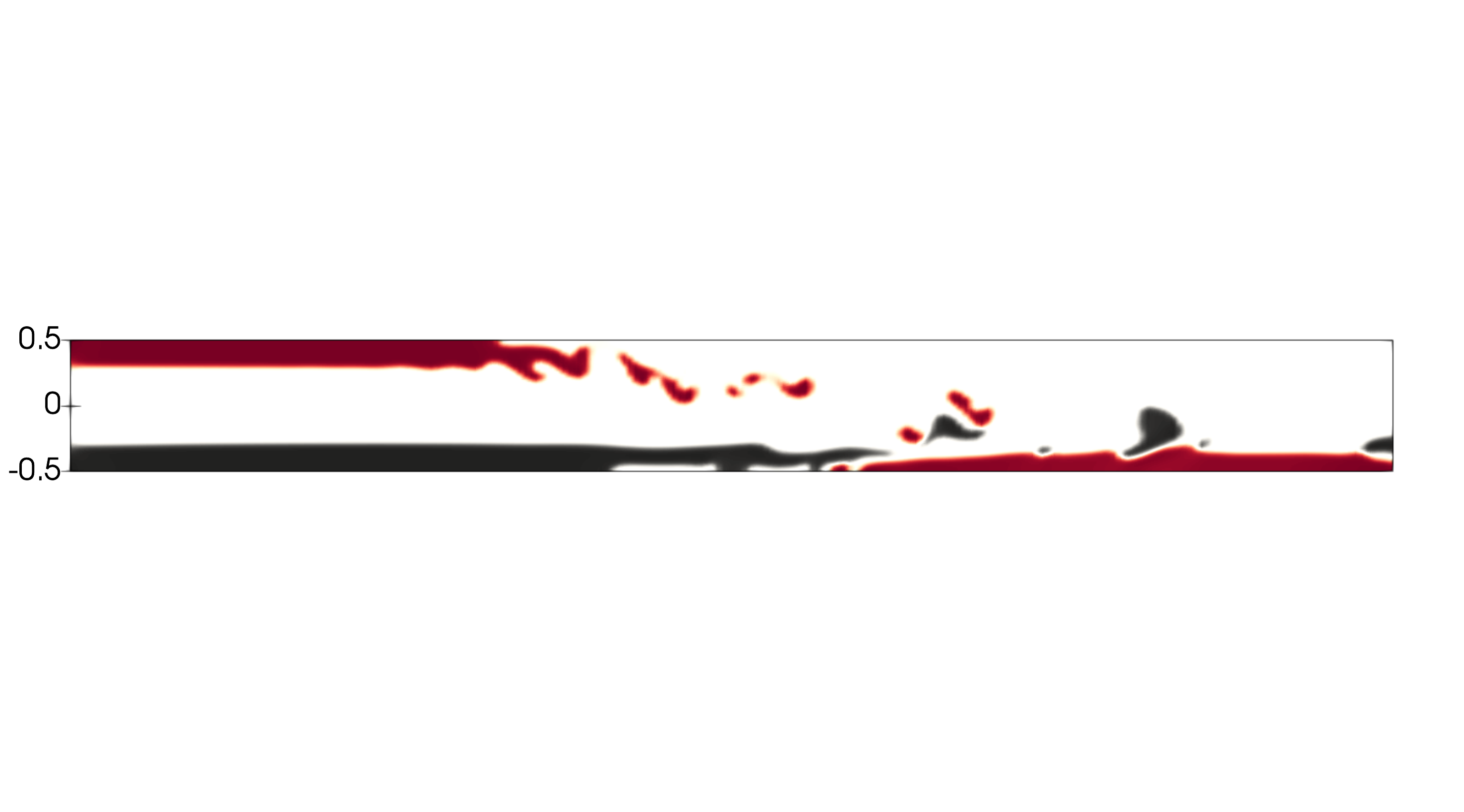}
	\includegraphics[clip,width=\textwidth,trim=0cm 11cm 1cm 11cm]{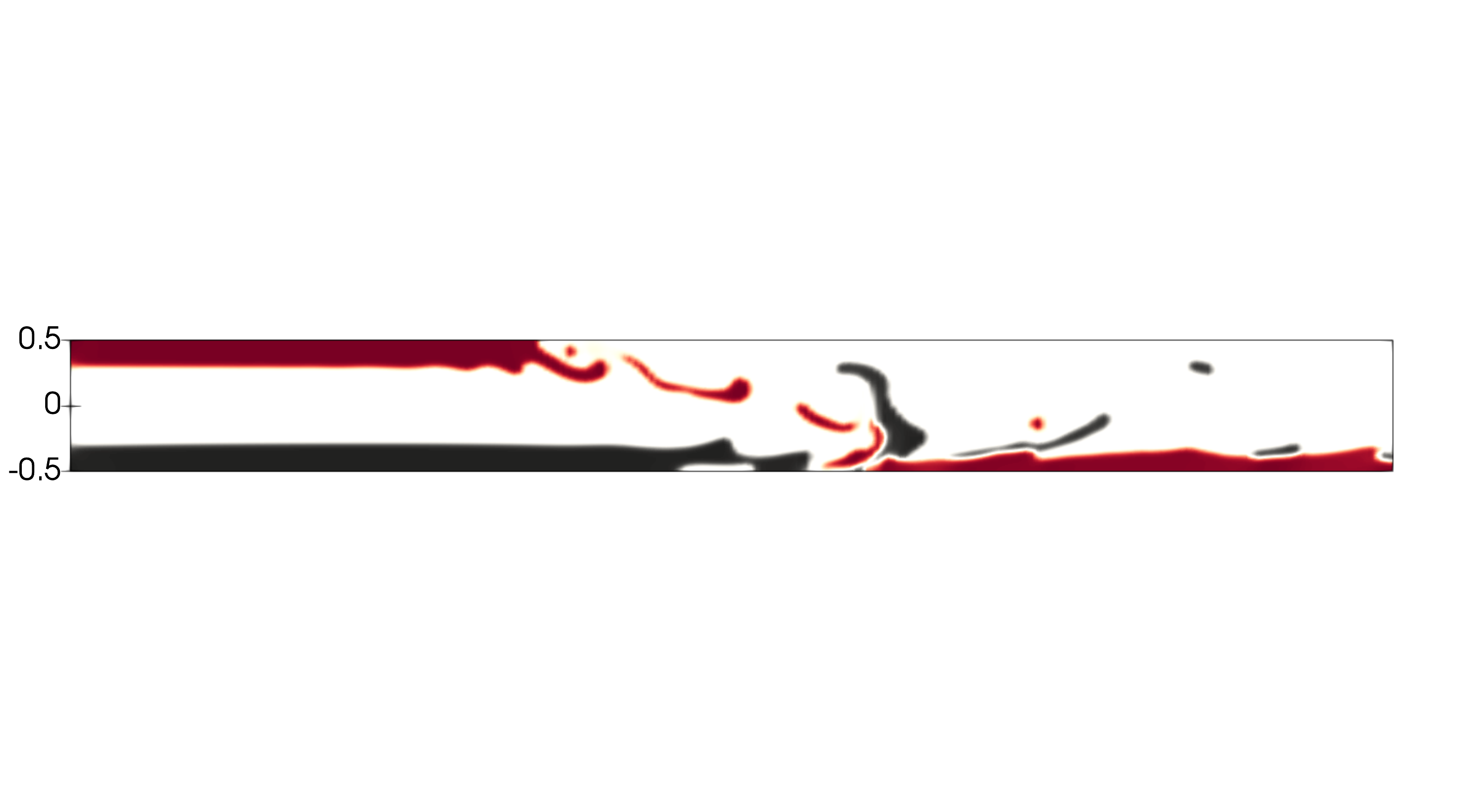}
	\caption{Three--phase solver: snapshots of the fluid evolution of the three--phase channel flow. We have represented the heavy fluid (Phase~2) in red, the light fluid (Phase~3) in black, and the medium density fluid (Phase~1) in white}			
	\label{fig-xpipe:num3ph:2D-channel}		
\end{figure}

In Fig.~\ref{fig-xpipe:num3ph:2D-channel}
we have represented 10 snapshots of the fluid configuration, corresponding to intervals of $0.3$~s. In red, we have represented the heavy fluid (Phase~2), and the light fluid (Phase~3) in black, whereas the medium density fluid (Phase~1) is white. We see that Phase~2 is prone to fall faster at the early stages because of the higher density ratio ($\rho_2/\rho_1=5$), while the rise of Phase~3 is more subtle (although it is more pronounced in the first snapshot), since the density ratio is lower ($\rho_3/\rho_1=0.8$). As a result, Phase~3 forms a film in the lower part of the channel, which is then broken by droplets of Phase~2 between $x=4$~m and $x=6$~m. These droplets are a result of the Rayleigh--Taylor instability between fluids 1 and 2 in the upper part, which periodically creates droplets of Phase~2 which fall to the lower part of the channel as a result of gravitational forces. As a result, we find that the heavy fluid accumulates in the lower part of the channel, while the light fluid is scattered into bubbles which leave the domain if they are big enough, whereas the smaller ones dissolve into Phase~1.

\subsection{Three--dimensional annular flow}\label{sec-xpipe:num3ph:annular}

In this section, we include a third phase to the annular flow simulation performed for the two--phase simplification in Sec.~\ref{sec-xpipe:2PH-pipe}. 
Therefore, to the initial configuration with $\rho_2=1$~kg/m$^3$ and $\rho_3=5$~kg/m$^3$, we add a lighter third phase with $\rho_1=0.5$~kg/m$^3$. The inflow configuration follows the description in Fig.~\ref{fig-xpipe:inflow-conf}. The physical parameters including superficial and slip velocities are given in Table~\ref{tab-xpipe:num3ph:3Dannular:param}.
\begin{table}[h]
		\centering
		\caption{Three--phase solver: list of the parameter values used for the three--phase annular flow}
		\label{tab-xpipe:num3ph:3Dannular:param}
		\begin{tabular}{llllllll}
			\hline
			$\rho_1$& $\rho_2$ & $\rho_3$ ($\text{kg}/\text{m}^3$)   & $\eta_1$ & $\eta_2$ & $\eta_3$  (Pa$\cdot$s) & $\varepsilon$ (m) \\ \hline
			0.5 & 1.0 & 5.0 & 1.0E-3 & 5.0E-3 & 1.0E-2 & $0.0424$   \\
			\hline \\ \hline
			&$M_0$ (m/s) & $c_0$ (m/s$^2$) & $\sigma_{12}$& $\sigma_{13}$ & $\sigma_{23}$ (N/m)\\ \hline
			&9.428E-5 & 1.0E3 &  2.5E-4 & 2.5E-4 & 2.5E-4\\ \hline\\
			\hline				
			&$V_{s,1}$&$V_{s,2}$&$V_{s,3}$ & $V_{s,12}$ & $V_{s,23}$ (m/s)  \\ \hline
			&1.0 & 3.9 & 0.06 & 0.0 & 10.0\\ \hline
		\end{tabular}
\end{table}
 We maintain the mesh used for the two--phase pipe simulations (see Sec.~\ref{sec-xpipe:2PH-pipe}), and we use  order $N=3$ polynomials. We use the IMEX time integrator with $S_0=8$ and $\Delta t=10^{-5}$~s.
 
\begin{figure}[h]
	\centering
	\subfigure[Three--dimensional view]{\includegraphics[width=\textwidth]{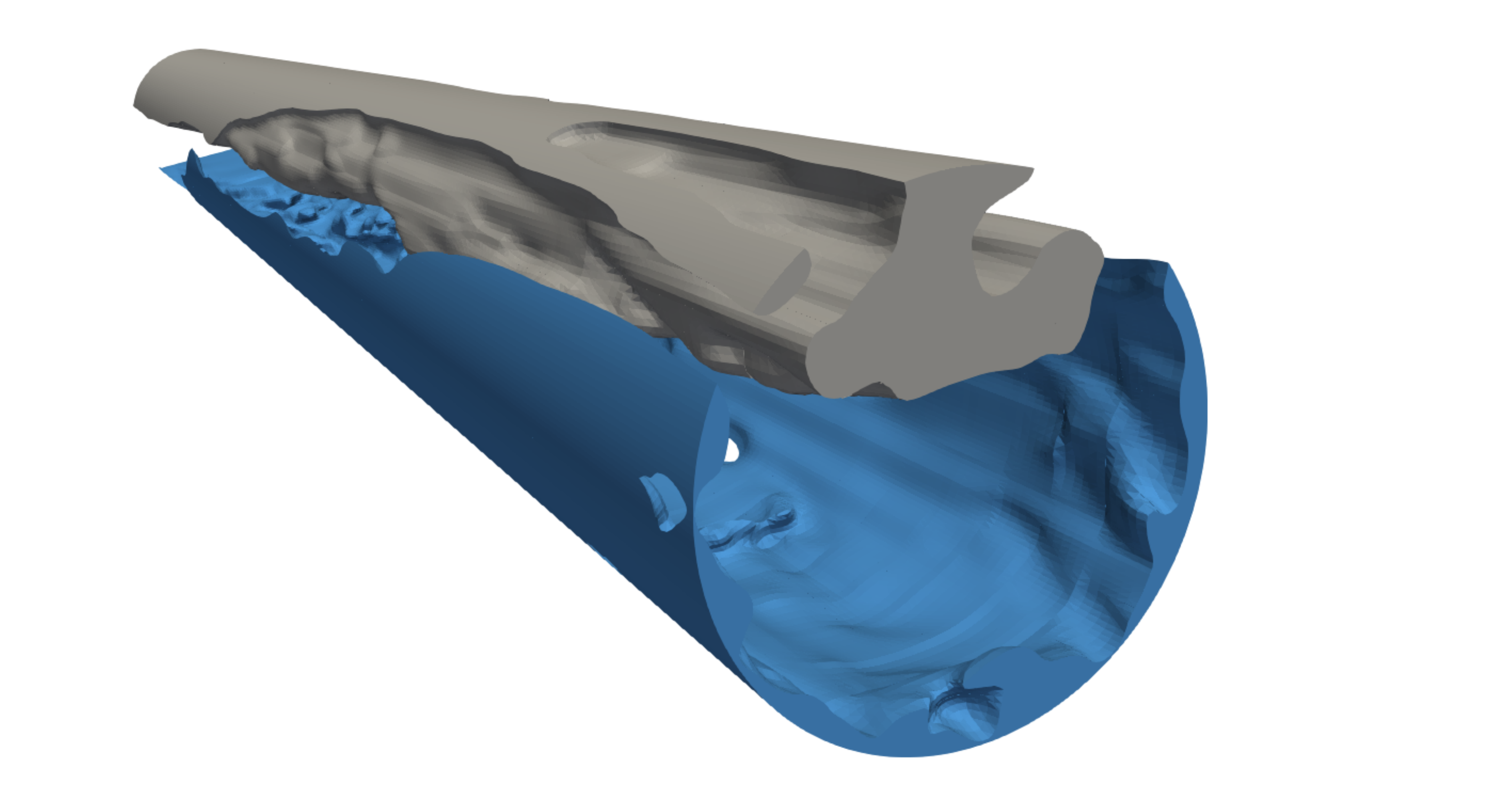}}
	\subfigure[Rear view]{\includegraphics[clip,height=0.4\textwidth,trim=14cm 1cm 14cm 1cm]{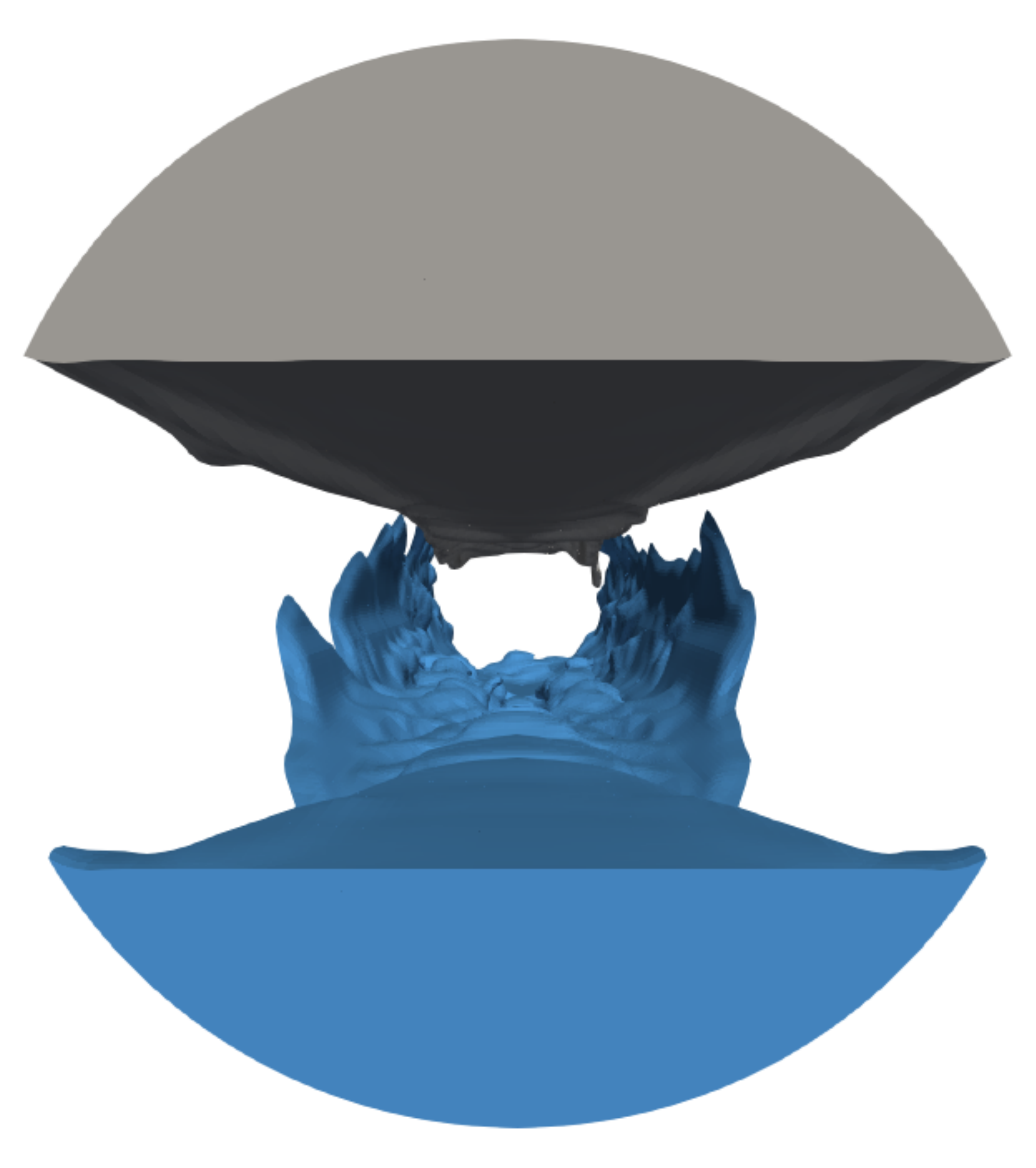}}	
	\subfigure[Front view]{\includegraphics[clip,height=0.4\textwidth,trim=13cm 1cm 13cm 1cm]{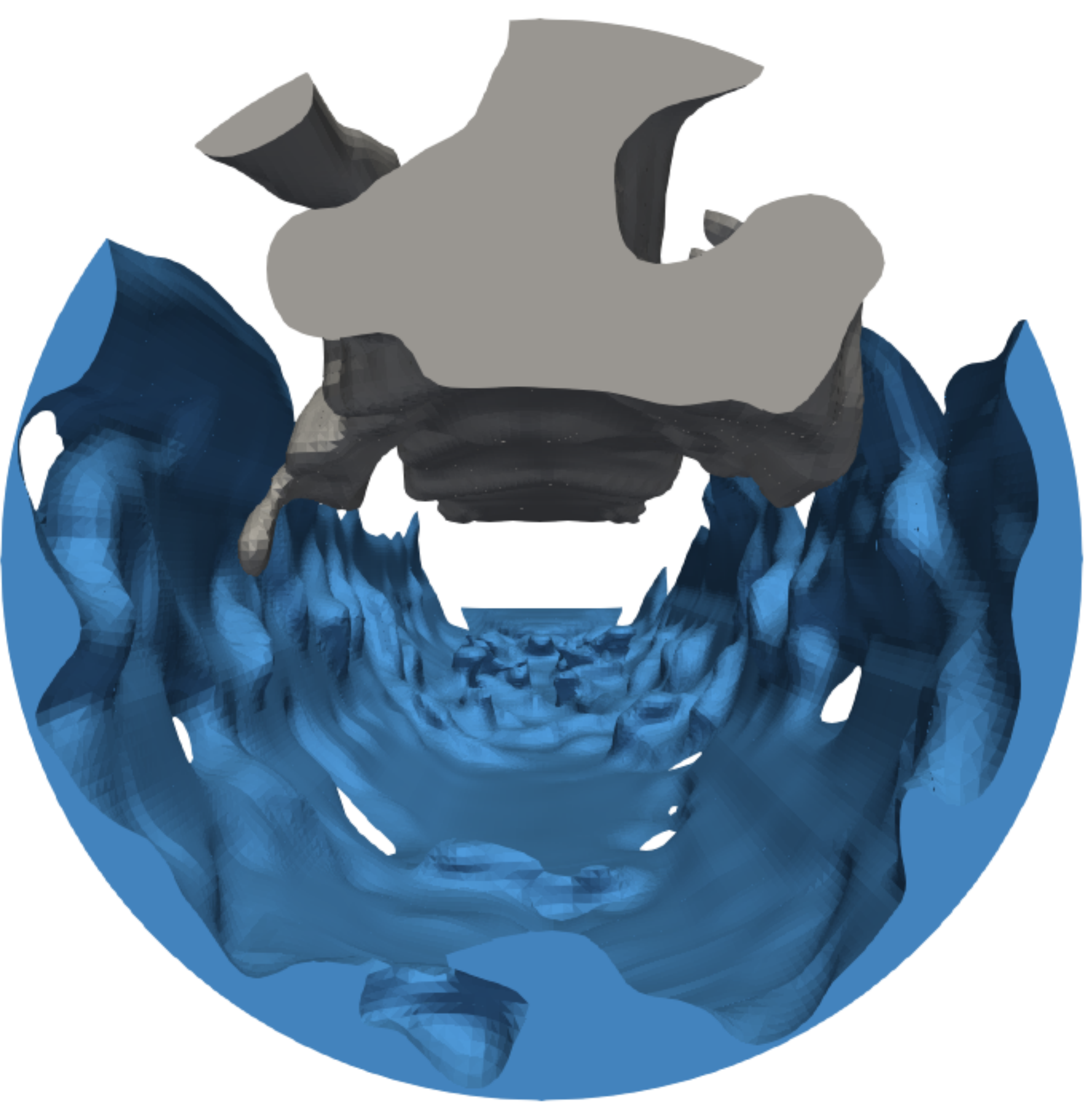}}		
	\subfigure[Side view]{\includegraphics[clip,width=\textwidth,trim=0cm 11cm 0cm 11cm]{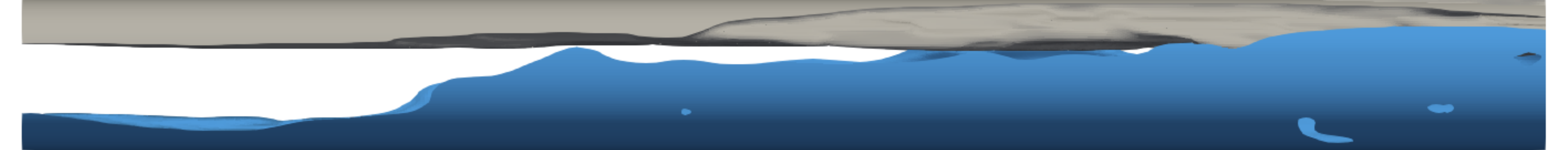}}			
	\subfigure[Top view]{\includegraphics[clip,width=\textwidth,trim=0cm 11cm 0cm 10cm]{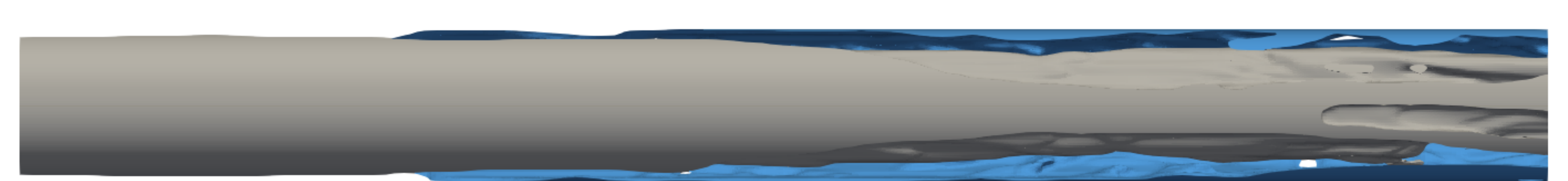}}	
	\caption{Three--phase solver: three--phase annular flow simulation: the fluid configuration is represented in $t=3$. The heavy fluid (Phase~3) is represented in blue, while the light fluid (Phase~1) is represented in gray}
	\label{fig-xpipe:num3ph:annular}
\end{figure}

We represent the flow configuration in $t=3$~s in Fig.~\ref{fig-xpipe:num3ph:annular}.
The light fluid (Phase~1) on top is represented in gray, and the heavy fluid (Phase~3) is represented in blue, both immersed in Phase~2 (not represented for clarity). We see that Phase~3 describes an annular flow regime very similar to that seen in the two--phase problem. The additional phase fills the upper space, which does not fall to the lower part because of its lower density.

\subsection{Three--dimensional T--shaped pipe intersection}\label{sec-xpipe:num3ph:Tee}

\begin{figure}[h]
  \centering
  \includegraphics[width=0.8\textwidth]{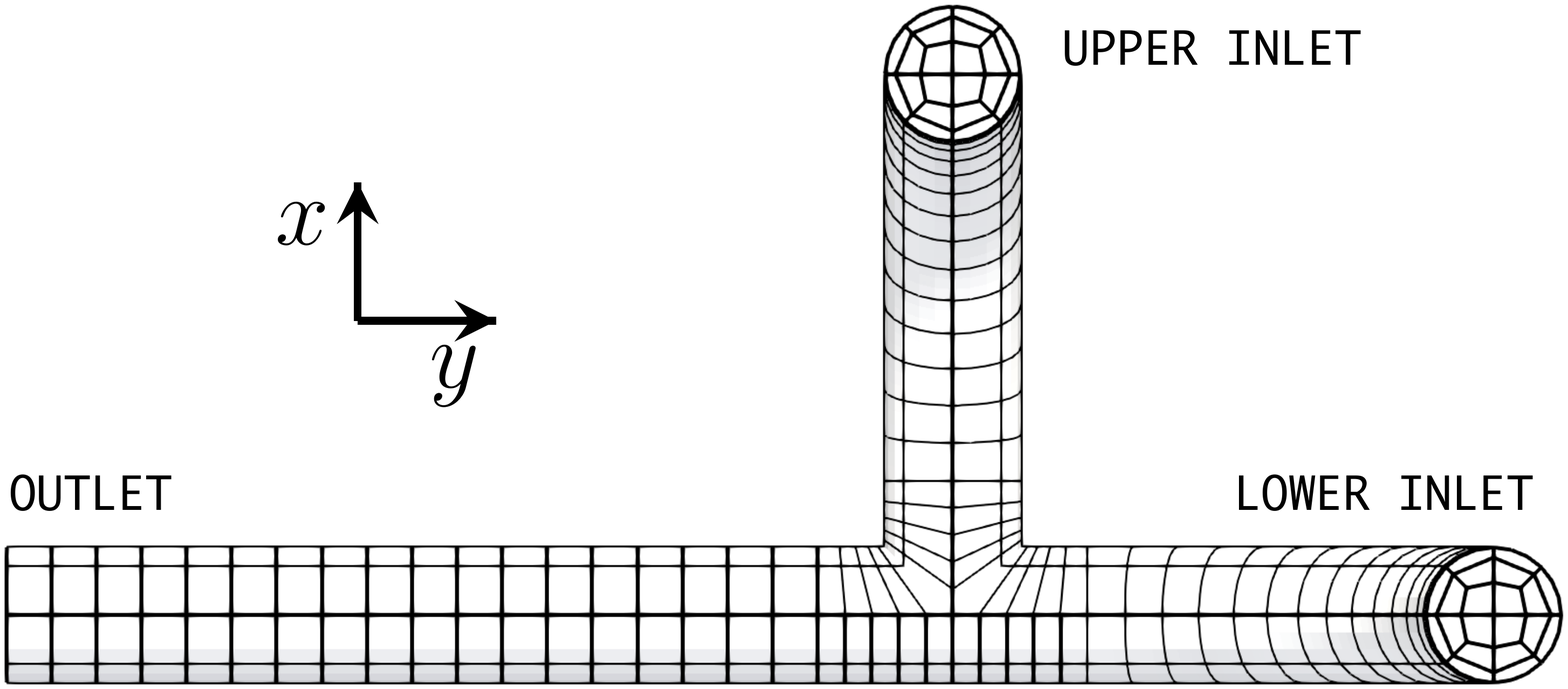}
  \caption{Three--phase solver: computational mesh for the T--shaped junction pipe domain, with 1700 elements}
  \label{fig-xpipe:num3ph:Tmesh}
\end{figure}

Finally, we solve a T--shaped pipe junction configuration, with two inlets and one outlet. The domain features a straight upper inlet whose length is $3$~m, which is then coupled to a $90^{\circ}$ bend whose radius is $3$~m. Additionally, the second inlet has a straight $5$~m section, and then another $90^{\circ}$/3~m bend. Finally, the outlet after the T--shaped junction is a straight pipe whose length is 6~m. The diameter of the pipe
is $D=1$~m. The computational mesh used, with 1700 elements, is represented in Fig.~\ref{fig-xpipe:num3ph:Tmesh}.\\

The physical parameters are given in Table~\ref{tab-xpipe:num3ph:3D-T:param}.
\begin{table}[h]
		\centering
		\caption{Three--phase solver: list of the parameter values used for the three--phase T--shaped pipe intersection simulation}
		\label{tab-xpipe:num3ph:3D-T:param}
		\begin{tabular}{llllllll}
			\hline
			$\rho_1$& $\rho_2$ & $\rho_3$ ($\text{kg}/\text{m}^3$)   & $\eta_1$ & $\eta_2$ & $\eta_3$  (Pa$\cdot$s) & $\varepsilon$ (m) \\ \hline
			5.0 & 1.0 & 0.2 & 1.0E-5 & 2.5E-5 & 5.0E-5 & $0.03$   \\
			\hline \\ \hline
			& $M_0$ & $c_0$ (m/s$^2$) & $\sigma_{12}$& $\sigma_{13}$ & $\sigma_{23}$ (N/m)\\ \hline
		    & 1.8856E-2 & 1.0E3 &  2.5E-4 & 2.5E-4 & 2.5E-4\\ \hline
		\end{tabular}
\end{table}
In the upper inlet, we only inject Phases~1~and~2 with superficial velocities $V_{s,1}=V_{s,2}=4$~m/s. In the lower inlet, we only inject Phase~3 with superficial velocity $V_{s,3}=2$~m/s. Additionally, the gravity acceleration is $\svec{g}=-1$~m/s$^2$ in x-direction. We use order $N=3$ polynomials and the IMEX scheme uses $S_0=8$ with a time--step size $\Delta t=5\cdot 10^{-5}$~s.
The initial condition is a steady--state with uniform pressure $p=0$, and with 
the pipe filled with Phase~3 ($c_1=c_2=0$). 

\begin{figure}[h]
  \centering
  \subfigure[$t=0$~s]{\includegraphics[clip,width=0.49\textwidth,trim=11cm 4cm 6cm 2cm]{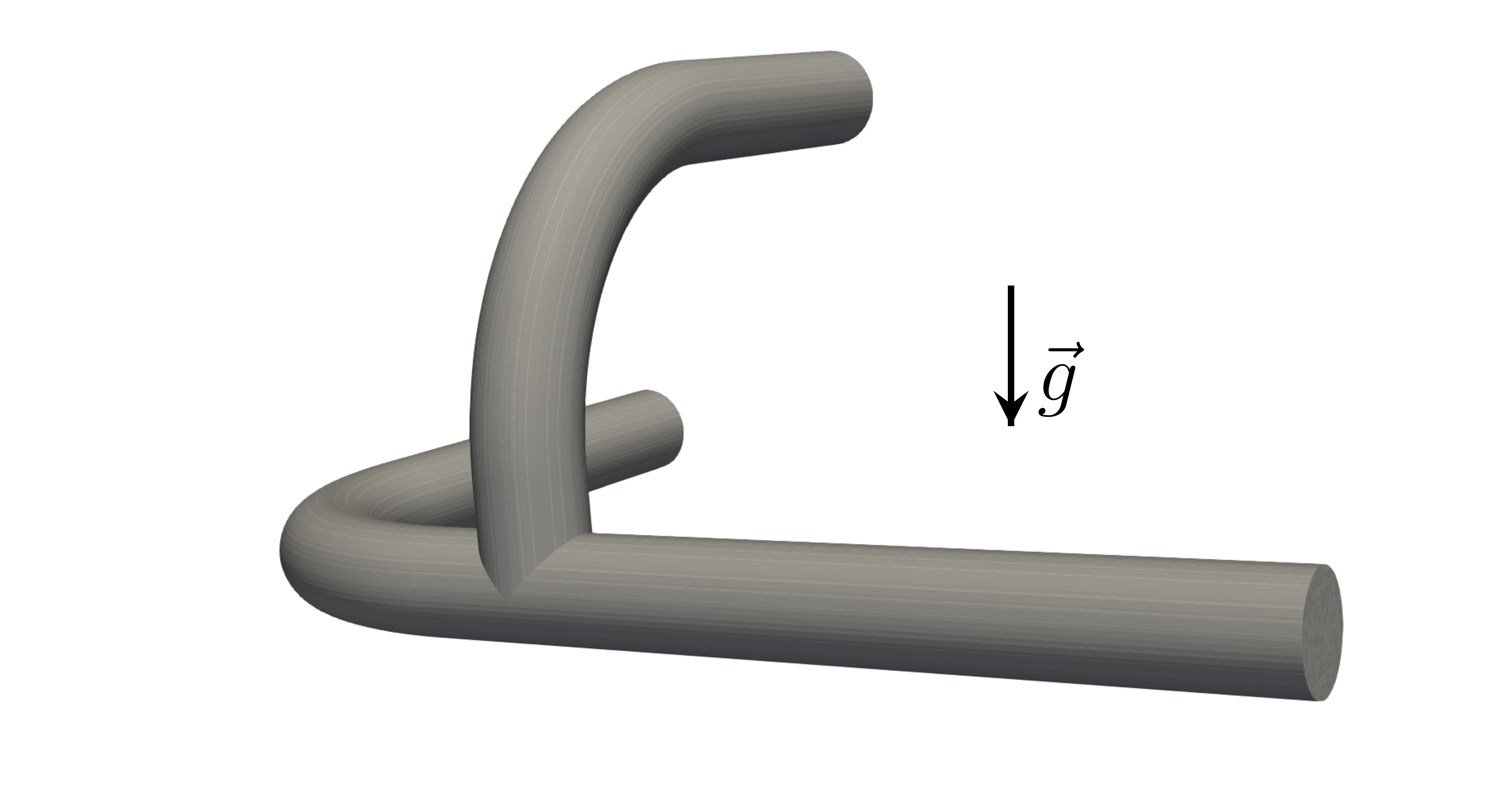}}
  \subfigure[$t=0.5$~s]{\includegraphics[clip,width=0.49\textwidth,trim=11cm 4cm 6cm 2cm]{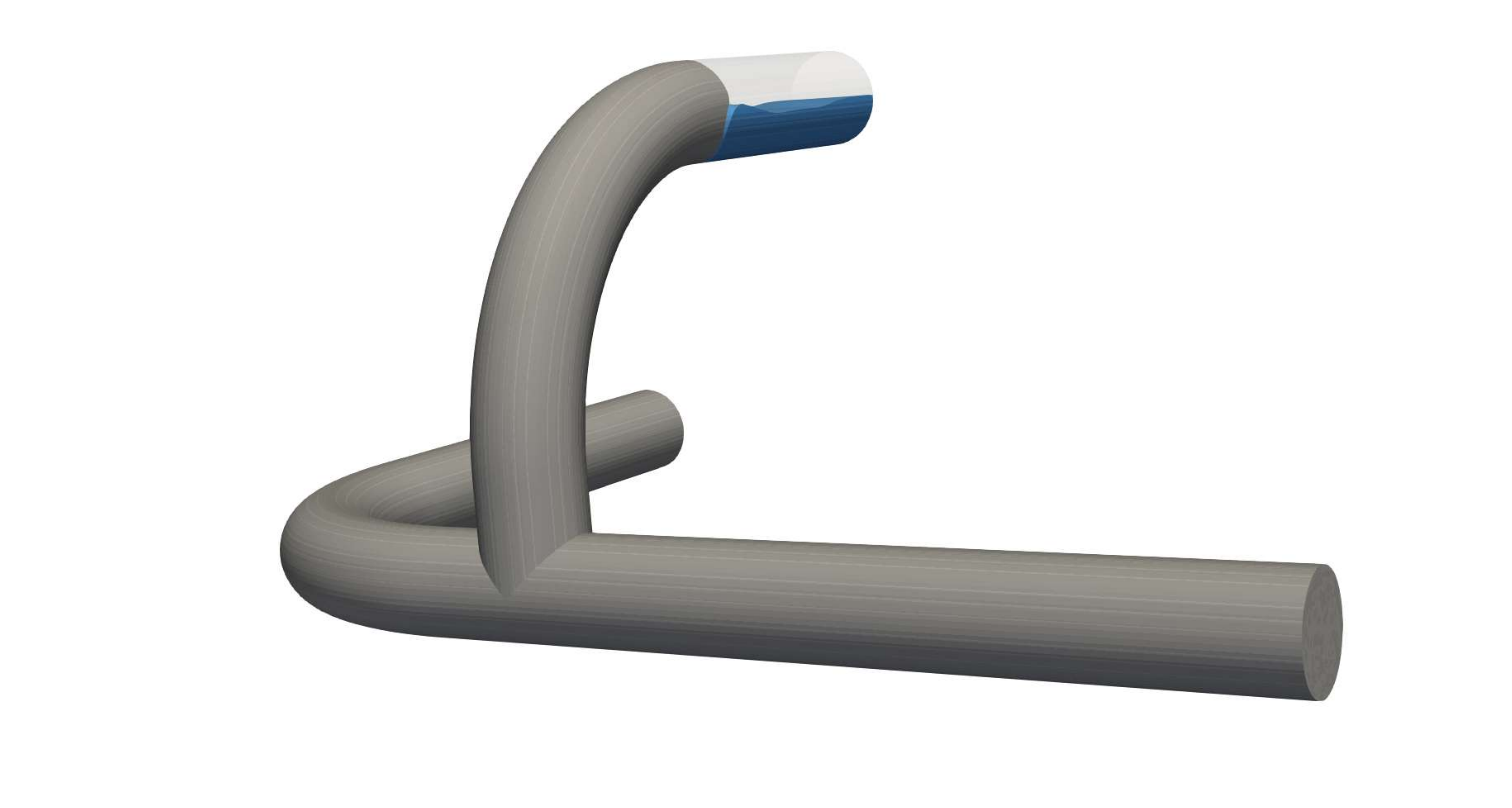}\label{fig-xpipe:num3ph:Tee:evolution-1}}
\end{figure}

\begin{figure}[h]
  \centering
  \subfigure[$t=1.0$~s]{\includegraphics[clip,width=0.49\textwidth,trim=11cm 4cm 6cm 2cm]{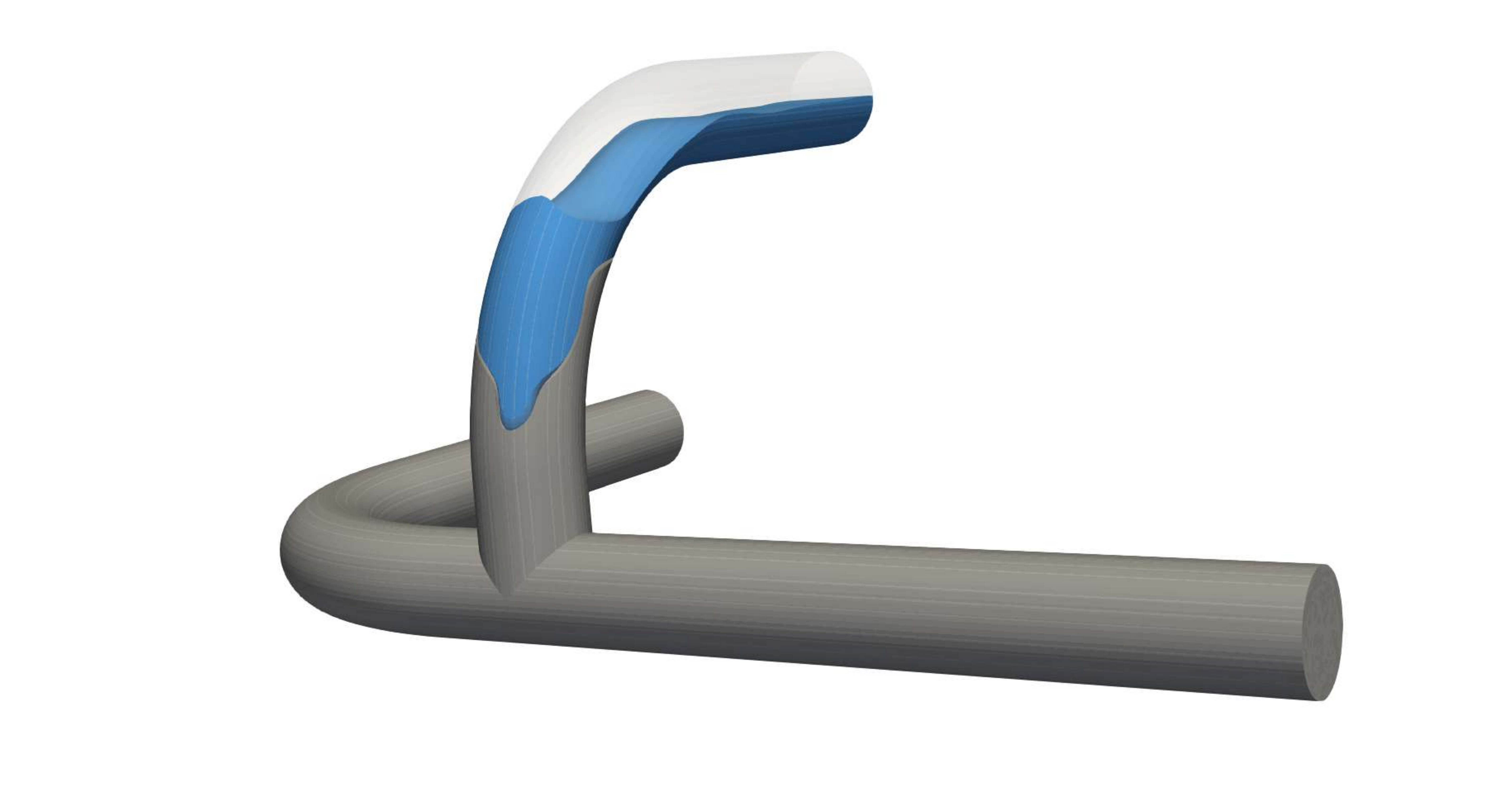}}
  \subfigure[$t=1.5$~s]{\includegraphics[clip,width=0.49\textwidth,trim=11cm 4cm 6cm 2cm]{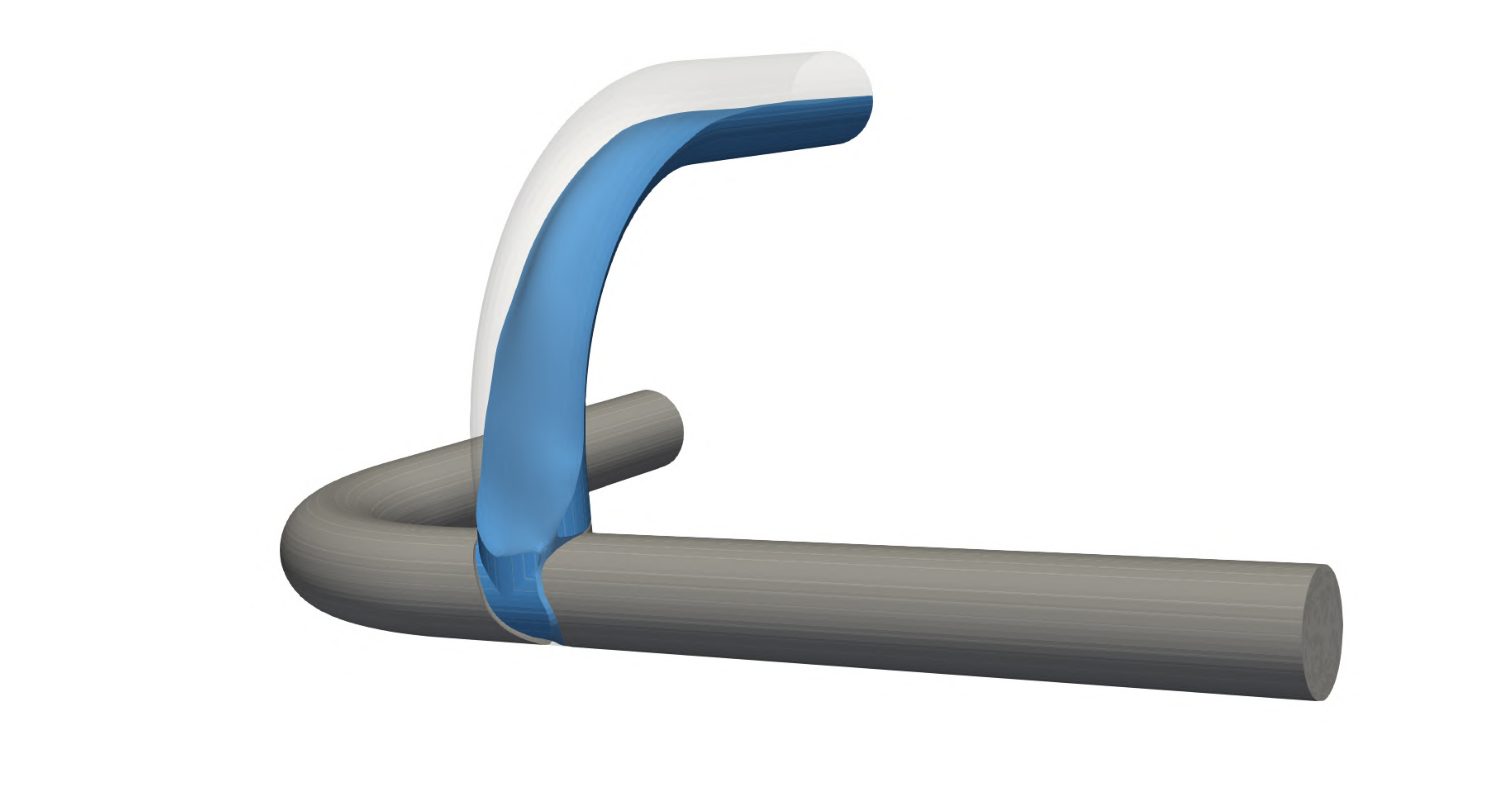}\label{fig-xpipe:num3ph:Tee:evolution-3}}
\end{figure}

\begin{figure}[h]
  \centering
  \subfigure[$t=4.0$~s]{\includegraphics[clip,width=0.49\textwidth,trim=11cm 4cm 6cm 2cm]{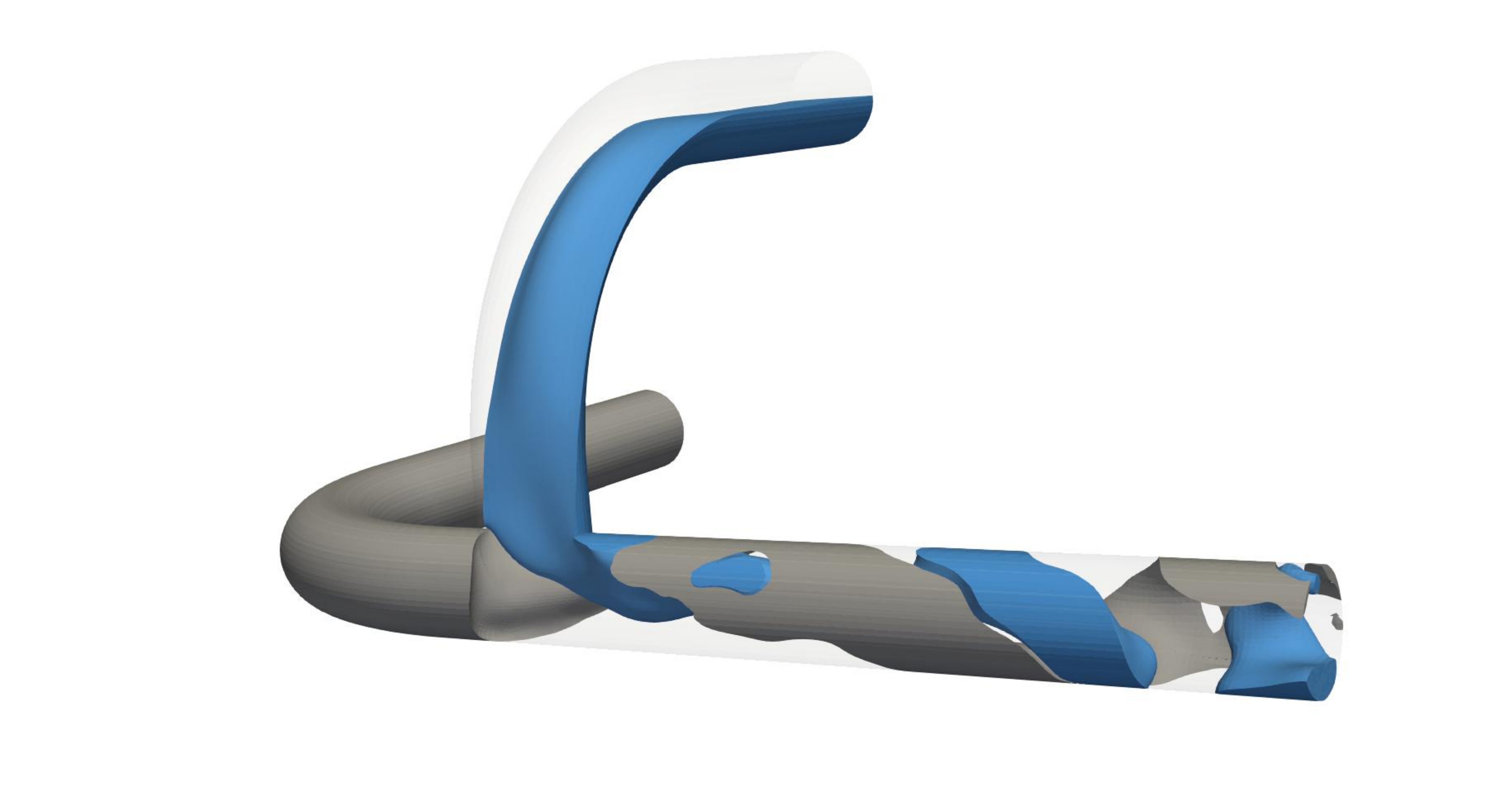}\label{fig-xpipe:num3ph:Tee:evolution-8}}
  \subfigure[$t=4.5$~s]{\includegraphics[clip,width=0.49\textwidth,trim=11cm 4cm 6cm 2cm]{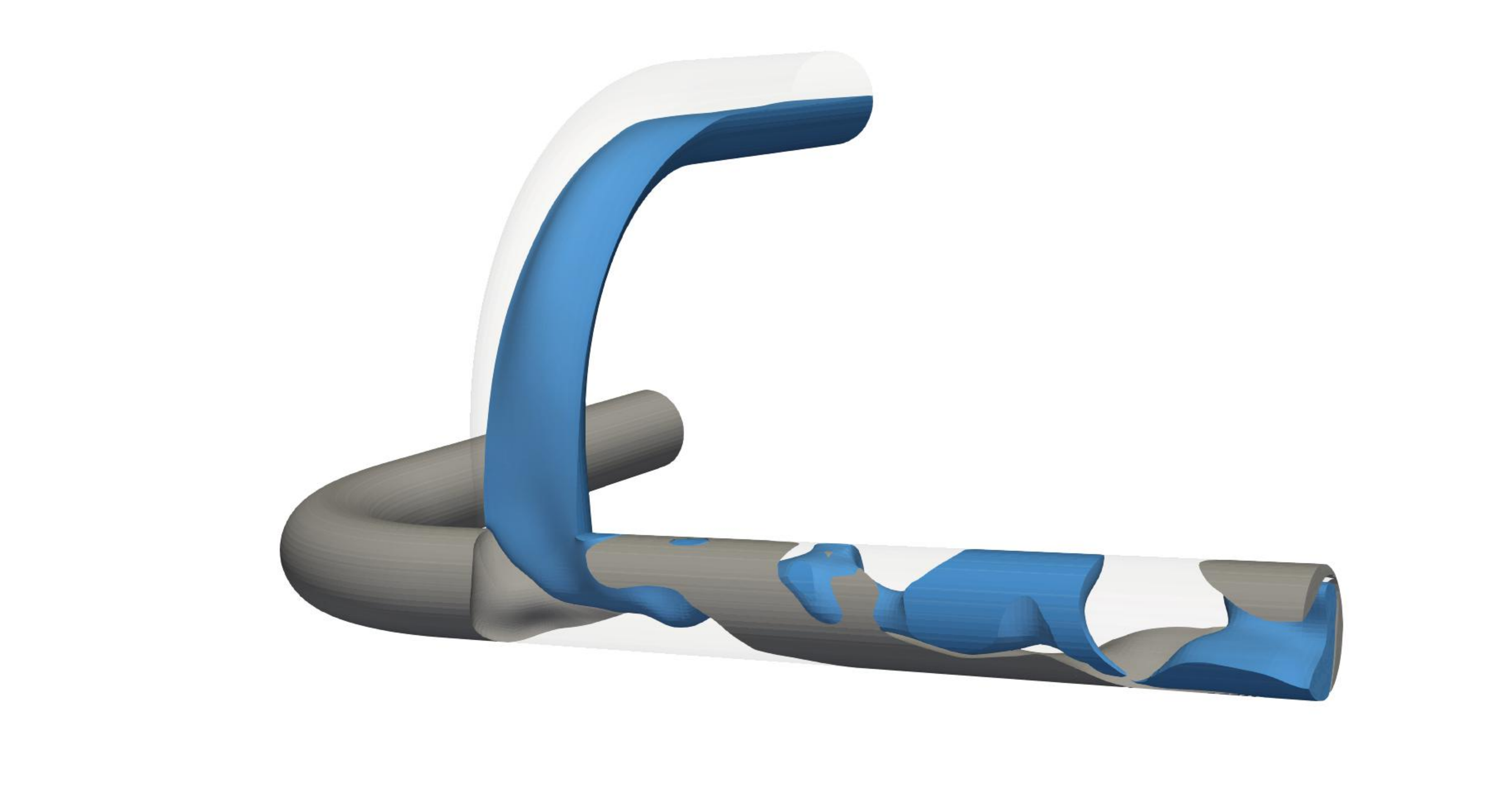}\label{fig-xpipe:num3ph:Tee:evolution-9}}
  \caption{Three--phase solver: evolution of phase two (blue) and phase three (gray) for the first $4.5$ seconds}
\label{fig-xpipe:num3ph:Tee:evolution}  
\end{figure}

We represent the evolution of the phases in Fig.~\ref{fig-xpipe:num3ph:Tee:evolution}, where we represent Phase~2 in blue, Phase~3 in gray and the space left is occupied by Phase~1. Initially the pipe is filled with Phase~3, which was chosen 
because it has the minimum density of the three--phase (therefore is easier for the other two phases to displace 
it). At the initial stages (see Fig.~\ref{fig-xpipe:num3ph:Tee:evolution-1}), we observe the advancing front at the upper inlet. We see that Phase~2 (blue) overtakes 
Phase~1 at the elbow, and then both phases arrive at the main pipe at $t\approx 1.5$ 
(see Fig.~\ref{fig-xpipe:num3ph:Tee:evolution-3}). Then, Phases~1~and~2 enter 
the principal pipe, and they restrict the flow of Phase~3 after the T--shaped pipe intersection. The lower density Phase~3 is then confined to the wall, and phases 1 
and 2 intermittently occupy the bulk of the pipe (see Fig.~\ref{fig-xpipe:num3ph:Tee:evolution-8}).

\begin{figure}[h]
  \centering
  \subfigure[$t=2.0$~s]{\includegraphics[clip,width=0.49\textwidth,trim=11cm 4cm 6cm 2cm]{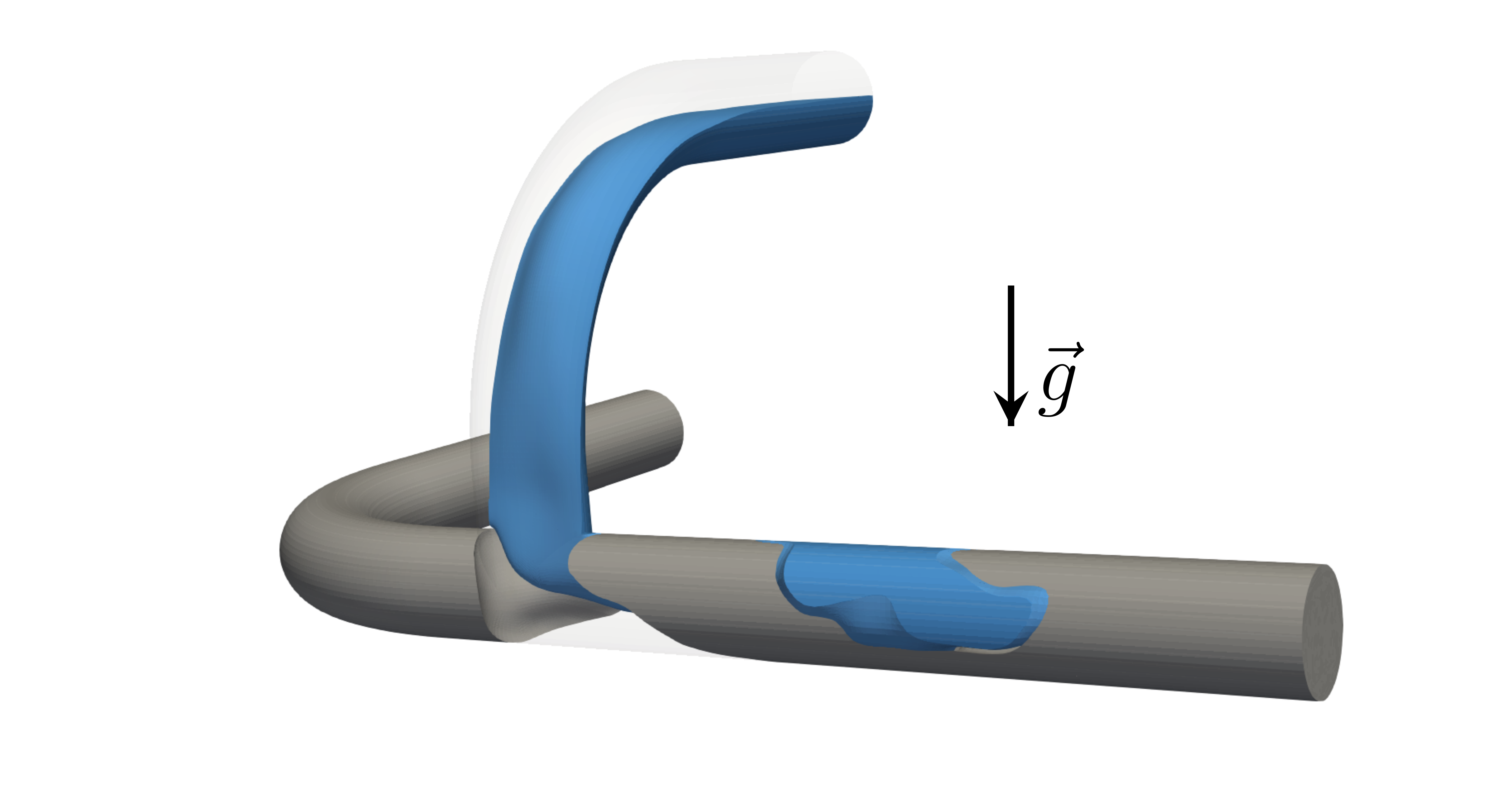}}
  \subfigure[$t=2.5$~s]{\includegraphics[clip,width=0.49\textwidth,trim=11cm 4cm 6cm 2cm]{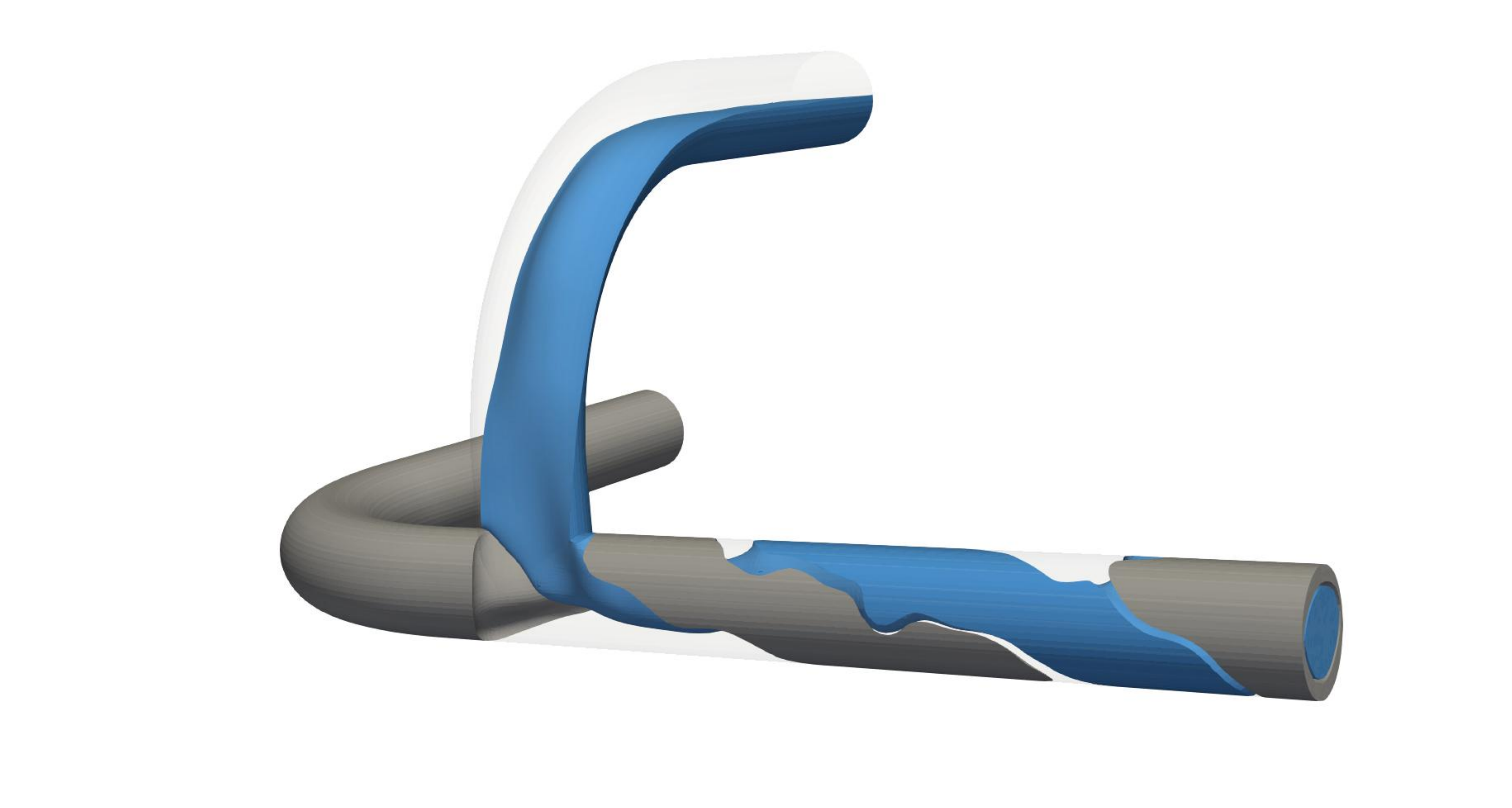}}
\end{figure}
\begin{figure}[h]
  \centering
  \subfigure[$t=3.0$~s]{\includegraphics[clip,width=0.49\textwidth,trim=11cm 4cm 6cm 2cm]{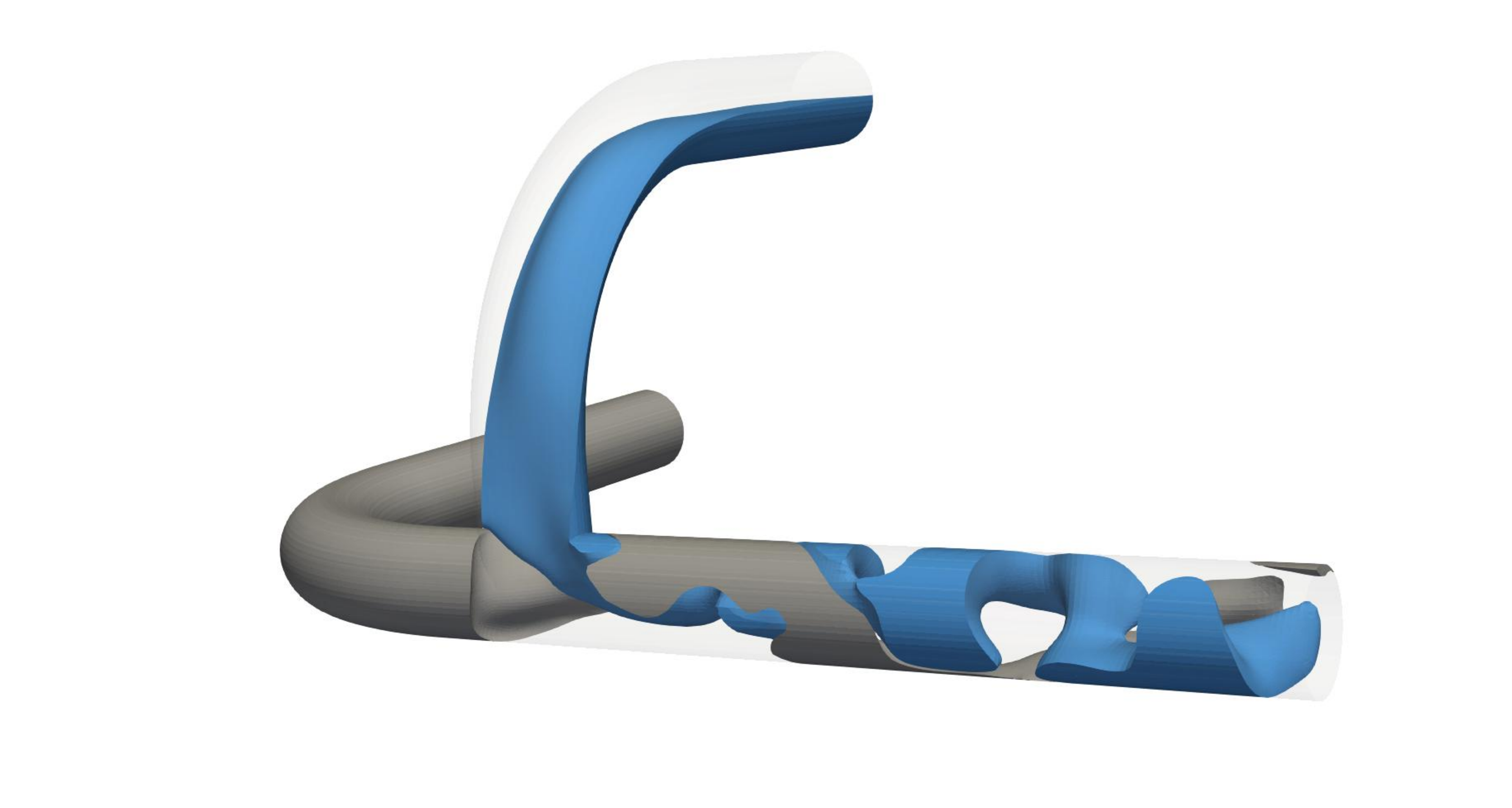}}
  \subfigure[$t=3.5$~s]{\includegraphics[clip,width=0.49\textwidth,trim=11cm 4cm 6cm 2cm]{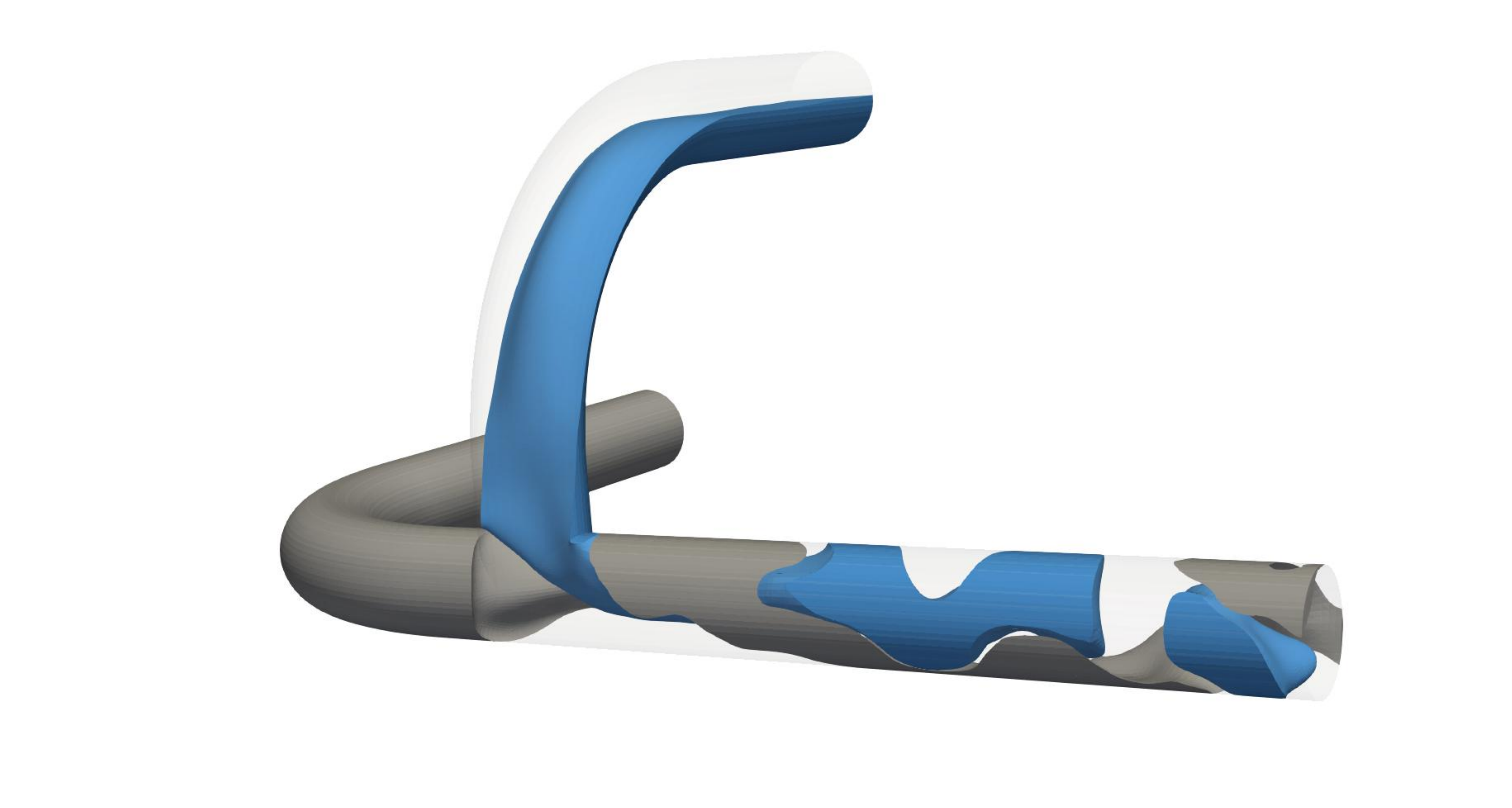}}
\end{figure}

\begin{figure}[h]
\centering
\subfigure[Phase~1]{\includegraphics[clip,width=0.9\textwidth,trim=8cm 2cm 1cm 
6cm]{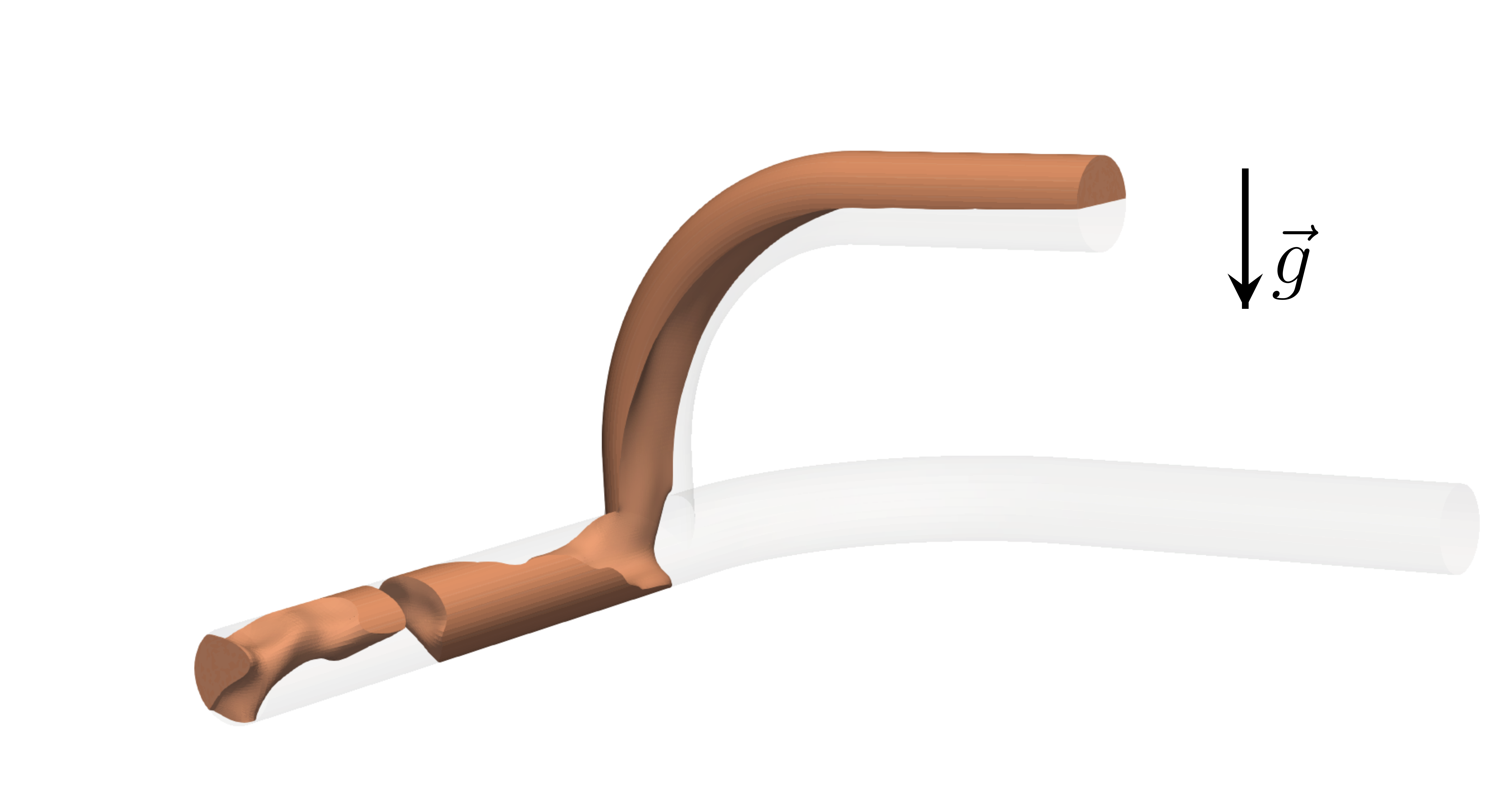}\label{fig-xpipe:num3ph:Tee:final-1}}
\subfigure[Phase~2]{\includegraphics[clip,width=0.9\textwidth,trim=8cm 2cm 1cm 
6cm]{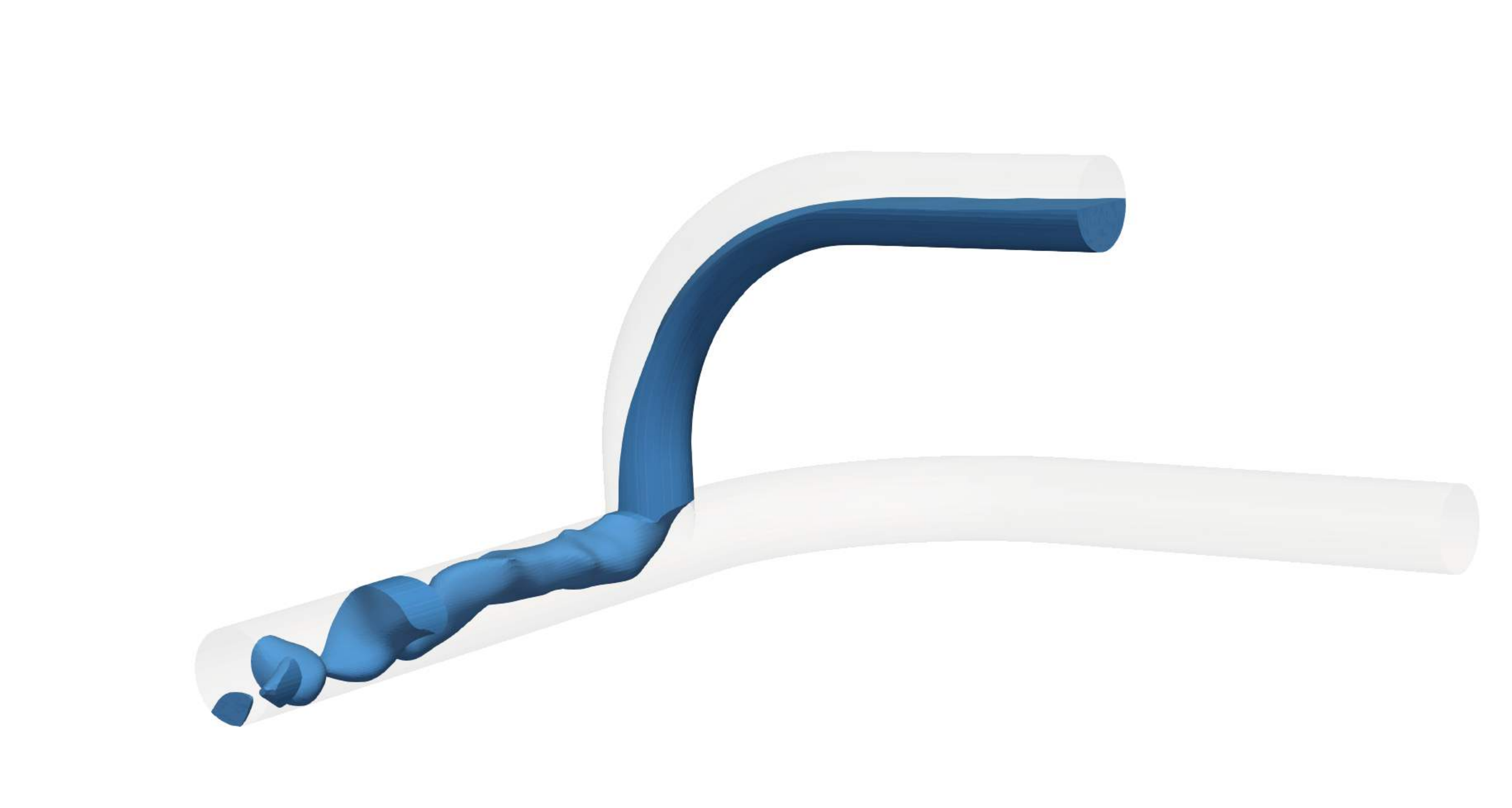}\label{fig-xpipe:num3ph:Tee:final-2}}
\subfigure[Phase~3]{\includegraphics[clip,width=0.9\textwidth,trim=8cm 2cm 1cm 
6cm]{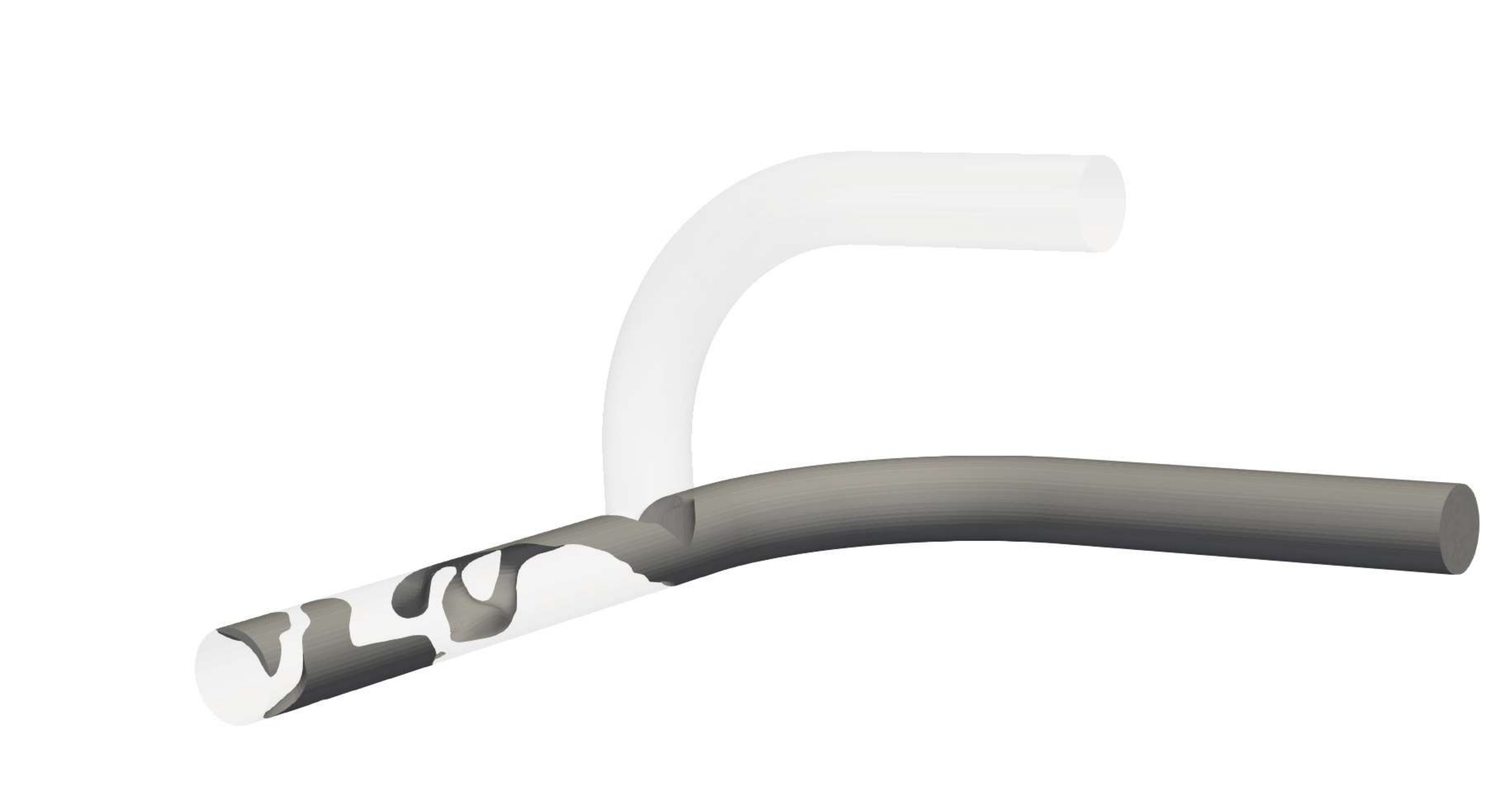}\label{fig-xpipe:num3ph:Tee:final-3}}
\caption{Three--phase solver: representation of the three--phase at the final time $t=7.5$}
\label{fig-xpipe:num3ph:Tee:final}
\end{figure}

In Fig.~\ref{fig-xpipe:num3ph:Tee:final} we represent the configuration of the 
three--phase at the final simulation time at ${t=7.5}$~s. Both Phases~1 (Fig.~\ref{fig-xpipe:num3ph:Tee:final-1}) and 2 (Fig.~\ref{fig-xpipe:num3ph:Tee:final-2}) 
represent the bulk of the pipe, whereas Phase~3 (Fig.~\ref{fig-xpipe:num3ph:Tee:final-3}) 
is forced to coat the pipe walls (similar to an annular flow regime). Due to the rupture of the flow of Phase~3 by the 
Phases~1~and~2, Phase~3 gets a counter--clockwise swirl motion around the pipe.

\begin{figure}[h]
  \centering
  \includegraphics[width=0.8\textwidth]{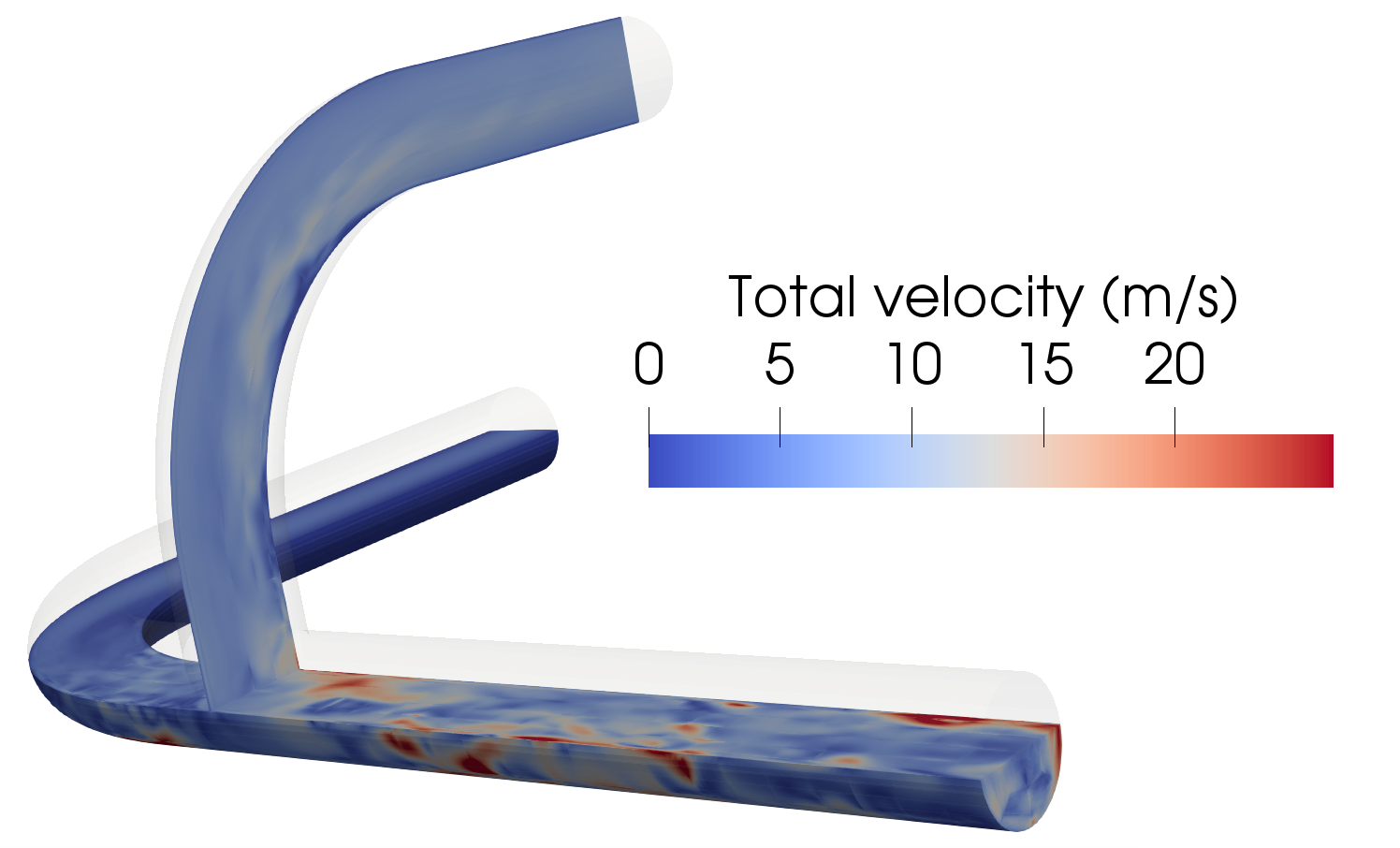}
  \caption{Three--phase solver: total velocity contour at the final time $t=7.5$~s}
  \label{fig-xpipe:num3ph:Tee:vel}
\end{figure}

Finally, we represent the velocity 
contours at the final time $t=7.5$~s in Fig.~\ref{fig-xpipe:num3ph:Tee:vel}. We can see the detachment due to the low viscosity of Phase~2 at the upper elbow, and also in the lower elbow for Phase~3. Then, after the T--shaped intersection, the flow 
becomes under--resolved with large velocity spots as a result of the interaction between the three phases.

Overall, we confirm the validity of the solver to compute simulations in complex three--dimensional geometries, as those found in the oil and gas industry.

\section{Conclusions}\label{sec-xpipe:Conclusions}


We present a three--phase incompressible Navier--Stokes/Cahn--Hilliard 
system, and its discontinuous Galerkin implementation. The model uses the 
three--phase Cahn--Hilliard model of \cite{boyer2006study}, and the incompressible Navier--Stokes with artificial compressibility of 
\cite{2019:Manzanero-iNS}. 

We construct a discontinuous Galerkin approximation of the equations, where we combine the scheme 
used for the three--phase Cahn--Hilliard model in \cite{2020:Manzanero-UR-CaF} and 
that used for the entropy--stable incompressible Navier--Stokes equations of \cite{2019:Manzanero-iNS}. 

We validate the solver in the two--phase simplification with a manufactured 
solution, and solving two--phase pipe regimes. Then, it is used to solve three--phase 
flows: a manufactured 
solution, a two--dimensional channel
and a three--dimensional T--shaped pipe intersection.
We find that the solver has not crashed in any of the simulations once the 
time--step size has been appropriately chosen. We highlight the ease in the 
configuration of the solver and the scheme for a user, as it only requires an appropriate 
choice of the physical parameters and conditions, plus the choice of the polynomial 
order of the simulation. The rest of the numerical parameters have been proven valid
in a vast range of flow conditions, and the boundary conditions are
automatically set--up by the algorithm that computes the inlet profile for given 
superficial/slip velocities.
The enhancement of the robustness added by the split--form scheme, although
formally does not satisfy a discrete entropy inequality, has been addressed with numerical experiments. Although the solutions of the flows presented are under--resolved, a better resolution can be achieved by increasing the polynomial order, which avoids re--meshing the geometry.


	\section*{Acknowledgement}
	
	The authors acknowledge the computer resources and technical assistance provided by the Centro de Supercomputaci\'on y Visualizaci\'on de Madrid (CeSViMa). The authors acknowledge Repsol Technology Lab and Universidad Polit\'ecnica de Madrid for their support and permission to publish this work. Gonzalo Rubio and Eusebio Valero acknowledge the funding received
by the project SIMOPAIR (Project No. REF: RTI2018-097075-B-I00) from the Ministry of Innovation of Spain.
	Authors also thank Gabriel Rucabado from Repsol Technology Lab for his assistance during the execution of this work.

\clearpage

	\bibliographystyle{elsarticle-num-names}
	\bibliography{mybibfile.bib}
	
\end{document}